\documentclass[11pt]{article}
\usepackage{amscd,amsfonts,amssymb,amscd,amsmath,latexsym}
\usepackage[hypertex]{hyperref}
\usepackage[all]{xy}

\makeatletter
\renewcommand{\subsubsection}{\@startsection
{subsubsection}
{1}
{0mm}
{0mm}
{0mm}
{\normalfont\normalsize\itshape}}
\makeatother

\begin{titlepage}

\title{Scattering theory for geometrically finite groups}

\author{Ulrich Bunke\thanks{NFW I--Mathematik, Universit\"at Regensburg, D-93040 Regensburg, Germany, E-mail: ulrich.bunke@mathematik.uni-regensburg.de} and
Martin
Olbrich\thanks{Universit\'e du 
Luxembourg,
UR en Math\'ematiques,
162 a, avenue de la Fa\"{\i}encerie,
L-1511 Luxembourg,
Grand-Duchy of Luxembourg, E-mail: martin.olbrich@uni.lu  }
}

\end{titlepage}

\newcommand{\proof}{{\it Proof.$\:\:\:\:$}}
 \newcommand{\dist}{{\rm dist}}

\newcommand{\paaa}{{\mathfrak{p}}}
\newcommand{\taaa}{{\mathfrak{t}}}
\newcommand{\haaa}{{\mathfrak{h}}}
\newcommand{\R}{{\mathbb{R}}}
\newcommand{\Hh}{{\mathbb{H}}}

\newcommand{\Z}{{\mathbb{Z}}}
\newcommand{\bF}{{\mathbb{F}}}
\newcommand{\C}{{\mathbb {C}}}

\newcommand{\gaaa}{{\mathfrak{g}}}
\newcommand{\zaaa}{\mathfrak{z}}

\newcommand{\ausg}{{\rm end}}

\newcommand{\maaa}{{\mathfrak{m}}}
\newcommand{\aaaa}{{\mathfrak{a}}}
\newcommand{\naaa}{{\mathfrak{n}}}

\newcommand{\Pol}{{\rm Pol}}
\newcommand{\Tr}{{\rm Tr}}

\newcommand{\cZ}{{\mathcal{Z}}}

\newcommand{\cE}{{\mathcal{E}}}

\newcommand{\cV}{{\cal V}}
\newcommand{\cY}{{\cal Y}}
\newcommand{\cW}{{\mathcal{ W}}}

\newcommand{\cA}{{\cal A}}

\newcommand{\cU}{{\mathcal{U}}}
\newcommand{\colim}{\mathrm{colim}}
\renewcommand{\lim}{\mathrm{lim}}
\newcommand{\Hom}{{\mbox{\rm Hom}}}
\newcommand{\vol}{{\rm vol}}
\newcommand{\cO}{{\mathcal{O}}}
\newcommand{\End}{{\mathrm{End}}}
\newcommand{\Ext}{{\mathrm{ Ext}}}

\newcommand{\im}{{\mbox{\rm im}}}

\newcommand{\spann}{{\mathrm{span}}}

\newcommand{\cF}{{\mathcal{F}}}

\newcommand{\Ree}{{\rm Re }}

\newcommand{\Imm}{{\rm Im}}

\newcommand{\ee}{{\rm e}}

\newcommand{\tr}{{\mbox{\rm tr}}}

\newcommand{\ad}{{\mbox{\rm ad}}}

\newcommand{\Gr}{{\mathrm{Gr}}}

\newcommand{\id}{{\mathrm {id}}}

\newcommand{\nat}{{\mathbb {N}}}
\newcommand{\supp}{{\mathrm{supp}}}

\newcommand{\aca}{{\aaaa_\C^\ast}}

\newcommand{\cB}{{\mathcal{B}}}

\newcommand{\cR}{{R}}
\def\hB{\hspace*{\fill}$\Box$ \newline\noindent}

\newcommand{\cS}{{S}}

\def\hB{\hspace*{\fill}$\Box$ \\[0.5cm]\noindent}
\newcommand{\cL}{{\cal L}}

\newcommand{\cP}{{\mathcal{P}}}

 \newcommand{\cX}{{\mathcal{X}}}

\newtheorem{prop}{Proposition}[section]
\newtheorem{lem}[prop]{Lemma}
\newtheorem{ddd}[prop]{Definition}
\newtheorem{theorem}[prop]{Theorem}
\newtheorem{kor}[prop]{Corollary}
\newtheorem{ass}[prop]{Assumption}

\newcommand{\vp}{{\varphi}}
\newcommand{\RH}{{H\mathbb{R}}}
\newcommand{\CH}{{H\mathbb{C}}}
\newcommand{\HH}{{H\mathbb{H}}}
\newcommand{\OH}{{H\mathbb{O}}}

\begin{document}
\maketitle
\nocite{*}
\tableofcontents
     \section{Introduction and statement of the results} 

\subsection{Introduction}

\subsubsection{}

Scattering theory is a collection of methods and results for studying the continuous spectrum of operators.
Classical fields of application of scattering theory are the propagation of waves in mechanics, electro-dynamics and quantum mechanics. Often one considers a model that is close to some ideal case in which the wave equation can be solved explicitely. Scattering theory provides the tools for
comparison.

\subsubsection{}

 In many cases it is possible to separate the time and space variables. In this case one can apply methods from stationary scattering theory. This branch of scattering theory investigates the generalized eigenfunctions contributing to the continuous spectrum. In good situations these eigenfunctions come in families which extend meromorphically. Often they are studied via a meromorphic continuation of the distribution kernel
of the resolvent of the spacial part of the operator.

\subsubsection{}

The problem of stationary scattering theory that is considered in the present paper has mainly an internal mathematical motivation.  The model situation is a Laplace-type operator on a globally symmetric space 
of negative curvature. Due to the presence of a large symmetry group it is possible to get an essentially complete description of the generalized eigenfunctions.

The real problem is to understand the generalized eigenfunctions on  associated locally symmetric spaces.
The special case of quotients of the globally symmetric space by arithmetic subgroup groups has many applications in number theory and arithmetic. While arithmetic quotients have finite volume the main emphasis in the present work is on spaces of infinite volume. The question is, how far the interplay between global and local symmetries of the problem can be used in order to get a complete picture.

\subsubsection{}

There is lot of literature on the stationary scattering theory for the function Laplacian on spaces which become close to the global symmetric space asymptotically at infinity. The major part is devoted to the asymptotically hyperbolic case. The main results are meromorphic continuation of the resolvent kernel, 
the scattering matrix and the parametrization of the relevant generalized eigenfunctions. 
Most of this work is based on a fine study of the resolvent kernel. We refer to \cite{MR2153454} for one of the latest works and the discussion of the literature therein.

\subsubsection{}

In the more rigid case of a locally symmetric space one can obtain better results.
In this case one can approach stationary scattering theory by pushing the analysis to the boundary.
Using this method we obtained in  \cite{MR1749869} the spectral decomposition of all locally invariant differential operators on vector bundles in the  convex-cocompact case.
We refer to this paper for a review of the literature about the convex-cocompact case.

\subsubsection{}

The goal of the present paper is to generalize this method to geometrically finite spaces.
These are locally symmetric spaces of negative curvature which at infinity look like the symmetric space or a cusp. From the point of view of analysis the part at infinity which can be compared with the globally symmetric space is easy with the results of  \cite{MR1749869} at hand.
The complications are due to the presence of cusps, in particular of those of non-maximal rank.

\subsubsection{}
Let $G$ be a real simple linear connected Lie group of real rank one.
It covers the connected component 
 of the group of isometries of the associated symmetric space $X$.
We consider a geometrically finite group $\Gamma\subset G$ which determines the locally symmetric space $\Gamma\backslash X$. For technical reasons we exclude the exceptional symmetric space $X=\OH^2$
from our considerations.

The classical problems of meromorphic continuation of the Eisenstein series, the scattering matrix,
and their functional equations  was previously
adressed by Guillop\'e \cite{guillope92} in the special case
that $X$ is the hyperbolic plane $\RH^2$ and by
Froese/Hislop/Perry \cite{MR1111571} for $X=\RH^3$
and spherical Eisenstein series. 
Cusps of non-maximal rank in the three-dimensional hyperbolic case have been investigated first in 
\cite{MR1128217} via the resolvent. For the meromorphic continuation of the resolvent on 
$\Gamma\backslash\RH^n$ see \cite{MR1710792},\cite{MR2209763}. All these papers are restricted to the case
of so-called rational cusps.

\subsubsection{}

The main result of the present paper is, in the geometrically finite case, the meromorphic continuation of the
Eisenstein series, i.e. of the families of generalized eigensections for all locally invariant differential operators on bundles, see Corollary \ref{main} below. Note that this result also includes the case of irrational cusps.

At the moment we are not able to obtain results, which are as complete as in the convex cocompact case.
In particular,
we are still quite far from showing a spectral decomposition.

\subsubsection{}

The class of geometrically finite discrete subgroups $\Gamma\subset G$
subsumes cocompact, convex-cocompact
groups, and  subgroups of finite covolume as well as various kinds
of combinations of these cases with cusps of non-maximal rank.
The paper can be considered as a continuation
of  \cite{MR1749869}, where the special case of a convex-cocompact
group $\Gamma$ is considered. We will frequently refer to this paper
for notations and conventions as well as for technical results.

\subsubsection{}

In the present paper the analysis takes place on the boundary $\partial X$ of
the symmetric space.
The group $\Gamma$ acts on vector bundles over this boundary.
As in \cite{MR1749869} the main players  of the geometric
scattering theory are the push-down $\pi^\Gamma_*$ (average of sections over $\Gamma$),
the extension $ext^\Gamma$ (the adjoint ot the push-down), the restriction $res^\Gamma$, and the
scattering matrix $S^\Gamma$.
Our task is to construct these as families of maps between suitable function/distribution spaces
depending meromorphically on the spectral parameter.
Due to the presence of cusps of non-maximal rank the detailed description
of these spaces turns out to be quite difficult.


\subsection{Notion of geometrical finiteness}\label{geomf}

\subsubsection{}
In the present paper the symmetric space $X$ belongs to one of the series
of hyperbolic spaces $\RH^n$, $\CH^n$, and $\HH^n$
(real, complex, and quaternionic) of real dimension $n$, $2n$, $4n$.
We exclude the exceptional symmetric space of rank one, the Cayley hyperbolic plane $\OH^2$,
for technical reasons since we are going to employ the fact that
$X$ belongs to a series in several places.

\subsubsection{}
By $G$ we denote a real simple linear connected Lie group of real rank one
covering the connected component 
 of the group of isometries of $X$.
By $\partial X$ we denote the geodesic boundary
of $X$.
The union $\bar X:= X\cup \partial X$ has the structure of
a compact $G$-manifold with boundary. While $X$ parametrizes
the set of maximal compact subgroups the boundary $\partial X$
parametrizes the set of parabolic subgroups of $G$. 
If $P\subset G$ is a parabolic subgroup, then we denote by
$\infty_P\in \partial X$ the unique fixed point of $P$.
Let $\Omega_P:=\partial X\setminus \infty_P$.

\subsubsection{}
Let $P\subset G$ be parabolic, and $N\subset P$ be its nil-radical.
Then we can write $P$ as an extension
\begin{equation}
\label{eepp}
0\rightarrow N\rightarrow P\stackrel{l}{\rightarrow} L\rightarrow 0\ .\end{equation}
Here $L$ is a reductive Lie group which is canonically
isomorphic to a product $MA$ of a compact group $M$ and a group $A\cong \R_+^*$.

\subsubsection{}
Let $\Gamma\subset G$ be a torsion-free discrete subgroup.
\begin{ddd}\label{weih111}
A parabolic subgroup $P\subset G$ is called $\Gamma$-cuspidal
if $U_P:=\Gamma\cap P$ is an infinite subgroup such that
$l(U_P)\subset L$ is precompact, i.e. $l(U_P)\subset M$.
\end{ddd}
Let $p$ denote the $\Gamma$-conjugacy class of the
$\Gamma$-cuspidal parabolic subgroup $P$. We call such classes cusps
and say that the cusp $p$ has full rank (as opposed to
smaller rank) if $U_P\backslash \Omega_P$ is compact for one (and hence for any) $P\in p$.
For each cusp $p$ and $P\in p$
we form the manifold with boundary $\bar Y_{U_P}:= U_P\backslash  (X\cup \Omega_P)$.

\subsubsection{}\label{weih133}

The boundary $\partial X$ of $X$ has a $\Gamma$-equivariant decomposition
$\partial X=\Omega_\Gamma\cup\Lambda_\Gamma$
into a limit set $\Lambda_\Gamma$ and a domain of discontinuity
$\Omega_\Gamma$ (\cite{MR0336648}, Prop.~8.5). Let $\bar Y_\Gamma$ denote the manifold
with boundary $\bar Y_\Gamma:=\Gamma\backslash (X\cup \Omega_\Gamma)$.

\begin{ddd}\label{t67}
 The group $\Gamma$ is called geometrically finite if the following conditions hold:
\begin{enumerate}
\item The set $\cP_\Gamma$ of $\Gamma$-conjugacy classes of $\Gamma$-cuspidal parabolic subgroups
is finite.
\item There is a bijection $\ausg(\bar Y_\Gamma)\stackrel{\sim}{\rightarrow}
\cP_\Gamma$, where $\ausg(\bar Y_\Gamma)$ denotes the set of ends of
the manifold with boundary $\bar Y_\Gamma$. 
\item For all $p\in\cP_\Gamma$ and $P\in p$ there is a representative
$\bar Y_p$ of the end corresponding to $p$ and an isometric embedding
$e_P: \bar Y_p  \rightarrow \bar Y_{U_P}$ such that
its image $e_p(\bar Y_p)$ represents the end of $\bar Y_{U_P}$.
\end{enumerate}
\end{ddd}

By Bowditch \cite{bowditch951}, Corollary 6.3,  this definition
is equivalent to the slightly weaker definition \cite{bowditch951}, F1
(which can be derived from  \ref{t67} by replacing
``set $\cP_\Gamma$ of $\Gamma$-conjugacy classes of $\Gamma$-cuspidal parabolic subgroups'' by ``set $\ausg(\bar Y_\Gamma)$ of ends of $\bar Y_\Gamma$" in 1., ``bijection" by ``map $c$" in 2. and ``all $p\in\cP_\Gamma$" by 
``all $p\in\cP_\Gamma$ in the range of $c$" in 3.).

\subsubsection{}\label{weih132}

Let $Y_\Gamma:=\Gamma\backslash X$ denote the locally symmetric space associated to $\Gamma$.
By $\bar Y_0$ we denote the compact subset $\bar Y_\Gamma\setminus \bigcup_{p\in\cP_\Gamma} \bar Y_{p}$
of $\bar Y_\Gamma$.
The boundary of $\bar Y_\Gamma$ is the manifold
$B_\Gamma:=\Gamma\backslash \Omega_\Gamma$.
Let $\cP_\Gamma^<\subset \cP_\Gamma$ denote the subset of cusps of smaller rank.
The decompositions of $\bar Y_\Gamma$ into a compact piece and its ends induces
a decomposition $B_\Gamma=B_0\cup \bigcup_{p\in \cP^<_\Gamma} B_p$,
where $B_p:=B\cap \bar Y_p$, $p\in \cP_\Gamma^<\cup\{0\}$.

\subsubsection{}\label{holle}
Set $\cP_\Gamma^{max}:=\cP_\Gamma\setminus\cP_\Gamma^<$. It will be convenient to fix
a set $\tilde\cP$ of parabolic subgroups representing the elements of $\cP_\Gamma$.
It comes with a partition $\tilde\cP=\tilde\cP^<\cup\tilde\cP^{max}$.

\subsubsection{}\label{weih115}

Our analysis will require an additional assumption on the cusps.
By Lemma \ref{martin}  this assumption is automatically
satisfied in the cases $\RH^n$, $\CH^n$, but by Lemma \ref{contrmartin}
it is non-trivial in the case $\HH^n$. We are now going to describe
this assumption in detail.
Let $P\subset G$ be parabolic. A Langlands decomposition
$P=MAN$ is the same as a split $s:MA\rightarrow P$ of the extension (\ref{eepp}),
where we identify $M$ and $A$ with their images $s(M)$ and $s(A)$.

Let $p$ be a cusp of $\Gamma$ and $P\in p$. 
In Subsection \ref{cuspgeom} we will construct a Langlands
decomposition $P=MAN$, a discrete
subgroup $V\subset N$, and a homomorphism $m:N_V\rightarrow M$
from the Zariski closure $N_V$ of $V$ in $N$ such that
the group $U^0:=\{m(v)v|v\in V\}$ is a subgroup
of $U_P$ (see Def.~\ref{weih111}) of finite index.  
\begin{ddd}\label{t799}
The cusp $p$ is called regular if
the Langlands decomposition of $P$ above can be adjusted such that
$N_V$ is invariant with respect to conjugation by $A$.
We call such a Langlands decomposition adapted.
\end{ddd}

\begin{ass}
In the remainder of the present paper we assume that
$\Gamma$ is geometrically finite and that all its cusps are
regular.
\end{ass}

\subsection{Twists and bundles}

\subsubsection{}\label{weih112}

We fix a parabolic subgroup $P$ of $G$ and write $\partial X=G/P$.
We further fix a Langlands decomposition $P=MAN$ and let $\aaaa$
denote the Lie algebra of $A$.
Let $\naaa$ denote the Lie algebra of $N$, and let $\alpha\in\aaaa^*$
denote the short root of $\aaaa$ on $\naaa$. We define
$\rho\in\aaaa^*$ by $\rho(H):=\frac12 \tr(\ad(H)_{|\naaa})$, $H\in\aaaa$.

\subsubsection{}\label{weih135}
If $\lambda\in \aca$ and
$(\sigma,V_\sigma)$ is a representation of $M$, then we define the representation $(\sigma_\lambda,V_{\sigma_\lambda})$ of $P$
 by $V_{\sigma_\lambda}:=V_\sigma$ and $\sigma_\lambda(man):=a^{\rho-\lambda}\sigma(m)$.
By $1$ we denote the trivial one-dimensional representation of $M$.

\subsubsection{}
If $(\theta,V_\theta)$ is a representation of $P$, then we define
the $G$-homogeneous bundle $V(\theta):=G\times_P V_\theta$
and denote by $\pi^\theta$ the representation of $G$
on spaces of sections of $V(\theta)$.

\subsubsection{}

If $(\vp,V_\vp)$ is a finite-dimensional representation of $G$ (or $\Gamma$),
then we denote by $V(\theta,\vp)$ the tensor product
of $V(\theta)$ with the trivial bundle $\partial X\times V_\vp$,
and by $\pi^{\theta,\vp}$ the representation
of $G$ (or $\Gamma$) on spaces of sections of $V(\theta,\phi)$.
Note that we can identify
$$C^\infty(\partial X,V(\theta,\vp))\cong C^\infty(\partial X,V(\theta))\otimes V_\vp,\quad \pi^{\theta,\vp}=\pi^\theta\otimes\vp\ .$$

\subsubsection{}

The representation $\vp$  of $\Gamma$ is called twist.
Our analysis requires further assumptions on the twist going under the
name``admissible''.
Let $p$ be a cusp of $\Gamma$, $P\in p$, and $MAN$ be an adapted Langlands decomposition of $P$ (see Definition \ref{t799}). We define $M_U:=\overline{l(U_P)}$.
Then $P_U:=M_UN_V\subset P$ is a subgroup containing $U_P$.
Let $(\vp,V_\vp)$
be a twist.
\begin{ddd}\label{weih134}
The twist $(\vp,V_\vp)$ is called admissible at the cups $p$ if its restriction to $U_P$ extends to a continuous representation of $AP_U$ such that $A$ acts algebraically by
semisimple endomorphisms. A twist is called admissible
if it is admissible at all cusps of $\Gamma$.
\end{ddd}
Note that the  algebraic functions on $A$ are linear combinations of  $A\mapsto a^{n\alpha}$, $n\in \Z$
(see \ref{weih112} for the definition of $\alpha$).

\subsubsection{}

If $\Gamma$ has cusps, then the condition of admissibility excludes most unitary representations of $\Gamma$. Examples of admissible twists are restrictions
of finite-dimensional representations of $G$ to $\Gamma$.
As explained in Subsection \ref{ttw} twists are used 
as a technical device. 

\subsubsection{}

If $(\gamma,V_\gamma)$ is a representation of $K$, then we can form the
homogeneous bundle $V(\gamma):=G\times_KV_\gamma$ over $X$.
For a twist $\vp$ let $V(\gamma,\vp)$ denote the product
$V(\gamma)\otimes V_\vp$. Furthermore let $V_\Gamma(\gamma,\vp):=\Gamma\backslash V(\gamma,\vp)$
denote the corresponding bundle over $Y$.

\subsection{Exponents}

\subsubsection{}

The main numerical invariant associated to a discrete subgroup $\Gamma\subset G$
is its exponent $\delta_\Gamma\in \aaaa^*$. Note that each $g\in G$ can be written
as $g=ka_gh\in KA_+K$, where $a_g$ is uniquely determined. Here $A_+:=\{a\in A\mid a^\alpha\ge 1\}$.
\begin{ddd}
The exponent $\delta_\Gamma\in\aaaa^*$ is defined as the infimum
of the set $$\{\nu\in\aaaa^*\:|\:\sum_{g\in\Gamma} a_g^{-\nu-\rho}<\infty\}\ .$$
\end{ddd}
If follows from the discreteness of $\Gamma$ that 
$\delta_\Gamma\le\rho$.

\subsubsection{}
The critical exponent $\delta_\Gamma$ has been extensively studied, in particular by Patterson \cite{patterson762}, Sullivan \cite{sullivan79}, \cite{MR766265}, and Corlette \cite{corlette90}, Corlette and Iozzi \cite{corletteiozzi99}. From these papers we know that $\delta_\Gamma\in [-\rho,\rho]$, if $\Gamma$ is infinite. Moreover, we have $\delta_\Gamma=\rho$ if and only if $\Gamma$ has finite covolume.
If $\Lambda_\Gamma$ contains at least $2$ points, then $\delta_\Gamma+\rho=\dim_H(\Lambda)\alpha$, where $\dim_H(\Lambda)$ denotes the Hausdorff dimension of the limit set with respect to the natural class of sub-Riemannian metrics on $\partial X$. If $Y=\Gamma\backslash X$ is an infinite volume quotient of a quaternionic hyperbolic space or the
Cayley hyperbolic plane, then $\delta_\Gamma$ can not be arbitrary close to
$\rho$. In these cases we have $\delta_\Gamma\le (2n-1)\alpha$ and $\delta_\Gamma\le 5\alpha$, respectively \cite{corlette90},\cite{corletteiozzi99}.

\subsubsection{}

Let  $\vp$ be a twist.
\begin{ddd} We define the exponent $\delta_\vp\in\aaaa^*$
of $\vp$ as the infimum of the set $$\{\nu\in\aaaa\:|\:\sup_{g\in\Gamma} a_g^{-\nu} \|\vp(g)\|<\infty\}\ ,$$
where $\|.\|$ denotes any norm on $\End(V_\vp)$.
\end{ddd}

\subsection{Description of the main results}

\subsubsection{}
Let $\Gamma\subset G$ be a discrete, torsion-free, geometrically finite
subgroup such that all its cusps are regular. Furthermore let
$\vp$ be an admissible twist. Let $B_\Gamma:=\Gamma\backslash \Omega_\Gamma$ and
$V_{B_\Gamma}(\sigma_\lambda,\vp):=\Gamma\backslash V(\sigma_\lambda,\vp)$.
Under the name push-down we subsume several constructions related
to the average of elements of $C^\infty(\partial X,V(\sigma_\lambda,\vp))$
with respect to $\Gamma$. It is a crucial matter to construct
a space $B_\Gamma(\sigma_\lambda,\vp)$ which contains the result
of the average. If $\lambda\in\aca$ varies, then these spaces
assemble as a projective limit of locally trivial bundles
of Fr\'echet spaces in the sense of Subsection
\ref{llim}. Therefore we can speak of meromorphic families of continuous maps
from and to these spaces.

\subsubsection{}\label{weih123}

We define the bundle $B_\Gamma(\sigma_\lambda,\vp)$
as a direct sum
$$B_\Gamma(\sigma_\lambda,\vp):=B_\Gamma(\sigma_\lambda,\vp)_1 \oplus \cR_{\Gamma}(\sigma_\lambda,\vp)_{max}\ ,$$
where the first component satisfies
$$C_c^\infty(B_\Gamma,V_{B_\Gamma}(\sigma_\lambda,\vp))\subset B_\Gamma(\sigma_\lambda,\vp)_1\subset C^\infty(B_\Gamma,V_{B_\Gamma}(\sigma_\lambda,\vp))$$
(equality occurs if all cusps have full rank).
The second component is associated to the cusps of full rank and is finite-dimensional.

\subsubsection{}\label{neuj913}

The elements of $B_\Gamma(\sigma_\lambda,\vp)_1$ are characterized as smooth 
sections of $V_{B_\Gamma}(\sigma_\lambda,\vp)$ with a certain asymptotic expansion near the cusps.
We postpone the detailed description of the asymptotics until Subsection \ref{asz}.
The bundle $\cR_{\Gamma}(\sigma_\lambda,\vp)_{max}$ has a further decomposition
$$\cR_{\Gamma}(\sigma_\lambda,\vp)_{max}=\bigoplus_{P\in\tilde\cP^{max}} \cR_{U_P}(\sigma_\lambda,\vp)\ .$$
Recall from \ref{holle} that $\tilde \cP^{max}$ is in bijection with the set of cusps of full rank. In order to define
$\cR_{U_P}(\sigma_\lambda,\vp)$ we first consider the sheaf 
$\cE_{\infty_P}(\tilde\sigma,\tilde\vp)$ on $\aca$ of holomorphic families $\phi_\nu$, $\nu\in\aca$, of $U_P$-invariant distribution sections
of $V(\tilde\sigma_{\nu},\tilde\vp)$ supported in $\infty_P$.
Here $\tilde\sigma,\tilde\vp$ are the dual representations to
$\sigma,\vp$. The sheaf
$\cE_{\infty_P}(\tilde\sigma,\tilde\vp)$ is the sheaf of holomorphic sections of a trivial finite-dimensional holomorphic
vector bundle over $\aca$, and we denote by $E_{\infty_P}(\tilde\sigma_\lambda,\tilde\vp)$
the fibre of this bundle over $\lambda\in\aca$ 
(see Lemma \ref{vermaolbrich} for an alternative description of the space  
$E_{\infty_P}(\tilde\sigma_\lambda,\tilde\vp)$ and \ref{weih121} for the finite-dimensionality).
Then we define
$$\cR_{U_P}(\sigma_\lambda,\vp):=  E_{\infty_P}(\tilde\sigma_{-\lambda},\tilde\vp)^*\ .$$

\subsubsection{}

We can now state the definition of the push-down. 

\begin{ddd}\label{pushdowndef}
For $f\in C^\infty(\partial X,V(\sigma_\lambda,\vp))$
we define the push-down $\pi^\Gamma_*(f)\in B_\Gamma(\sigma_\lambda,\vp)$
as the direct sum of $\pi^\Gamma_*(f)_1\in B_\Gamma(\sigma_\lambda,\vp)_1$
and $\pi^\Gamma_*(f)_2\in \cR_{\Gamma}(\sigma_\lambda,\vp)_{max}$.
Here $\pi^\Gamma_*(f)_1$ is given by
$\pi^\Gamma_*(f)_1:=\sum_{g\in\Gamma} \pi^{\sigma_\lambda,\vp}(g) f_{|\Omega_\Gamma}$
(if the sum converges). The component
$\pi^\Gamma_*(f)_{2,P}\in \cR_{U_P}(\sigma_\lambda,\vp)$ for $P\in\tilde\cP^{max}$ is defined
by the condition that
$$\langle \phi,\pi^\Gamma_*(f)_{2,P}\rangle = \sum_{[g]\in\Gamma/U_P}\langle \pi^{\tilde\sigma_{-\lambda},\tilde\vp}(g)\phi,f \rangle$$
for all $\phi\in E_{\infty_P}(\tilde\sigma_{-\lambda},\tilde\vp)$
(if the sum converges).
\end{ddd}

\subsubsection{}

In \ref{pushdowndef} the push-down is defined if certain sums converge. The first part of the  following theorem gives a range of $\lambda\in \aca$ for which these convergence conditions are satisfied. 
The second part asserts its meromorphic continuation to all of $\aca$.

\begin{theorem}\label{t119}
If $f\in C^\infty(\partial X,V(\sigma_\lambda,\vp))$, then the push-down
$\pi^\Gamma_*(f)\in C^\infty(B_\Gamma,V_{B_\Gamma}(\sigma_\lambda,\vp))$ converges for $\Ree(\lambda)<-\delta_\Gamma-\delta_\vp$.
The push-down induces a meromorphic family on $\aca$ of continuous maps
$$\pi^\Gamma_*:C^\infty(\partial X,V(\sigma_\lambda,\vp))\rightarrow 
B_\Gamma(\sigma_\lambda,\vp)$$ with finite-dimensional singularities.
\end{theorem}

We will first prove the part of the theorem asserting the convergence of the
push-down for $\Ree(\lambda)<-\delta_\Gamma-\delta_\vp$. Then we consider
the adjoint of the push-down, the extension $ext^\Gamma$, 
and  obtain the meromorphic continuation of the latter (see Theorem
\ref{t110}).
The remainder of Theorem \ref{t119} then follows from this result by duality.

\subsubsection{}\label{weih126}

By $D_\Gamma(\sigma_\lambda,\vp)$ we denote the topological dual space
of $B_\Gamma(\tilde \sigma_{-\lambda},\tilde\vp)$.
It is a dual Fr\'echet and Montel space. 
By the latter property it is reflexive so that $$D_\Gamma(\sigma_\lambda,\vp)^\prime\cong
B_\Gamma(\tilde \sigma_{-\lambda},\tilde\vp)\ .$$
Varying $\lambda\in\aca$ the  spaces $D_\Gamma(\sigma_\lambda,\vp)$ 
form a direct limit of locally trivial holomorphic bundles of dual Fr\'echet
spaces in the sense of Subsection \ref{llim}. 

\begin{ddd}
For $\Ree(\lambda)>\delta_\Gamma+\delta_\vp$
the extension 
$$ext^\Gamma:D_\Gamma(\sigma_\lambda,\vp)\rightarrow 
C^{-\infty}(\partial X,V(\sigma_\lambda,\vp))$$ 
is defined to be the adjoint of 
$\pi^\Gamma_*:C^\infty(\partial X,V(\tilde \sigma_{-\lambda},\tilde\vp))\rightarrow B_\Gamma(\tilde \sigma_{-\lambda},\tilde\vp)$.
\end{ddd}

\begin{theorem}\label{t110}
The extension induces a meromorphic family of continuous maps on $\aca$
$$ext^\Gamma:D_\Gamma(\sigma_\lambda,\vp)\rightarrow 
C^{-\infty}(\partial X,V(\sigma_\lambda,\vp))$$
with finite-dimensional singularities and values in ${}^\Gamma C^{-\infty}(\partial X,V(\sigma_\lambda,\vp))$.
\end{theorem}

\subsubsection{}
The meromorphic continuation of the extension is the main goal
of the present paper. We obtain the continuation in close connection
with the meromorphic continuation of the scattering matrix
which will be defined below. In order to define this scattering matrix
we need the restriction map $res^\Gamma$ which is a left-inverse of the
extension. Due to the presence of cusps the definition of
the restriction map is more complicated than in \cite{MR1749869}.

\subsubsection{}
We first recall the definition of $res^\Gamma$ in the case
of convex-cocompact groups $\Gamma$ (i.e. in the case without cusps, compare Sec.~4 of \cite{MR1749869}). In this case
$$res^\Gamma:{}^\Gamma C^{-\infty}(\partial X,V(\sigma_\lambda,\vp))\rightarrow
C^{-\infty}(B_\Gamma,V_{B_\Gamma}(\sigma_\lambda,\vp))$$
is given by the composition of the restriction of distributions
to $\Omega_\Gamma$
and the identification ${}^\Gamma C^{-\infty}(\Omega_\Gamma,V(\sigma_\lambda,\vp))\cong  C^{-\infty}(B_\Gamma,V_{B_\Gamma}(\sigma_\lambda,\vp))$.
For the purpose of the present paper we employ a different
description.

\subsubsection{}

Let us still assume that $\Gamma$ is convex-cocompact.
We choose a cut-off function
$\chi^\Gamma\in C_c^\infty(\Omega_\Gamma)$ such that
$\sum_{g\in\Gamma} g^* \chi^\Gamma=1$.
Then we define $$\pi^{*}:
C^\infty(B_\Gamma,V_{B_\Gamma}(\tilde\sigma_{-\lambda},\tilde\vp))
\rightarrow C^\infty(\partial X,V(\tilde\sigma_{-\lambda},\tilde\vp))$$ by
$\pi^*(f)=\chi^\Gamma f$ (in order to understand the right-hand side properly one must identify sections of bundles over $B_\Gamma$ with $\Gamma$-equivariant sections on $\Omega_\Gamma$).
We define $$\widetilde{res}:C^{-\infty}(\partial X,V(\sigma_\lambda,\vp))
\rightarrow C^{-\infty}(B_\Gamma,V_{B_\Gamma}(\sigma_\lambda,\vp))$$
to be the adjoint of $\pi^*$. Then
$res^\Gamma$ coincides with the restriction of $\widetilde{res}$
to the subspace ${}^\Gamma C^{-\infty}(\partial X,V(\sigma_\lambda,\vp))$.
While $\widetilde{res}$ depends on the choice of $\chi^\Gamma$,
the restriction $res^\Gamma$ is independent of choices.
 
\subsubsection{}

If $\Gamma$ has cusps of smaller rank, then we can not assume
that $\chi^\Gamma$ has compact support. The definition of $\pi^*$ breaks down.
Our way arround that problem is as follows.
For each $k\in\nat_0$ we consider the
Banach spaces $C^k(\partial X,V(\tilde\sigma_{-\lambda},\tilde\vp))$
and define a meromorphic family of maps
$\pi^*_k:B_\Gamma(\tilde\sigma_{-\lambda},\tilde\vp)\rightarrow
C^k(\partial X,V(\tilde\sigma_{-\lambda},\tilde\vp))$.
The map $\pi^*_k$ is a right-inverse of the push-down $\pi^\Gamma_*$,
and its definition is similar in spirit to that
of $\pi^*$ above, but the details are more complicated
since we have to take into account the asymptotic expansions of the elements of
$B_\Gamma(\tilde\sigma_{-\lambda},\tilde\vp)$ near cusps.
One problem is that we can only deal with a finite number (depending on $k$) of these terms at a time in order to define $\pi_k^*$. The situation is similar to the problem of the construction
of smooth functions with given Taylor series at a given point. It is impossible to construct a continuous map from the space of formal power series to smooth functions which is right-inverse to the surjective map which takes the Taylor series.
Back to the definition of the restriction, let $C^{-k}(\partial X,V(\sigma_\lambda,\vp))$ denote the
distributions of order $k$, i.e. the topological dual of
$C^k(\partial X,V(\tilde\sigma_{-\lambda},\tilde\vp))$,
and define
$res_k^\Gamma:C^{-k}(\partial X,V(\sigma_\lambda,\vp))\rightarrow
D_\Gamma(\sigma_\lambda,\vp)$ as the adjoint of $\pi^*_k$.
The collection of maps $res_k^\Gamma$ play the role of the map
$\widetilde{res}$ above. These maps depend on choices
and are, in particular, not compatible if we vary $k$.

\subsubsection{}

Nevertheless, we have the following uniqueness property.
If $f_\nu\in {}^\Gamma C^{-\infty}(\partial X,V(\sigma_\nu,\vp))$,
$\nu\in\aca$, is a meromorphic family, and $W\subset \aca$
is compact, then for sufficiently large $k$
the meromorphic family $res^\Gamma_k(f_\nu)$, $\nu\in W$, is  well-defined. In fact,
due to the compactness of $W$ the order of $f_\nu$ as a distribution is uniformly bounded for $\nu\in W$.
In \ref{neuj3001} we will see that
$res^\Gamma_k(f_\nu)$ is independent of the choices. 

\subsubsection{}

In general the  restriction maps $res^\Gamma_k$ may have poles.
Moreover, not every element of ${}^\Gamma C^{-\infty}(\partial X,V(\sigma_\lambda,\vp))$
can be written as evaluation of a holomorphic family $f_\nu\in {}^\Gamma C^{-\infty}(\partial X,V(\sigma_\nu,\vp))$ defined near $\lambda$. But for generic $\lambda\in \aca$ the restriction
map $res^\Gamma_k$ is regular and every element $f\in {}^\Gamma C^{-\infty}(\partial X,V(\sigma_\lambda,\vp))$ extends to a family (in the present paper we will prove a weaker statement only, which suffices for our purposes).
In this case 
$res^\Gamma_k(f)$ is well-defined independent of the choice of (the sufficiently large) $k$,
and we can omit the subscript $k$ and write
$res^\Gamma(f)$ for $res^\Gamma_k(f)$.
We have the identity
$$res^\Gamma\circ ext^\Gamma=\id_{D_\Gamma(\sigma_\lambda,\vp)}\ ,$$
which holds true for generic $\lambda\in \aca$ or as an identity of  maps defined on meromorphic families
(see Lemma \ref{weih114}). 

\subsubsection{}

If $f\in {}^\Gamma C^{-k}(\partial X,V(\sigma_\lambda,\vp))$, $k\in\nat_0$,
then $f_{|\Omega_\Gamma}$ is a $\Gamma$-invariant distribution
on $\Omega_\Gamma$, hence can be considered as an element of
$C^{-\infty}(B_\Gamma,V_{B_\Gamma}(\sigma_\lambda,\vp))$. If $$i:D_\Gamma(\sigma_\lambda,\vp)
\rightarrow C^{-\infty}(B_\Gamma,V_{B_\Gamma}(\sigma_\lambda,\vp))$$
is the natural map
(adjoint to the inclusion $C_c^\infty(B,V_B(\tilde\sigma_{-\lambda},\tilde\vp))
\hookrightarrow B_\Gamma(\tilde\sigma_{-\lambda},\tilde\vp)$),
then $i\circ res^\Gamma(f)$ is defined (even though
$res^\Gamma$ may have a pole at $\lambda$)
and coincides with $f_{|\Omega_\Gamma}$.
In particular, the condition
$i\circ res^\Gamma(f)=0$ is  equivalent to $\supp(f)\subset \Lambda_\Gamma$.

If $\Gamma$ has cusps, then $i$ is not injective.
In this case  the condition $\supp(f)\in \Lambda_\Gamma$ does not imply
that $res^\Gamma(f)=0$ (provided the latter is defined).

\subsubsection{}
 
Vanishing of $res^\Gamma(f)$ is a rather strong condition
and should play the role of the condition $\supp(f)\subset \Lambda_\Gamma$ in the convex-cocompact case which was successfully exploited in \cite{MR1749869} and \cite{MR1689342}.
One application of this condition  is the following result
which is employed in proving the functional equation of the scattering matrix in the domain of convergence.
\begin{theorem}[Corollary \ref{wieder}]\label{micro}
Assume that $\Ree(\lambda)>\max\left(\{\delta_\Gamma\}\cup\{\rho_{U_P}\mid P\in\tilde\cP\}\right)+\delta_\vp$. If
$f_\mu\in {}^\Gamma C^{-\infty}(\partial X,V(1_\mu,\vp))$
is a germ of a meromorphic family at $\lambda$, then
$ext^\Gamma\circ res^\Gamma(f_\mu) = f_\mu$.
\end{theorem}

Further applications of the condition $res^\Gamma(f)=0$
are contained in \cite{MR1926489} (non-existence of such $f$ for $\Ree(\lambda)$ large) and  \cite{math.DG/0103144} (estimates of the regularity of $f$). In fact, the proof of 
Theorem \ref{micro} is based on the main result of  \cite{MR1926489}.

\subsubsection{}

We now turn to the scattering matrix. We assume that $\sigma$ is either
a Weyl-invariant irreducible representation of $M$ or
of the form $\sigma^\prime\oplus(\sigma^\prime)^w$, 
where $\sigma^\prime$ is irreducible and not Weyl-invariant,
and $(\sigma^\prime)^w(m):=\sigma^\prime(w^{-1}mw)$
is the Weyl-conjugate representation of $\sigma^\prime$.
Here $w\in N_K(A)$ is a representative of the non-trivial 
element of the Weyl group $W=N_K(A)/M\cong \Z/2\Z$.

\subsubsection{}\label{gaa101}

  Our starting point
is the scattering matrix
$$S^{\{1\}}_\lambda=J_\lambda:C^{-\infty}(\partial X,V(\sigma_\lambda,\vp)) \rightarrow
C^{-\infty}(\partial X,V(\sigma_{-\lambda},\vp))$$ associated to the trivial group $\{1\}$. In representation theory it firms  under the name 
Knapp-Stein intertwining operator.
Here we employ a  suitably normalized version as in the paper
\cite{MR1749869}, Sec. 5, to which we refer for further details.
Let $I_\aaaa\subset\aaaa^*$ be the $\Z$-module spanned by
$\alpha$ if $2\alpha$ is a root of $(\gaaa,\aaaa)$, and by
$\frac{1}{2}{\alpha}$ otherwise.
The operators $J_\lambda$ form a meromorphic family
of continuous maps with singularities in $I_\aaaa$
and satisfy the functional equation $J_\lambda\circ J_{-\lambda}=\id$.
 
\subsubsection{}

If $\Gamma$ is non-trivial, then we obtain the
scattering matrix $S_\lambda^\Gamma$ from $J_\lambda$
using restriction and extension.

\begin{ddd}
For $\Ree(\lambda)>\delta_\Gamma+\delta_\vp$ we define the scattering matrix
$$S^\Gamma_\lambda:D_\Gamma(\sigma_\lambda,\vp)\rightarrow D_\Gamma(\sigma_{-\lambda},\vp)$$
as the meromorphic
family of continuous maps
$$S^\Gamma_\lambda:=res^\Gamma\circ J_\lambda\circ ext^\Gamma\ .$$
\end{ddd}

\begin{theorem}\label{t113}
The scattering matrix has a meromorphic
continuation to all of $\aca$. It 
satisfies the functional equation $S^\Gamma_\lambda\circ S^\Gamma_{-\lambda}=\id$ and
the relation  $J_\lambda\circ ext^\Gamma=ext^\Gamma\circ S^\Gamma_\lambda$.
\end{theorem}

As noted above this theorem is proved in a multistep
procedure which involves the meromorphic continuation of $ext^\Gamma$
at the same time. The basic ideas are adapted
from \cite{MR1749869}. The non-compactness
of $B_\Gamma$ due to the presence of cusps of non-maximal
rank is responsible for the various complications which
have to be resolved on the way.

\subsubsection{}

The following application to the Eisenstein series is 
an easy consequence of the preceding
results and can be derived exactly as in \cite{MR1749869},
Cor.~10.2.
Let $\gamma$ be a finite-dimensional representation of $K$ and $T\in \Hom_M(V_\sigma,V_\gamma)$.
Then we have a holomorphic family of $G$-equivariant maps
$$P^T_\lambda: C^{-\infty}(\partial X,V(\sigma_\lambda,\vp))\rightarrow
C^\infty(X,V(\gamma,\vp))$$
which are called Poisson transformations. 
We refer to \cite{MR1749869}, Def. 4.8,  for a definition of the Poisson
transformation and to \cite{MR1749869}, Sec. 6, for a discussion of its properties.

\begin{ddd}\label{neuj1001}
If $\lambda\in\aca$, $f\in D_\Gamma(\sigma_\lambda,\vp)$,
$ext^\Gamma(f)$ is regular,
and $T\in \Hom_M(V_\sigma,V_\gamma)$, then
we define the Eisenstein series
$E(\lambda,f,T)\in C^\infty(Y,V_Y(\gamma,\vp))$
by $P^T_\lambda\circ ext^\Gamma(f)$.
\end{ddd}

In the special case that $\sigma$, $\gamma$, and $\vp$  are trivial, $T$ is
the identity, $\Ree(\lambda)>\delta_\Gamma$, and $f=\delta_b$ is the
delta-distribution located at $b\in B_\Gamma$, this definition
coincides with the classical definition of the Eisenstein series
as the $\Gamma$-average of the Poisson kernel, i.e. the
integal kernel of $P^{\id}_\lambda$.

Let $c_\sigma$ and $c_\gamma$ be the $c$-functions introduced in \cite{MR1749869}, Sec. 5.
We further employ the notation $(c_\gamma(\lambda) T)^w$ as introduced in the same section of \cite{MR1749869}.

\begin{kor}\label{main}
The Eisenstein series forms a meromorphic family of continuous maps
with finite-dimensional singularities
$$E(\lambda,.,T): D_\Gamma(\sigma_\lambda,\vp)    \rightarrow   C^\infty(Y,V_Y(\gamma,\vp))\ .$$
It satisfies the functional equation
$$E(\lambda,S_{-\lambda} f,T)=E(-\lambda,f,(\frac{c_\gamma(\lambda)}{c_\sigma(\lambda)}T)^w)\ .$$
\end{kor}

\section{Some machinery}
 
\subsection{Limits of bundles}\label{llim}

\subsubsection{}

In the present paper a main issue is the construction of
families of topological vector spaces $V_\lambda$, $\lambda\in\C$, and
the study of holomorphic (meromorphic) families of vectors 
$f_\lambda\in V_\lambda$ or holomorphic (meromorphic) families of 
homomorphisms $h_\lambda\in \Hom(V_\lambda,W_\lambda)$
between such families of spaces. In order to make the notion
of a holomorphic family precise we must relate the spaces $V_\lambda$
for neighbouring $\lambda$ in some holomorphic manner.
One way to do this is to equip the family of spaces with the 
structure of a locally trivial holomorphic bundle of topological vector
spaces. If we fix two local trivializations $U_0\times W_0$,
$U_1\times W_1$ of the family $V_\lambda$ with non-trivial overlap,
then the transition function is a holomorphic map from $U_0\cap U_1$ to
$\Hom(W_0,W_1)$ with values in the subspace of isomorphisms. 

\subsubsection{}

All topological vector spaces in the present paper will be locally convex.
Spaces of homomorphisms between topological vector spaces
will always be equipped with the topology of uniform convergence
on bounded sets. We refer to \cite{MR1689342}, Sec.2.2
for more details.

If we have two holomorphic families
$h_\lambda\in\Hom(W_0,W_1)$ and $g_\lambda\in\Hom(W_1,W_2)$,
then we can consider the composition $g_\lambda\circ h_\lambda\in \Hom(W_0,W_2)$. 
This composition is again holomorphic (see \cite{gloeckner}).
If the evaluation $W_0\to W^{\prime\prime}_0$ is continuous, then the adjoint
$h_\lambda^\prime\in \Hom(W_1^\prime,W_0^\prime)$ is holomorphic, too (see again \cite{gloeckner}).

In order to show that the composition of two holomorphic families of maps
is again holomorphic in \cite{MR1689342}, Sec.2.2, we assumed that $W_0$ is a Montel space. The discussion of \cite{gloeckner} shows that this assumption can be omitted. 

\subsubsection{}

In this paper we do not show that the spaces $B_\Gamma(\sigma_\lambda,\vp)$
($D_\Gamma(\sigma_\lambda,\vp)$) form locally trivial bundles of Fr\'echet
(dual Fr\'echet) spaces. These families of spaces arise instead as projective (symbol $\lim$) and inductive limits (symbol $\colim$) of locally trivial bundles, respectively. 
E.g. the bundle $B_U(\sigma_\lambda,\vp)$ comes as a limit $\lim B_{\Gamma,k}(\sigma_\lambda,\vp)$
of trivial bundles $B_{\Gamma,k}(\sigma_\lambda,\vp)$, 
$k\in \nat_0$. The connecting maps (in the example $B_{\Gamma,k}(\sigma_\lambda,\vp)\to B_{\Gamma,k+1}(\sigma_\lambda,\vp)$) of the diagrams are holomorphic families of continuous maps
which are injective and have dense range (resp. are injective).
But they are not compatible with the trivializations.

\subsubsection{}\label{weih211}

We must extend the  notions of a holomorphic or meromorphic
family of continuous of maps between such limits of locally trivial holomorphic bundles.
Let $$\dots\subset  E_{n+1}\subset E_{n} \subset\dots\ ,$$
$$\dots\subset  F_{n+1}\subset F_{n} \subset\dots$$
be decreasing families of locally trivial bundles of Fr\'echet spaces
defined over some open subset $U\subset \C$
such that all inclusion maps are holomorphic.
Let $E:=\lim E_n$, $F:=\lim  F_n$.
\begin{ddd}
A family of maps $\phi_z:E_z\rightarrow F_z$, $z\in U$,
is called holomorphic (meromorphic) if
for each $x\in U$ and $n\in \nat_0$ there exists
a neighbourhood $U_{x,n}\subset U$ of $x$, $m(n)\in\nat$, and 
a holomorphic (meromorphic) bundle map
$\Phi_n: (E_{m(n)})_{|U_{x,n}}\rightarrow (F_{n})_{|U_{x,n}}$ such that
$\phi_y$ is the restriction of $\Phi_n(y)$ to $E_y$
for all $y\in U_{n,x}$.
\end{ddd}
The composition of two holomorphic families is again a holomorphic family
(compare \cite{MR1749869}, Subsection 2.2).

\subsubsection{}

We now consider the dual situation.
Let $$\dots \rightarrow E^\prime_n\rightarrow E^\prime_{n+1}\rightarrow\dots\ ,$$
$$\dots \rightarrow F^\prime_n\rightarrow F^\prime_{n+1}\rightarrow\dots$$
be  direct systems of holomorphic locally trivial bundles
of dual Fr\'echet spaces defined over some open
subset $U\subset \C$. Let $E^\prime:=\colim  E^\prime_n$,
$F^\prime:=\colim F^\prime_n$.
\begin{ddd}
A family of maps $\phi_z:F_z^\prime\rightarrow E_z^\prime$, $z\in U$,
is called holomorphic (meromorphic) if
for each $x\in U$ and $n\in\nat_0$ 
 there exists
a neighbourhood $U_{x,n}\subset U$ of $x$, $m(n)\in\nat$,
and a holomorphic (meromorphic) bundle map
$\Phi_n^\prime: (F^\prime_{n})_{|U_{x,n}}\rightarrow (E^\prime_{m(n)})_{|U_{x,n}}$
such that for each $y\in U_{x,n}$
the restriction of $\phi^\prime_y$ to $(F^\prime_n)_y$
is the composition of $(\Phi_n^\prime)_y$ with the
natural map $(E^\prime_{m(n)})_y\rightarrow E^\prime_y$ .
\end{ddd}

Again the composition of two such families is a holomorphic (meromorphic)  family. The adjoint of a holomorphic family is again a holomorphic  (meromorphic) family.

\subsection{Embedding}\label{embedd}

\subsubsection{}

The symmetric spaces
$X$ considered in the present paper belong to a series of symmetric spaces.
Let $X^n$ denote the $n$'th space of the series. 
We will use the superscript
${}^n$ as a decoration of symbols for other objects in order to indicate
that they are associated to $X^n$. For a number of arguments we need
that $\delta_\Gamma^n$ is sufficiently negative. Note that
$\delta_\Gamma^n\to-\infty$ as $n\to\infty$. E.g., we will first obtain
a meromorphic continuation of $ext^{\Gamma,n+k}$ for sufficiently
large $k$, and then use the propositions below in order to conclude that
$ext^{\Gamma,n}$ has a meromorphic continuation, too.

\subsubsection{}

The main point of this subsection is to explain that the concept of embedding is compatible with the function and distribution spaces introduced in \ref{weih123} and \ref{weih126}. Furthermore we need compatibility with the push-down and extension maps. Note that the concept of embedding is only applied to the spherical case $\sigma=1$.

\subsubsection{}\label{neuj700}

We have embeddings $X^n\hookrightarrow X^{n+1}$ and
$i:\partial X^n\hookrightarrow \partial X^{n+1}$.
In order to have compatible embeddings $G^n\hookrightarrow G^{n+1}$
of the groups we assume
at this point that $G^n$ is one of $$\{Spin(1,n),SO(1,n)_0,SU(1,n),Sp(1,n)\}\ .$$
If we realize the last three groups of the list
using $\bF$-valued matrices, $\bF\in\{\R,\C,\Hh\}$, then the
embedding $G^n\hookrightarrow G^{n+1}$ is the usual embedding
into the left upper corner. We have compatible Iwasawa decompositions
such that $A=A^n=A^{n+1}$, and in particular $P^n\hookrightarrow P^{n+1}$.
Let $\zeta$ be $\alpha/2$, $\alpha$, and $2\alpha$ in the cases
$X=\RH^n$, $X=\CH^n$, and $X=\HH^n$, respectively.
Then $\zeta=\rho^{n+1}-\rho^n$, and we have
$(V_{1^{n+1}_\lambda})_{|P^n}=V_{1^n_{\lambda-\zeta}}$, and
hence $V(1^{n+1}_\lambda,\vp)_{|\partial X^n}=V(1^n_{\lambda-\zeta},\vp)$.

\subsubsection{}\label{weih301}

Let
$$i^*:C^\infty(\partial X^{n+1},V(1^{n+1}_\lambda,\vp))\rightarrow
C^\infty(\partial X^n, V(1^n_{\lambda-\zeta},\vp))$$
be the restriction of sections. It is a $\Gamma$-equivariant and surjective
continuous linear map.

\begin{prop}\label{rrttee}
There exists a meromorphic family of maps
$$i^*_\Gamma:B_\Gamma(1^{n+1}_\lambda,\vp)\rightarrow
B_\Gamma(1^n_{\lambda-\zeta},\vp)$$
such that the following diagram commutes:
$$\begin{array}{ccc} C^\infty(\partial X,V(1^{n+1}_\lambda,\vp)) &\stackrel{i^*}{\rightarrow}&C^\infty(\partial X^n, V(1^n_{\lambda-\zeta},\vp))\\
\downarrow \pi^{\Gamma,n+1}_*&&\downarrow \pi^{\Gamma,n}_*\\
B_\Gamma(1^{n+1}_\lambda,\vp)&\stackrel{i_\Gamma^*}{\rightarrow }& B_\Gamma(1^n_{\lambda-\zeta},\vp)
\end{array}\ .$$
If $\Gamma$ does not have cusps of full rank, then $i^*_\Gamma$ is holomorphic.
\end{prop}
This proposition will be proved in various stages. First we consider  a version for the  spaces associated to the cusps, see Lemma \ref{klopp}. The global result is stated in  Lemma \ref{comppp}.  

\subsubsection{}

In the following we explain the construction of $i_\Gamma^*$.
If all cusps of $\Gamma$ as a subgroup of $G^n$ have smaller rank, then
$i^*_\Gamma$ is induced
by the ususal restriction $i^*_\Gamma:
C^\infty(B_\Gamma,V_{B_\Gamma}(1_\lambda^{n+1},\vp))\rightarrow C^\infty(B_\Gamma,V_{B_\Gamma}(1^n_{\lambda-\zeta},\vp))$,
and both $\pi^{\Gamma,n}_*$
and $\pi^{\Gamma,n+1}_*$ are given as the average of sections
over $\Gamma$. The relation $i^*_\Gamma\circ \pi^{\Gamma,n+1}_*=\pi^{\Gamma,n}_*\circ i^*$
is obvious in the domain of convergence $\Ree(\lambda)>\delta^{n+1}_\Gamma+\delta_\vp$,
and follows by meromorphic continuation for all $\lambda\in\aca$.
We postpone the proof of the fact that $i^*_\Gamma$
really maps $B_\Gamma(1^{n+1}_\lambda,\vp)$
to $B_\Gamma(1^{n}_{\lambda-\zeta},\vp)$
until these spaces are defined.

\subsubsection{}

If $\Gamma$ has cusps of full rank (when considered as a subgroup of $G^n$), then $i^*_\Gamma$
is the sum of two maps $(i^*_\Gamma)_1:B_\Gamma(1^{n+1}_\lambda,\vp)
\rightarrow B_\Gamma(1^n_{\lambda-\zeta},\vp)_1$ and $(i^*_\Gamma)_2:
B_\Gamma(1^{n+1}_\lambda,\vp)\rightarrow \cR_\Gamma(1^n_{\lambda-\zeta},\vp)_{max}$.
The map $(i^*_\Gamma)_1$ is the restriction of smooth sections, and
we can apply the argument above in order to show commutativity of
the part of the commutative
diagram involving $(\pi^{\Gamma,n}_*)_1$ and $(i^*_\Gamma)_1$.

\subsubsection{}

We define $(i^*_\Gamma)_2$ such that the corresponding diagram
is commutative. We choose $k\in\nat_0$
so large that $E_{\infty_{P^n}}(1^n_{-\lambda+\zeta},\tilde\vp)
\subset C^{-k}(\partial X^n,V(1^n_{-\lambda+\zeta},\tilde\vp))$
for all $P^n\in\tilde\cP^{max,n}$. For $f\in B_\Gamma(1^{n+1}_\lambda,\vp)$
the component $i^*_\Gamma(f)_{2,P^n}
\in \cR^{U_{P^n}}(1^n_{\lambda-\zeta},\vp)$ is characterized
by $$\langle\phi,i^*_\Gamma(f)_{2,P^n}\rangle=\langle res^{U_{P^{n+1}},n+1}_k\circ i_*(\phi), T_{P^{n+1}}(f) \rangle$$
for all $\phi\in E_{\infty_{P^n}}(1^n_{-\lambda+\zeta},\tilde\vp)$,
where $P^{n+1}\subset G^{n+1}$ is the unique parabolic subgroup
containing $P^n$, and
$$T_{P^{n+1}}:C^\infty(B^{n+1}_\Gamma,V_{B_\Gamma}(1^{n+1}_\lambda,\vp))\rightarrow C^\infty(B^{n+1}_{U_{P^{n+1}}},V_{B^{n+1}_{U_{P^{n+1}}}}(1^{n+1}_\lambda,\vp))$$
is the natural map defined using the map $e_{P^{n+1}}$
and a cut-off function  $\chi_p\in C^\infty(B_\Gamma)$
which is supported in $B_p$ and one on a neighbourhood of infinity
of $B_p$.  
This definition is independent of the choice of $k$.
In order to show that $(i^*_\Gamma)_2\circ \pi^{\Gamma,n+1}_*=(\pi^{\Gamma,n}_*)_2\circ i^*$
we compute in the domain of convergence for generic $\lambda$ (so that the restrictions are defined)
\begin{eqnarray*}
\langle\phi,(i^*_\Gamma)_2\circ  \pi^{\Gamma,n+1}_*(f)\rangle&=&
\langle res^{U_{P^{n+1}},n+1}_k\circ i_*(\phi),T^*_{P^{n+1}}\circ \pi^{\Gamma,n+1}_*(f)\rangle\\
&=&\sum_{g\in\Gamma^{P^{n+1}}} \langle res^{U_{P^{n+1}},n+1}_k\circ i_*(\phi),\chi_{P^{n+1}}\circ \pi^{U_{P^{n+1}},n+1}_*(\pi^{1^{n+1}_\lambda,\vp}(g)f)\rangle\\
&\stackrel{(*)}{=}&\sum_{g\in\Gamma^{P^{n+1}}} \langle i_*(\phi), \pi^{1^{n+1}_\lambda,\vp}(g)f\rangle\\
&=&\sum_{g\in\Gamma^{P^{n+1}}} \langle \pi^{1^{n+1}_{-\lambda+\zeta},\tilde\vp}(g^{-1}) \phi, i^*(f)\rangle\\
&=&\langle \phi,(\pi^{\Gamma,n}_*)_2\circ i^*(f)\rangle\ ,
\end{eqnarray*} 
where $\Gamma^{P^{n+1}}$ denotes any system of representatives
of $U_{P^{n+1}}\backslash \Gamma$, and $\chi_{P^{n+1}}:= e_{P^{n+1}}^* \chi_p$ (see Def.~\ref{t67} for the map $e_p$).
In order to get the equality marked by $(*)$ we use that
$$(1- \chi_{P^{n+1}})\circ res^{U_{P^{n+1}},n+1}_k \circ i_*(\phi)=0$$ and
$$ext^{U_{P^{n+1}},n+1}\circ res_k^{U_{P^{n+1}},n+1}\circ i_*(\phi)=i_*(\phi)\ .$$
\hB

\subsubsection{}

Let $i_*$ and $i_*^\Gamma$ denote the dual maps to
$i^*$ and $i^*_\Gamma$. Dualizing Prop.~\ref{rrttee} we obtain the following corollary.
\begin{kor}\label{rrttee1}
We have a commutative diagram of meromorphic families of maps
$$\begin{array}{ccc} C^{-\infty}(\partial X,V(1^{n}_\lambda,\vp)) &\stackrel{i_*}{\rightarrow}&C^{-\infty}(\partial X^{n+1}, V(1^{n+1}_{\lambda-\zeta},\vp))\\
\uparrow ext^{\Gamma,n}_*&&\uparrow ext^{\Gamma,n+1}_*\\
D_\Gamma(1^{n}_\lambda,\vp)&\stackrel{i^\Gamma_*}{\rightarrow }&D_\Gamma(1^{n+1}_{\lambda-\zeta},\vp)
\end{array}\ .$$
\end{kor}

\subsubsection{}


\begin{lem}\label{jjjjjjjjjjjtwzew}
There exists holomorphic families of continuous maps
$$j^*:C^{-\infty}(\partial X^{n+1},V(1^{n+1}_{\lambda-\zeta},\vp))\rightarrow C^{-\infty}(\partial X^n, V(1^n_\lambda,\vp))$$
and $$j_*:
C^{\infty}(\partial X^n, V(1^n_{-\lambda},\tilde\vp))\rightarrow 
C^{\infty}(\partial X^{n+1},V(1^{n+1}_{-\lambda+\zeta},\tilde\vp))$$
such that $j^*\circ i_* =\id$, $i^*\circ j_*=\id$.
\end{lem}
\proof
We choose a tubular neighbourhood $T:\bF\times\partial X^n\hookrightarrow \partial X^{n+1}$
of $\partial X^n$ such that $T(\{0\}\times \partial X^n)=\partial X^n$.
Furthermore we choose a cut-off function $\chi\in C_c^\infty(\bF)$
such that $\chi(0)=1$. Then we define the map
$t:
C^{\infty}(\partial X^n, V(1^n_{\rho^n},\tilde\vp))\rightarrow
C^{\infty}(\partial X^{n+1},V(1^{n+1}_{\rho^{n+1}},\tilde\vp))$
setting $(tf)(T(u,x)):=\chi(u)f(x)$ and extending
this function by zero outside the tubular neighbourhood.
Here we use the canonical identifications
$C^{\infty}(\partial X^n, V(1^n_{\rho^n},\tilde\vp))=C^{\infty}(\partial X^n)\otimes V_{\tilde\vp}$
and $C^{\infty}(\partial X^{n+1},V(1^{n+1}_{\rho^{n+1}},\vp))=C^{\infty}(\partial X^{n+1})\otimes V_{\tilde\vp}$.

We choose a positive section
$s_{n+1}\in  C^{\infty}(\partial X^{n+1},V(1^{n+1}_{\rho^{n+1}+\alpha}))$ and put $i^* s_{n+1}=:s_n\in  C^{\infty}(\partial X^n, V(1^n_{\rho^n+\alpha}))$.
Then we define
$\Phi_{\rho^n}:C^{\infty}(\partial X^n, V(1^n_{-\lambda},\tilde\vp))
\rightarrow C^{\infty}(\partial X^n, V(1^n_{\rho^n},\tilde\vp))$
and
$\Phi_{\rho^{n+1}}:
C^{\infty}(\partial X^{n+1},V(1^{n+1}_{-\lambda+\zeta},\tilde\vp))\rightarrow
C^{\infty}(\partial X^{n+1},V(1^{n+1}_{\rho^{n+1}},\tilde\vp))$
as multiplication by $s_n^{(\rho^n+\lambda)/\alpha}$ and
$s_{n+1}^{(\rho^{n+1}+\lambda-\zeta)/\alpha}$, respectively.
We define $j_*:=\Phi_{\rho^{n+1}}^{-1}\circ t\circ \Phi_{\rho^n}$,
and $j^*$ as the adjoint of $j_*$. It is now easy to check
that $j_*$ and $j^*$ have the required properties. \hB

\subsubsection{}\label{ohu}

If $G^n$ does not belong to the list
$\{Spin(1,n),SO(1,n)_0,SU(1,n),Sp(1,n)\}$, then there is a finite covering
$p:\tilde G^n\rightarrow G^n$ with $\tilde G^n\in \{Spin(1,n),SO(1,n)_0,SU(1,n),Sp(1,n)\}$.
In this case can find a normal subgroup
$\Gamma^0\subset \Gamma$ of finite index and a discrete subgroup $\tilde \Gamma^0\subset \tilde G^n$ such that $p$ induces an isomorphism from $\tilde \Gamma^0$ to $\Gamma^0$. Indeed, using Selberg's Lemma we can take
a torsion-free subgroup $\tilde\Gamma^0$  of $p^{-1}(\Gamma)$ of finite index
and set $\Gamma^0:=p(\tilde\Gamma^0)$.

We can apply the concept of embedding to
the subgroup $\tilde\Gamma^0$. In order to transfer results for
$\tilde\Gamma^0$ to $\Gamma$ we use averages over the finite group $\Gamma/\Gamma^0$.

\subsection{Twisting}\label{ttw}
  
\subsubsection{}

Twisting is an important technical device of the present paper. We explain in \ref{weih125} and \ref{weih1266} how this concept is used. But first we must introduce some notation.

\subsubsection{}\label{weih122}

If $(\pi,V_\pi)$ is a finite-dimensional representation
of $G$, then we can form
the bundles $V(\sigma_\lambda\otimes\pi,\vp)$ and
$V(\sigma_\lambda,\pi\otimes \vp)$.
There is an isomorphism $T:V(\sigma_\lambda\otimes\pi,\vp)\stackrel{\sim}{\rightarrow} 
V(\sigma_\lambda,\pi\otimes \vp)$, which is given on the level of sections by
 $T:C^\infty(\partial X,V(\sigma_\lambda\otimes\pi,\vp))
\stackrel{\sim}{\rightarrow} C^\infty(\partial X,V(\sigma_\lambda,\pi\otimes\vp))$, $T(f)(g):=\pi(g)f(g)$.
 
\subsubsection{}\label{weih144}
Given an irreducible representation $\sigma$ of $M$
and $\mu_0\in\aaaa^*$   there exists a finite-dimensional representation $(\pi_{\sigma,\mu},V_{\pi_{\sigma,\mu}})$
of $G$ with highest $A$-weight $\mu\ge\mu_0$ such that $V_\sigma=V_{\pi_{\sigma,\mu}}(\mu)$
as representations of $M$. Here $V_{\pi_{\sigma,\mu}}(\mu)$ denotes the subspace
of $V_{\pi_{\sigma,\mu}}$ on which $A$ acts with weight $\mu$.
Note that $V_{\pi_{\sigma,\mu}}(-\mu)= V_{\sigma^w}$
as representations of $M$.

\subsubsection{}\label{weih241}

If $\sigma$ is a representation of $M$ of the form
$\sigma^\prime\oplus (\sigma^\prime)^w$, where $\sigma^\prime$
irreducible and not equivalent to its Weyl conjugate, then given $\mu_0\in\aaaa^*$
there exist
$(\pi_{\sigma^\prime,\mu},V_{\pi_{\sigma^\prime,\mu}})$ and $(\pi_{(\sigma^\prime)^w,\mu},V_{\pi_{(\sigma^\prime)^w,\mu}})$
for suitable $\mu\ge \mu_0$. In this case we set
$\pi_{\sigma,\mu}:=\pi_{\sigma^\prime,\mu} \oplus \pi_{(\sigma^\prime)^w,\mu}$.
In connection with twisting without further notice we will always assume that
$\sigma$ is Weyl-invariant and either irreducible or of the form
$\sigma^\prime\oplus (\sigma^\prime)^w$, where $\sigma^\prime$
is irreducible and not Weyl invariant. 
Note that  $\pi_{\tilde\sigma,\mu}=\tilde\pi_{\sigma,\mu}$.

\subsubsection{}\label{weih142}

We have an embedding of $P$-modules
$V_{\sigma_\lambda}\hookrightarrow V_{1_{\lambda+\mu}}\otimes V_{\pi_{\sigma,\mu}}$
and a corresponding embedding of $\Gamma$-equivariant bundles
$V(\sigma_\lambda,\vp)\hookrightarrow V(1_{\lambda+\mu}\otimes\pi_{\sigma,\mu},\vp)$.
Composing this embedding with the isomorphism $T$ (introduced in \ref{weih122})
we obtain the $\Gamma$-equivariant embedding
$$i_{\sigma,\mu}:C^{\pm\infty}(\partial X,V(\sigma_\lambda,\vp))\hookrightarrow
C^{\pm\infty}(\partial X,V(1_{\lambda+\mu},\pi_{\sigma,\mu}\otimes\vp))\ .$$ 

Similarly we have a projection of $P$-modules
$V_{1_{\lambda-\mu}}\otimes V_{\pi_{\sigma,\mu}} \rightarrow V_{\sigma_{\lambda}}$ and a corresponding projection of $\Gamma$-equivariant bundles
$V(1_{\lambda-\mu}\otimes\pi_{\sigma,\mu},\vp)\rightarrow V(\sigma_{\lambda},\vp)$. Composing this projection with the inverse of $T$ 
we obtain the $\Gamma$-equivariant projection
$$p_{\sigma,\mu}:C^{\pm\infty}(\partial X,V(1_{\lambda-\mu},\pi_{\sigma,\mu}\otimes\vp))
\rightarrow   C^{\pm\infty}(\partial X,V(\sigma_{\lambda},\vp))\ .$$

\subsubsection{}\label{neuj1000}

Since $i_{\sigma,\mu}$ and $p_{\sigma,\mu}$ are induced by homomorphisms
of $\Gamma$-equivariant bundles restriction to $\Gamma$-equivariant
sections over $\Omega_\Gamma$ provides
the embedding 
$$i^\Gamma_{\sigma,\mu}:C^{\infty}(B_\Gamma,V_{B_\Gamma}(\sigma_\lambda,\vp))\hookrightarrow
C^{\infty}(B_\Gamma,V_{B_\Gamma}(1_{\lambda+\mu},\pi_{\sigma,\mu}\otimes\vp))$$ 
and the projection
$$p^\Gamma_{\sigma,\mu}:C^{\infty}(B_\Gamma,V_{B_\Gamma}(1_{\lambda-\mu},\pi_{\sigma,\mu}\otimes\vp))
\rightarrow   C^{\infty}(B_\Gamma,V_{B_\Gamma}(\sigma_{\lambda},\vp))\ .$$  
In Lemma \ref{compat3} will show that $i^\Gamma_{\sigma,\mu}$ (resp. $p^\Gamma_{\sigma,\mu}$) 
maps $B_\Gamma(\sigma_{\lambda},\vp)_1$ (resp. $B_\Gamma(1_{\lambda-\mu},\pi_{\sigma,\mu}\otimes\vp)_1$) 
to $B_\Gamma(1_{\lambda+\mu},\pi_{\sigma,\mu}\otimes\vp)_1$
(resp. $B_\Gamma(\sigma_{\lambda},\vp)_1$).
Hence we can consider the
adjoint of
$$p^\Gamma_{\tilde\sigma,\mu}:B_\Gamma(1_{-\lambda-\mu},\pi_{\tilde\sigma,\mu}\otimes\tilde\vp)_1\rightarrow B_\Gamma(\tilde\sigma_{-\lambda},\tilde\vp)_1$$ which will be denoted by 
by $$(i^\Gamma_{\sigma,\mu})_1:D_\Gamma(\sigma_{\lambda},\vp)_1
\rightarrow D_\Gamma(1_{\lambda+\mu},\pi_{\sigma,\mu}\otimes\vp)_1\ .$$

In a similar manner the adjoint of
$$i^\Gamma_{\tilde\sigma,\mu}:B_\Gamma(\tilde\sigma_{-\lambda},\tilde\vp)_1\rightarrow B_\Gamma(1_{-\lambda+\mu},\pi_{\tilde\sigma,\mu}\otimes\tilde\vp)_1$$
will be  a map
$$(p^\Gamma_{\sigma,\mu})_1:D_\Gamma(1_{\lambda-\mu},\pi_{\sigma,\mu}\otimes\vp)_1\rightarrow D_\Gamma(\sigma_{\lambda},\vp)_1\ .$$

\subsubsection{}

In the presence of cusps of full rank we must take into account their contribution to the function spaces (see \ref{weih123}). It is clear that $i_{\sigma,\mu}$ (resp. $p_{\sigma,\mu}$)  maps
$E_{\infty_P}(\sigma_\lambda,\vp)$ (resp. $E_{\infty_P}(1_{\lambda-\mu},\pi_{\sigma,\mu}\otimes\vp)$) to  $E_{\infty_P}(1_{\lambda+\mu},\pi_{\sigma,\mu}\otimes\vp)$ (resp.
$E_{\infty_P}(\sigma_{\lambda},\vp)$).  
We get second components
\begin{eqnarray*}
(i^\Gamma_{\sigma,\mu})_2&:&\cR_{\Gamma}(\sigma_\lambda,\vp)_{max}\to \cR_{\Gamma}(1_{\lambda+\mu},\pi_{\sigma,\mu}\vp)_{max}\\(p^\Gamma_{\sigma,\mu})_2&:&\cR_\Gamma(1_{\lambda-\mu},\pi_{\sigma,\mu}\otimes\vp)_{max} \rightarrow \cR_\Gamma(\sigma_{\lambda},\vp)_{max}\end{eqnarray*}
Combining all these maps we obtain 
maps
\begin{eqnarray*}
i^\Gamma_{\sigma,\mu}&:&D_\Gamma(\sigma_{\lambda},\vp)\rightarrow D_\Gamma(1_{\lambda+\mu},\pi_{\sigma,\mu}\otimes\vp)\\
p^\Gamma_{\sigma,\mu}&:&D_\Gamma(1_{\lambda-\mu},\pi_{\sigma,\mu}\otimes\vp)\rightarrow D_\Gamma(\sigma_{\lambda},\vp)
\end{eqnarray*}
and dually corresponding maps between spaces
of functions.

\subsubsection{}

Using the isomorphism $T$ we transfer the action  
$\pi^{1_{-\lambda-\mu}\otimes\pi_{\tilde\sigma,\mu},\id_{\tilde\vp}}$ of $\cZ(\gaaa)$ to $$C^\infty(\partial X,V(1_{-\lambda-\mu},\pi_{\tilde\sigma,\mu}\otimes \tilde\vp))\ .$$ Since this action commutes with $\Gamma$ and
is implemented by local operators we obtain
an action of $\cZ(\gaaa)$ on 
$C^\infty(B_\Gamma,V_{B_{\Gamma}}(1_{-\lambda-\mu},\pi_{\tilde\sigma,\mu}\otimes \tilde\vp))$ as well.
We will show that $B_\Gamma(1_{-\lambda-\mu},\pi_{\tilde\sigma,\mu}\otimes\tilde\vp)_1$
is a $Z(\gaaa)$-invariant subspace of $C^\infty(B_\Gamma,V_{B_{\Gamma}}(1_{-\lambda-\mu},\pi_{\tilde\sigma,\mu}\otimes \tilde\vp))$. By duality we obtain an action of $Z(\gaaa)$
on $D_\Gamma(1_{\lambda+\mu},\pi_{\sigma,\mu}\otimes\vp)_1$.
Since $\cZ(\gaaa)$ also acts on 
$E_{\infty_P}(1_{\lambda+\mu},\pi_{\sigma,\mu}\otimes\vp)$, $P\in \tilde\cP^{max}$, the space
$D_\Gamma(1_{\lambda+\mu},\pi_{\sigma,\mu}\otimes\vp)$
has the structure of a $\cZ(\gaaa)$-module.

\subsubsection{}\label{weih146}

 Let $\Omega\in\cZ(\gaaa)$ denote the Casimir operator.
Let $\{\sigma^i | i\in I_\nu\}$, be the set of Weyl-invariant
representations of $M$ (irreducible or sum of two non-Weyl invariant
irreducible)  occuring in the restriction
of $V_{\pi_{\sigma,\mu}}(\nu)$ to $M$, where we set $I_\mu:=\{1\}$.
Let
$\nu_i(\lambda):=\pi^{\sigma^i_{\lambda+\nu}}(\Omega)\in\C$
be the eigenvalue of the Casimir operator on the principal series representation
$\pi^{\sigma^i_{\lambda+\nu}}$. We set $I:=\bigcup_{\nu\not=\mu} I_\nu$
and define the meromorphic function 
$$\aca\ni\lambda\mapsto \Pi(\lambda):=1-\prod_{i\in I} \frac{\Omega-\nu_i(\lambda)}{\nu_1(\lambda)-\nu_i(\lambda)}\in\cZ(\gaaa)\ .$$
This function has a finite number of poles all contained in the subset
$I_\aaaa\subset \aaaa^*$.
Furthermore, we define
$$Z(\lambda):=\pi^{1_{\lambda+\mu}\otimes\pi_{\sigma,\mu},\id_{\vp}}(\Pi(\lambda))\ .$$
If $\lambda\not\in I_\aaaa$ (or more precisely, if $\Pi(\lambda)$ is regular),
then $Z(\lambda)$ is a projection. 

\subsubsection{}

For $\lambda\in  I_\aaaa$ we will show in Lemma \ref{compat6}
that 
$$i_{\sigma,\mu}^\Gamma:D_\Gamma(\sigma_{\lambda},\vp)\rightarrow 
\ker \{Z(\lambda):D_\Gamma(1_{\lambda+\mu},\pi_{\sigma,\mu}\otimes\vp)\rightarrow
D_\Gamma(1_{\lambda+\mu},\pi_{\sigma,\mu}\otimes\vp)\}=:\ker_\Gamma(Z(\lambda))$$
is an isomorphism. We then define the continuous map
$$j_{\sigma,\mu}^\Gamma:D_\Gamma(1_{\lambda+\mu},\pi_{\sigma,\mu}\otimes\vp)
\stackrel{1-Z(\lambda)}{\longrightarrow} \ker_\Gamma(Z(\lambda))\stackrel{(i_{\sigma,\mu}^\Gamma)^{-1}}{\longrightarrow}
D_\Gamma(\sigma_{\lambda},\vp)\ .$$
In Lemma  \ref{compat6} we will furthermore show that as a function of $\lambda$ the maps
$j_{\sigma,\mu}^\Gamma$ form a meromorphic family of continuous maps
which is holomorphic on $\aca\setminus I_\aaaa$.
If $\Gamma$ is trivial, then we omit the superscript
and write $j_{\sigma,\mu}$ for $j_{\sigma,\mu}^{\{1\}}$.

\subsubsection{}\label{weih125}

We now explain by examples how the concept of twisting is applied.
Assume that we have obtained a meromorphic continuation of the extension
$ext^\Gamma:D_\Gamma(1_\lambda,\vp)\rightarrow C^{-\infty}(\partial X,V(1_\lambda,\vp))$ to some half-plane 
$W=\{\Ree(\lambda)>\lambda_0\}$ for all admissible twists $\vp$.
Then we can obtain a meromorphic continuation of 
$ext^\Gamma:D_\Gamma(\sigma_\lambda,\vp)\rightarrow
C^{-\infty}(\partial X,V(\sigma_\lambda,\vp))$ to some larger half plane
$W-\mu_0$, $\mu_0\ge 0$, and for all $\sigma$ as follows. We choose $(\pi_{\sigma,\mu}, V_{\pi_{\sigma,\mu}})$ for some $\mu\ge \mu_0$.
Then we can define a meromorphic family 
$\widetilde{ext}^\Gamma:D_\Gamma(\sigma_\lambda,\vp)\rightarrow
C^{-\infty}(\partial X,V(\sigma_\lambda,\vp))$
for $\lambda\in W-\mu_0$ by the diagram
$$\begin{array}{ccc}
D_\Gamma(\sigma_\lambda,\vp)&\stackrel{i^\Gamma_{\sigma,\mu}}{\rightarrow}&D_\Gamma(1_{\lambda+\mu},\pi_{\sigma,\mu}\otimes\vp))\\
\downarrow \widetilde{ext}^\Gamma&&\downarrow ext^\Gamma\\
C^{-\infty}(\partial X,V(\sigma_\lambda,\vp))&\stackrel{j_{\sigma,\mu}}{\leftarrow}& 
C^{-\infty}(\partial X,V(1_{\lambda+\mu},\pi_{\sigma,\mu}\otimes\vp))
\end{array}\ .$$
We then show that $ext^\Gamma=\widetilde{ext}^\Gamma$
for all $\lambda$ with sufficiently large real part.
Thus $\widetilde{ext}^\Gamma$ provides a meromorphic continuation
of $ext^\Gamma$.

\subsubsection{}\label{weih1266}

Assume that we have defined the scattering matrix
$S^\Gamma_\lambda:D_\Gamma(1_\lambda,\vp)\rightarrow D_\Gamma(1_{-\lambda},\vp)$
in the spherical case for all admissible twists $\vp$.
Then we could define the scattering matrix 
$S^\Gamma_\lambda:D_\Gamma(\sigma_\lambda,\vp)\rightarrow D_\Gamma(\sigma_{-\lambda},\vp)$ for all $\sigma$ using the diagram
\begin{equation}\label{qw11}
\begin{array}{ccc}
D_\Gamma(\sigma_\lambda,\vp)&\stackrel{i^\Gamma_{\sigma,\mu}}{\rightarrow}&D_\Gamma(1_{\lambda+\mu},\pi_{\sigma,\mu}\otimes\vp)\\
\downarrow S^\Gamma_\lambda&&\downarrow S^\Gamma_{\lambda+\mu}\\
D_\Gamma(\sigma_{-\lambda},\vp)&\stackrel{p^\Gamma_{\sigma,\mu}}{\leftarrow}&D_\Gamma(1_{-\lambda-\mu},\pi_{\sigma,\mu}\otimes\vp)
\end{array}\ .
\end{equation}
For generic $\lambda\in\aca$ we have a well-defined map
$res^\Gamma:{}^\Gamma C^{-\infty}(\partial X,V(1_{\lambda+\mu},\pi_{\sigma,\mu}\otimes\vp)) 
\rightarrow D_\Gamma(1_{\lambda+\mu},\pi_{\sigma,\mu}\otimes\vp)$,
and we can define
$res^\Gamma:{}^\Gamma C^{-\infty}(\partial X,V(\sigma_\lambda,\vp)) 
\rightarrow D_\Gamma(\sigma_\lambda,\vp)$ using the diagram
$$\begin{array}{ccc}
{}^\Gamma C^{-\infty}(\partial X,V(\sigma_\lambda,\vp)) &\stackrel{i_{\sigma,\mu}}{\rightarrow}&{}^\Gamma C^{-\infty}(\partial X,V(1_{\lambda+\mu},\pi_{\sigma,\mu}\otimes\vp))\\
\downarrow res^\Gamma &&\downarrow res^\Gamma \\
D_\Gamma(\sigma_{\lambda},\vp)&\stackrel{j^\Gamma_{\sigma,\mu}}{\leftarrow}&D_\Gamma(1_{\lambda+\mu},\pi_{\sigma,\mu}\otimes\vp)
\end{array}\ .$$
We then check that 
$S^\Gamma_\lambda$ defined by (\ref{qw11}) coincides with
$res^\Gamma\circ J_\lambda\circ ext^\Gamma$.

\subsection{Holomorphic vector bundles over $\C$}

\subsubsection{}

It is a general fact
that each finite-dimensional
holomorphic vector bundle over $\C$ or a half-plane is trivial.
Let $U\subset \C$ be open and $E$ be a finite-dimensional holomorphic vector bundle
over $U$. Furthermore let $U\times V\rightarrow U$ be a trivial
bundle of Fr\'echet spaces over $U$. By $\cE$ and $\cV$ we denote
the corresponding sheaves of holomorphic sections.
Let $\phi:E\rightarrow V$ be a meromorphic family of linear maps.
The assertions of the following lemma are well-known in the case where $V$ is finite-dimensional.  
\begin{lem}\label{bunle}
\begin{enumerate}
\item
There exists a unique finite-dimensional subbundle $F\subset U\times V$
such that $\phi$ factors over a meromorphic family of maps $\psi:E\rightarrow F$
which is surjective for generic $z\in U$.
\item
If we have  meromorphic families of bundle maps
$Z_E:E\rightarrow E$ and $Z_V:V\rightarrow V$ such that $\phi\circ Z_E=Z_V\circ \phi$, then $Z_V$ restricts to a meromorphic family of
bundle maps of $F$.
\item
Furthermore, if $U=\C$ or a half-plane, then there exists a meromorphic right-inverse
$\eta:F\rightarrow E$ such that $\phi\circ \eta=\id_F$.
\end{enumerate}
\end{lem}
\proof
The proof of this lemma will occupy the remainder of the present subsection.
The main point which we will explain in detail is the reduction to the finite-dimensional case.

\subsubsection{}
We first show that $F$ exists.
Let $\phi$ have singularities in the discrete subset
$A\subset U$, and let $n_z$, $z\in A$,
denote the order of the corresponding singularity.
We consider the divisor $D:=\sum_{z\in A} n_z z$, the associated line bundle
$L(D)$, and its sheaf of sections $\cL(D)$.
Let $\cV(D)$ be the sheaf of sections of $V\otimes L(D)$.
Then $\phi$ induces a holomorphic map
$\phi_D:E\rightarrow V\otimes L(D)$. Let $\cW$ be the quotient sheaf
$$\cE \stackrel{\phi_D}{\rightarrow}\cV(D)\rightarrow \cW\rightarrow 0$$
and denote by $Tor(\cW)$ its torsion subsheaf. We define
$\cF(D)$ as the kernel
$$0\rightarrow \cF(D)\rightarrow \cV(D)\rightarrow \cW/Tor(\cW)\rightarrow 0\ .$$
There is a natural factorization of $\phi_D$ over $\psi_D:\cE\rightarrow \cF(D)$.

\subsubsection{}\label{weih129}

We claim that $\cF(D)$ is a coherent sheaf.
Let $x\in U$.
If $v\in \cV(D)_x$, then let $lp(v)\in V$ denote
the leading part of $v$ at $x$. Note that $lp(v)=0$ implies that $v=0$ since a Laurent series must start somewhere.

We define the subspace
$Z\subset V$ as the set of all leading parts $lp(\phi_D(e))$, $e\in \cE_x$.

\subsubsection{}

We now show that $\dim(Z)<\infty$.
We choose a connected  neighbourhood $U_x\subset U$
of $x$ and a holomorphic trivialization
$E_{|U_x}=U_x\times E_x$. Let $\phi=\sum_{n\ge m} \phi_n z^n$
denote the Laurent expansion of $\phi$ in this trivialization, where
$\phi_n\in \Hom(E_x,V)$. Then using that $\dim E_x<\infty$ we have
$$Z=\phi_m(E_x)+\phi_{m+1}(\ker\phi_m)+\phi_{m+2}(\ker\phi_m\cap\ker\phi_{m+1})+\dots + \phi_{m+k}(\bigcap_{i=0}^{k-1} \ker\phi_{m+i})\ ,$$
where $k\in\nat_0$ is sufficiently large such that
$\bigcap_{i=0}^{k-1} \ker\phi_{m+i}=\bigcap_{i=0}^{l-1} \ker\phi_{m+i}$
for all $l\ge k$. This proves that $\dim(Z)<\infty$.

\subsubsection{}

We choose some closed subspace $Y\subset V$ of finite codimension
such that $Z\cap Y=\{0\}$.
We consider the composition 
$\overline{\phi}_{D}:E\stackrel{\phi}\rightarrow V \rightarrow V/Y$.
We form the quotient of sheaves
$$\cE \stackrel{\overline{\phi}_D}{\rightarrow}\cV(D)/\cY(D)\rightarrow\tilde \cW\rightarrow 0\ ,$$
and we define $\tilde \cF(D)$ as the kernel
$$0\rightarrow\tilde\cF(D)\rightarrow \cV(D)/\cY(D)\rightarrow\tilde \cW/Tor(\tilde \cW)\rightarrow 0\ .$$

\subsubsection{}
Consider the following commutative diagram of sheaves on $U_x$:
$$\begin{array}{ccccccccc}
&&0&&0&&0&&\\
&&\downarrow&&\downarrow&&\downarrow&&\\
0&\rightarrow&0&\rightarrow&\cY(D)&\stackrel{\cong}{\rightarrow}& \cY(D)&\rightarrow &0\\
&&\downarrow&&\downarrow&&\downarrow&&\\
0&\rightarrow&\cE/\ker \phi_D&\rightarrow&\cV(D)&\rightarrow& \cW&\rightarrow &0\\
&&\downarrow\cong &&\downarrow&&\downarrow&&\\
0&\rightarrow&\cE/\ker \overline{\phi}_D&\rightarrow&\cV(D)/\cY(D)&\rightarrow& \tilde\cW&\rightarrow &0\\
&&\downarrow&&\downarrow&&\downarrow&&\\
&&0&&0&&0&&\end{array}
$$
The two lower rows are exact by construction.
We now show that the left column is exact.
We have a sequence
$$\ker(\phi_D)\subset \ker(\overline{\phi_D})\subset \cE$$ of inclusioins of torsion-free sheaves.
In particular, $\ker(\phi_D)$ and $\ker(\overline{\phi_D})$ are sheaves of holomorphic sections of vector bundles on $U_x$. If we show that the inclusion induces an isomorphism
\begin{equation}\label{weih128}
\ker(\phi_D)_x\cong  \ker(\overline{\phi_D})_x\ ,
\end{equation}
then we conclude an isomorphism of sheaves
$\ker(\phi_D) \cong \ker(\overline{\phi_D})$ after shrinking $U_x$, if necessary.
To see (\ref{weih128}) consider  $h\in \ker (\overline{\phi}_D)_x$. Then $lp(\phi_D(h))\in Y$.
Since on the other hand $lp(\phi_D(h))\in Z$ we conclude that
$lp(\phi_D(h))=0$, hence $h\in\ker (\phi_D)_x$.


 We now conclude that
the last column is exact, too. 

\subsubsection{}

We now consider the diagram
$$\begin{array}{ccccccccc}
&&0&&0&&0&&\\
&&\downarrow&&\downarrow&&\downarrow&&\\
0&\rightarrow&0&\rightarrow&Tor \cW&\stackrel{\cong}{\rightarrow}& Tor\tilde \cW &\rightarrow &0\\
&&\downarrow&&\downarrow&&\downarrow&&\\
0&\rightarrow&\cY(D)&\rightarrow&\cW&\rightarrow&\tilde \cW&\rightarrow &0\\
&&\downarrow\cong &&\downarrow&&\downarrow&&\\
0&\rightarrow&\cY(D)&\rightarrow&\cW /Tor \cW&\rightarrow&\tilde \cW/ Tor\tilde\cW&\rightarrow &0\\
&&\downarrow&&\downarrow&&\downarrow&&\\
&&0&&0&&0&&\end{array}\ .
$$
Since $\cY(D)$ is torsion-free we have $\cY(D)\cap Tor\cW
=0$ This implies the isomorphism in the upper row.
It follows that all rows and columns of this diagram are exact.
%
%

\subsubsection{}

We now consider the following diagram of sheaves on $U_x$:
$$\begin{array}{ccccccccc}
&&0&&0&&0&&\\
&&\downarrow&&\downarrow&&\downarrow&&\\
0&\rightarrow&0&\rightarrow&\cY(D)&\stackrel{\cong}{\rightarrow}& \cY(D)&\rightarrow &0\\
&&\downarrow&&\downarrow&&\downarrow&&\\
0&\rightarrow&\cF(D)&\rightarrow&\cV(D)&\rightarrow& \cW/Tor \cW&\rightarrow &0\\
&&\downarrow\cong &&\downarrow&&\downarrow&&\\
0&\rightarrow&\tilde\cF(D)&\rightarrow&\cV(D)/\cY(D)&\rightarrow& \tilde\cW/Tor \tilde\cW&\rightarrow &0\\
&&\downarrow&&\downarrow&&\downarrow&&\\
&&0&&0&&0&&\end{array}
$$
The rows and the middle column are exact. We just have shown
that the right column is exact. 

\subsubsection{}

We conclude that the natural map $\cF(D)\rightarrow \tilde \cF(D)$
is an isomorphism of sheaves. 
This proves the claim \ref{weih129} since
$\tilde \cF(D)$ is obviously coherent.

\subsubsection{}

Since $\cF(D)$ is torsion-free it is the sheaf of sections
of a holomorphic vector bundle $F(D)$ (here we use the fact that
the base space is smooth and one-dimensional).
By construction $\cF(D)/\psi_D(\cE)$ is torsion and therefore
$\psi_{D,x}:E_x\rightarrow F(D)_x$ is surjective for
generic $x\in U$.
We define $F:=F(D)\otimes L(-D)$ and let $\psi:E\rightarrow F$
be the corresponding meromorphic family of maps.
The uniqueness part and assertion 2.
are left to the reader.

\subsubsection{}

We now construct the meromorphic right-inverse $\eta$.
Note that we can assume that $E$ and $F$ are trivial.
We fix trivializations and a constant Hermitian metric
on $E$. Let $z\in \C$ be such that $\phi:E_z\rightarrow F_z$
is regular and surjective. Let $P$ be the orthogonal projection onto the
orthogonal complement of $\ker\phi_z$.
We extend $P$ constantly over $U$ and let $P(E)$ be the
range of $P$.
The composition $\phi\circ P:P(E)\rightarrow F$ is invertible at $z$
and therefore has a meromorphic family of inverses
$\eta:F\rightarrow P(E)\subset E$.
\hB

\subsubsection{}

\begin{ddd}\label{imagebundle}
The bundle $F$ constructed in Lemma \ref{bunle} is called the
image bundle of $\phi$.
\end{ddd}

\section{Pure cusps}

\subsection{Geometry of cusps}\label{cuspgeom}

\subsubsection{}

In this subsection we analyze the geometry of cusps.
First we recall the following theorem of Auslander.
\begin{theorem}[\cite{auslander61}, \cite{MR1490024}] \label{ausl}
Let $N$ be a connected, simply-connected nilpotent Lie group and $M$
be a compact group of automorphisms of $N$.
Furthermore let $U\subset N\rtimes M$ be a discrete subgroup
and $U^*:=\overline{NU}_0\cap U=p_{|U}^{-1}(\overline{p(U)}_0)$,
where $p:N\rtimes M\rightarrow M$ is the projection,
$()_0$ stands for connected component of the identity,
and "$\bar{.}$" means the closure of the set in the argument.
Then
\begin{enumerate}
\item $U^*\subset U$ is a normal subgroup of  finite index.
\item There exists $b\in N$ and a connected subgroup
$N_V\subset N$ such that $(U^*)^b:=bU^*b^{-1}$ acts
effectively and cocompactly by translations on 
$N_V$, i.e. there is lattice $V^*\subset N_V$ and an isomorphism
$\theta:(U^*)^b\rightarrow V^*$ such that
$u.x = \theta(u)x$, $x\in N_V$, where $(N\rtimes M)\times N\ni (a,x)\mapsto a.x\in N$
denotes the natural action of $N\rtimes M$ on $N$. 
\item The element $b\in N$ of 2. can be choosen
such that $U^b$ leaves the space $N_V$ invariant.
\item $M_{U^0}:=\overline{p(U^*)}\subset M$ is a torus.
\end{enumerate}
\end{theorem}
The third assertion is due to Apanasov \cite{MR1490024}.

\subsubsection{}

We now apply Theorem \ref{ausl}
to a discrete torsion-free subgroup $U\subset P$
of a parabolic subgroup $P\subset G$ such that 
$P$ is $U$-cuspidal. We have the exact sequence
\begin{equation}\label{yy77}
0\rightarrow N\rightarrow P\stackrel{l}{\rightarrow} L\rightarrow 0\ ,\end{equation}
where $N$ is the nil-radical of $P$. The group
$L$ decomposes as $L=M\times A$, and we denote
by $l_M:P\rightarrow M$ the composition of
$l$ with the projection from $L$ to $M$.
Let $E:=l^{-1}(M\times\{1\})$.
If we choose a split $s_1:M\rightarrow E$ of the sequence
$0\rightarrow N\rightarrow E\stackrel{l_M}{\rightarrow} M\rightarrow 0$,
then $m\in M$ acts by the automorphism $s_1(m)$ on $N$, and 
we have an isomorphism $T_{s_1}:E\rightarrow N\rtimes_{s_1}M$.
Since $P$ is $U$-cuspidal we have $U\subset E$.
Applying  Theorem \ref{ausl} to $T_{s_1}(U)$ we obtain a subgroup
$T_{s_1}(U)^*\subset T_{s_1}(U)$ of finite index, $b\in N$, a connected
subgroup $N_V$, a lattice $V^*\subset N_V$, and an isomorphism
$\theta:bT_{s_1}(U)^*b^{-1}\rightarrow V^*$
such that $bT_{s_1}(U)b^{-1}$ leaves $N_V$ invariant and acts by translations on $N_V$ via $\theta$.

\subsubsection{}

 If we replace the split $s_1$ and $N_V$ by $b^{-1}N_Vb$ and
$s:=b^{-1}s_1b$, then $T_s(U)$  itself leaves $N_V$
invariant, and $T_s(U)^*$ acts by translations on $N_V$.

From now on we use the split $s$ in order to identify $M$ with the
subgroup $s(M)\subset E\subset P$, to write $E$ as the product $NM$,
and to indentify $T_s(U)$ with $U$.

\subsubsection{}\label{weih131}

Any element $u\in U^*$ ($u\in U$) can be written as $n_um_u\in NM$
such that $n_u\in N_V$ and $m_u$ centralizes (normalizes) $N_V$.

We have $V^*:=\{n_u|u\in U^*\}\subset N_V$.
The map $V^*\ni n_u \mapsto m_u\in M$ defines a homomorphism
$m^*:V^*\rightarrow M_{U^0}$.
\begin{lem}
There exists a subgroup $V\subset V^*$ of finite index
such that the restriction $m:=m^*_{|V}$ extends to a homomorphism
$m:N_V\rightarrow M_{U^0}$.
\end{lem}
\proof
Let $\pi:V^*\rightarrow N_V/[N_V,N_V]$ be the projection.
Then $\pi(V^*)\subset N_V/[N_V,N_V]$ is a lattice (see
e.g. \cite{raghunathan72}, proof of Thm. 2.10).
Let $g_1,\dots g_n\in V^*$, $n=\dim(N_V/[N_V,N_V])$
be such that $\pi(g_1),\dots,\pi(g_n)$ generate the lattice
$\pi(V^*)$. We define $V:=\langle g_1,\dots,g_n\rangle$.
Furthermore we put $X_i:=\log(g_i)\in\naaa_V$
and consider $\pi(X_i)\in \naaa_V/[\naaa_V,\naaa_V]$.
Then $\pi(X_i)$ is a basis of $\naaa_V/[\naaa_V,\naaa_V]$,
and $\{X_1,\dots,X_n\}$ generate the Lie algebra $\naaa_V$.
We conclude that $N_V$ is the smallest connected group
containing $V$, and by \cite{raghunathan72}, Ch.2,  $V$ is a lattice 
in $N_V$. Therefore $V\subset V^*$ has finite index.

We claim that $V\cap [N_V,N_V]=[V,V]$.
Let $g\in V\cap [N_V,N_V]$.
Then there are finite sequences $i_k\in\{1,\dots,n\}$ and $e_k\in\{1,-1\}$, $k=1,\dots r$,
such that $g=g_{i_1}^{e_1}\dots g_{i_r}^{e_r}$.
Applying $\pi$ we obtain
$1=\pi(g_{i_1})^{e_1}\dots \pi(g_{i_r})^{e_r}$.
Using the fact that $N_V/[N_V,N_V]$ is free abelian we conclude
that there is a permutation $\sigma\in S_r$ such that
$1=g_{i_{\sigma(1)}}^{e_{\sigma(1)}}\dots g_{i_{\sigma(r)}}^{e_{\sigma(r)}}$.
We conclude that
$1=g$ modulo $[V,V]$, i.e. $g\in [V,V]$. This proves the claim.

We first define the derivative $dm:\naaa_V\rightarrow \maaa_{U^0}$
of $m$ as the composition 
$$\naaa_V\stackrel{\pi}{\rightarrow} \naaa_V/[\naaa_V,\naaa_V]\stackrel{q}{\rightarrow} \maaa_{U^0}\ ,$$
where $q$ is given by $q(\pi(X_i)):= Z_i$, $i=1,\dots n$,
where $Z_i\in\maaa_{N_V}$ is any element satisfying $\exp(Z_i)=m(g_i)$.
Integrating the derivative we obtain
a representation $\tilde m:N_V\rightarrow M_{U^0}$
such that $m(g_i)=\tilde m(g_i)$, $i=1\dots,n$.
Since $M_{U^0}$ is a torus we see that
$m$ must factor over $[V,V]$ and conclude that $m=\tilde m$.
This finishes the proof of the lemma.
\hB 
 
\subsubsection{}\label{weih251}

We consider the subgroup of finite index $U^1:=\{v m(v)| v\in V\}\subset U$.
\begin{ddd}\label{tu0}
We define $U^0$ to be the largest normal subgroup of $U$
such that $U^0\subset U^1\subset U$, i.e.
$$U^0:=\bigcap_{u\in U} (U^1)^u\ .$$
\end{ddd}
$U^0\subset U$ has finite index, too.
The torus $M_{U^0}$ coincides with $\overline{m(N_V)}$.
\begin{ddd}\label{weih1422}
We set $P_{U^0}:=N_VM_{U^0}$. Furthermore we define
$M_U:=\overline{l_M(U)}$ and $P_U:=N_VM_U$. 
\end{ddd}
Then $U\subset P_U$, $M_{U^0}\subset M_U$ has finite index, and $M_U$
normalizes $N_V$.

\subsubsection{}

We have an exact sequence
$$0\rightarrow E\rightarrow P\stackrel{l_A}{\rightarrow} A\rightarrow 0 \ ,$$
where the decomposition $E=NM$ is already fixed by the split $s$.
The split extends uniquely to a split
$s:MA\rightarrow P$ of (\ref{yy77}) and induces a split
 of the sequence above.
In particular we
obtain a Langlands decomposition $P=NAM$, where
we identify $A$ with its image by $s$. The split
furthermore defines an action of $A$ on $N$ by
automorphisms commuting with the action of $M$.

\subsubsection{}

We call the cusp associated to $U\subset P$ regular
if we can choose $s$ such that $N_V$ is $A$-invariant (compare Definition \ref{t799}).

\begin{lem}\label{martin}
If $X$ is a real or complex hyperbolic space,
then every cusp is regular.
\end{lem}
\proof
If $X$ is real-hyperbolic, then for any split $s$ the group $A$
acts on $\naaa$ as multiplication by scalars. Any subspace,
in particular $\naaa_V$, is invariant with respect to $A$.

We now consider the case that $X=\CH^n$.
Assume that we have chosen a split $s$. Then we can decompose
$\naaa=\naaa_{\alpha}\oplus\naaa_{2\alpha}$ with respect to the action
of $A$ such that $a\in A$ acts on $\naaa_{i\alpha}$ as multiplication
by $a^{i\alpha}$. Here $\naaa_{\alpha}$ is a symplectic
vector space of dimension $2(n-2)$,
where the symplectic form 
with values in the one-dimensional space $\naaa_{2\alpha}$
is given by the commutator.
Taking $M_{U^0}$-invariants we obtain a decomposition
$\naaa^{M_{U^0}}=\naaa_1\oplus\naaa_2$,
where $\naaa_1\subset \naaa_{\alpha}$ is a 
symplectic subspace (\cite{guilleminsternberg77}, Prop. 4.2.1). If $\naaa_1=\{0\}$,
then $\naaa_V=\naaa_2$ is invariant with respect to $A$,
and we are done. We now assume that
$\naaa_1\not=\{0\}$. 
We have   two cases. If $\naaa_2\subset \naaa_V$,
then $\naaa_V=\naaa_V\cap\naaa_1\oplus \naaa_2$,
and this space is invariant with respect to $A$.
Thus we assume that $\naaa_2\not\subset\naaa_V$.
Since $\dim(\naaa_2)=1$ we then have  $\naaa_2\cap \naaa_V=\{0\}$.
Let $p_i:\naaa_V\rightarrow \naaa_i$ denote the projections.
Then $p_1:\naaa_V\rightarrow \naaa_1$ is injective.
Let $\lambda:\naaa_1\rightarrow \naaa_2$ be a linear extension
of $p_2\circ p_1^{-1}:p_1(\naaa_V)\rightarrow \naaa_2$. 
Since the symplectic form is non-degenerate 
there exists a unique $Y\in \naaa_1$ such that
$\lambda(X)=[Y,X]$ for all $X\in\naaa_1$.

Let $h:=\exp(-Y)\in N^{M_{U^0}}$.
We claim that $\naaa_V^h\subset \naaa_1$.
Indeed, if $X\in \naaa_V$,
then $p_2(X^h)=p_2(X-[Y,X])=\lambda(p_1(X))-[Y,X]=[Y,X]-[Y,X]=0$.

If we replace the split $s:M\times A\rightarrow P$
by $hsh^{-1}$, then $\naaa_V$ gets replaced by $\naaa_V^h$.
Thus by an appropriate choice of the split
we can assume that $\naaa_V\subset\naaa_1$, and we are done.
\hB
 
\subsubsection{}

The next lemma shows that regularity of cusps is a proper restriction in the case of $\HH^n$.

\begin{lem}\label{contrmartin}
If $X=\HH^n$, $n\ge 2$, then there exist non-regular
cusps.
\end{lem}
\proof
It suffices to provide an example for $n=2$.
We can identify
$\naaa=\Hh\oplus \Imm_\Hh(\Hh)$ ($\Imm_\Hh(\Hh)$ denoting the imaginary quaternions),
and the commutator 
of $X,Y\in\Hh$ is given by 
$[X,Y]= \bar X Y-\bar Y X \in \Imm_\Hh(\Hh)$.
Implicitly we have fixed some split $s\times t$
such that $a\in A$ acts on $\Hh$ as multiplication by
$a^\alpha$ and on $\Imm_\Hh(\Hh)$ by $a^{2\alpha}$. 
Let $1,I,J,K$ be a base of the copy of $\Hh$ and $i,j,k$
be the base of $\Imm_\Hh(\Hh)$.
Then we consider $\naaa_V=\spann_\R\{1,I+j,i\}$
and let $U\subset N_V$ be any lattice.
It is not possible to conjugate this subspace
into an $A$-invariant one.
\hB

\subsubsection{}

Assume that $U\subset P$ defines a regular cusp.
\begin{ddd}\label{weih215}
We define $\rho_U\in\aaaa^*$ by $\rho_U(H):=\frac12 \tr(\ad(H)_{|\naaa_V})$,
where $\naaa_V$ is the Lie algebra of $N_V$. Furthermore, we set
 $\rho^U:=\rho-\rho_U$.
\end{ddd}

\subsection{Schwartz spaces}\label{schw}

\subsubsection{}

In this subsection we start the description of our function spaces. Here we introduce the Schwartz spaces $\cS_{U}(\sigma_\lambda,\vp)$. We will show that these spaces form trivial holomophic bundles of Fr\'echet spaces for $\lambda\in \aca$, and that they are compatible with twisting. The Schwartz spaces are basic building blocks of the function spaces $B_U(\sigma_\lambda,\vp)$ which we will define later
(see Subsection \ref{asz}) .

\subsubsection{}

We fix a representative $w\in K$ of the non-trivial element of the Weyl group $W(\gaaa,\aaaa)$ such that $w^2=1$. Let $\infty_P\in\partial X$
be the fixed point of $P$ and $0_P:=w\infty_P$. 
The subset  
$\Omega_U=\Omega_P=\partial X\setminus\infty_P$
is the $N$-orbit of $0_P$ (see Subsection \ref{geomf} for notation). Indeed, the map
$N\ni x\mapsto xw\infty_P\in \Omega_U$ is a diffeomorphism from $N$ to $\Omega_U$.

\subsubsection{}

Let $(\sigma,V_\sigma)$ be a finite-dimensional unitary representation 
of $M$. Given $\lambda\in \aca$ we form the representation $(\sigma_\lambda,V_{\sigma_\lambda})$ of $P$ (see \ref{weih135}). Furthermore,
let $(\vp,V_\vp)$ be an admissible twist (see Definition \ref{weih134}).

\subsubsection{}\label{weih243}

The action of $A$ on $\cU(\naaa)$ induces a grading.
We choose a basis $\{A_i\}$ of $\cU(\naaa)$ such that the subset
$\{A_i\:|\: \deg(A_i)\le d\alpha\}$ spans $\cU(\naaa)^{\le d\alpha}$.
Furthermore, we choose a norm  $|.|$  on $V_\vp$. 
For any $d\in\nat_0$, $k\in \R$, and compact subset $W\subset N\setminus N_V$  (see Theorem \ref{ausl} for the definition of $N_V$) we define
the seminorm $q_{W,d,k}$ on $C_c^\infty(B_{U},V_{B_{U}}(\sigma_\lambda,\vp))$
by
$$q_{W,d,k}(f):=\sup_{\{i|\deg(A_i)\le d\alpha\}} \sup_{x\in W}\sup_{a\in A_+}
a^{k\alpha-2(\lambda-\rho^U)} |\vp(a)^{-1}  f(x^aA_iw)|\ .$$
Here we consider a section of $V_{B_{U}}(\sigma_\lambda,\vp)$ as a
$U$-invariant function on $G\setminus P$ with values in $V_{\sigma_\lambda}\otimes V_\vp$
satisfying the corresponding invariance conditions with respect to the
right $P$-action, and where $u\in U$ acts by
$(\pi^{\sigma_\lambda,\vp}(u)f)(xw)=\vp(u) f(u^{-1}xw)$.

\subsubsection{}\label{weih259}

We fix $W$ and another compact subset $W^\prime\subset N$ such that
$W^\prime\cup W^{A_+}$ projects surjectively onto $B_{U}$.
We add an arbitrary $C^d$-norm over $W^\prime$ to $q_{W,d,k}$
in order to obtain a norm $\|.\|_{k,d}$.
We first define $S_{U,k,d}(\sigma_\lambda,\vp)$ to be
the Banach space closure of $C_c^\infty(B_{U},V_{B_{U}}(\sigma_\lambda,\vp))$
with respect to $\|.\|_{k,d}$. This Banach space is independent
(up to equivalent norms)
of the choices of $W$, $W^\prime$, $|.|$, and the base of $\cU(\naaa)$.
\begin{ddd}\label{weih136}
We define $\cS_{U,k}(\sigma_\lambda,\vp)$ to be the intersection of
the spaces $\cS_{U,k,d}(\sigma_\lambda,\vp)$ for all $d\in\nat_0$.
The space $\cS_{U,k}(\sigma_\lambda,\vp)$ is a Fr\'echet space which
is topologized by the countable set of norms $\|.\|_{k,d}$, $d\in\nat_0$.
\end{ddd}
The space $\cS_{U,k}(\sigma_\lambda,\vp)$ is a space of smooth sections with a fixed growth rate at infinity
of $B_U$ measured by $k\in \nat_0$.

\subsubsection{}
 
We now show that the family of spaces $\{\cS_{U,k}(\sigma_\lambda,\vp)\}_{\lambda\in\aca}$
forms a trivial holomorphic bundle of Fr\'echet spaces. To this end we construct a certain holomorphic family of functions $s^z$, $z\in \C$. In the proof of Lemma \ref{ytra} we use multiplication by these functions in order to identify
the Schwartz spaces for different $\lambda\in \aca$. 

\subsubsection{}

\begin{lem}\label{ingf}
There exists a positive $P_{U}$-invariant function
$s\in C^\infty(N)$ such
that for any compact subset $W\subset N\setminus N_V$ there is $a_0\in A$ with
$s(x^a)=a^{2\alpha}s(x)$ for all $a\ge a_0$ and $x\in W$.
\end{lem}
\proof 
Let $\bar{N}:=N^w$ and consider the decomposition
$$G\setminus wP=\bar NMAN\ , \quad g=\bar n(g) m(x)a(x)n(x)\ .$$
Note that $a(xw)$ satisfies $a(x^aw)=a^2a(xw)$ for all $a\in A$ and $x\in N\setminus\{1\}$ (see (\ref{weih141}) below).
We consider the function $N\setminus\{1\}\ni x\mapsto a(xw)^{\alpha}$ as a positive
function $s_1\in C^\infty(N\setminus \{1\})$.
Let $s_2\in C^\infty(N)$ be any positive function.
Using a partition of unity we glue $s_1$ and $s_2$
to obtain a positive function $s_3\in C^\infty(N)$
which coincides with $s_1$ outside of a compact subset $W_1$ of
$N$. We now choose $z\in\R$ such $z<-\rho_U/\alpha$
and define $s_4$ as the average
$$s_4(x):=\int_{P_{U}} s_3^z(y.x) dy$$
(see Definition \ref{weih1422} for $P_U$ and Theorem \ref{ausl}, (2)  for the action $(y,x)\mapsto y.x$).
The integral converges (compare Lemma  \ref{210}) and defines a positive $P_{U}$-invariant smooth function on $N$.
If $W\subset N\setminus N_V$ is compact, then there is $a^0\in A$
such that $(N_VW)^{aM_U}\cap W_1=\emptyset$ for all $a\ge a_0$.
If $a\ge a_0$ and $x\in W$, then we have
\begin{eqnarray*}
s_4(x^a)&=&\int_{P_{U}} s_3^z(y.x^a) dy\\
&=&\int_{M_U}\int_{N_V} s_3^z(vx^{au}) dv du\\
&=&a^{2\rho_U}\int_{M_U}\int_{N_V} s^z_3((vx)^{au}) dvdu\\
&=&a^{2(z\alpha+\rho_U)} s_4(x)\ .
\end{eqnarray*}
We define $s$ as the $\alpha/(z\alpha+\rho_U)$'th power of $s_4$. \hB

\subsubsection{}\label{neuj911}

We now consider the function $s$ constructed in Lemma \ref{ingf}
as a section $s\in C^\infty(\Omega_P,V(1_{\rho+\alpha}))$
by defining $s(xw):=s(x)$.

If $\lambda\in\aaaa^*$, then the bundle $V(1_\lambda)$ is a
complex line bundle which can be written as the complexification
of a trivial $G$-equivariant real line bundle $V(1_\lambda)^\R$.
It makes sense to speak of a positive
section of $V(1_\lambda)^\R$.
A section of $V(1_\lambda)$ is called positive if it is
a positive section of the real subbundle $V(1_\lambda)^\R$.
In particular, the section $s$ of $V(1_{\rho+\alpha})$
constructed above is positive.

\subsubsection{}\label{neuj910}

For each $z\in \C$ we can form $s^z$, which is a
section of $V(1_{\rho+z\alpha})$. There is a natural
identification $V(\sigma_\lambda)\otimes V(1_{\rho+z\alpha})\cong V(\sigma_{\lambda+z\alpha})$.
Multiplication by $s^z$
identifies $C^\infty(\Omega_P,V(\sigma_\lambda,\vp))$ with $C^\infty(\Omega_P,V(\sigma_{\lambda+z\alpha},\vp))$.

\begin{lem}\label{ytra}
For each $d\in\nat_0$ and $z\in\C$
multiplication by $s^z$ defines a continuous map from
$\cS_{U,k,d}(\sigma_\lambda,\vp)$ to $\cS_{U,k,d}(\sigma_{\lambda+z\alpha},\vp)$.
The family
$\{\cS_{U,k}(\sigma_\lambda,\vp)\}_{\lambda\in\aca}$ is a
trivial holomorphic bundle of Fr\'echet spaces.
\end{lem}
\proof
Since $s$ is $U$-invariant and non-vanishing multiplication by $s^z$
is an isomorphism of $C^\infty(B_{U},V_{B_{U}}(\sigma_\lambda,\vp))$
with $C^\infty(B_{U},V_{B_{U}}(\sigma_{\lambda+z\alpha},\vp))$
and  of the subspaces of sections with compact support.
It suffices to show that there
is a constant $C\in \R$ such that
for all $f\in C_c^\infty(B_{U},V_{B_{U}}(\sigma_\lambda,\vp))$
we have $\|s^z f\|_{k,d}\le C \|f\|_{k,d}$.
This follows from the Leibniz rule and the following estimate.
Let $A\in\cU(\naaa)$ be homogeneous of degree $d\alpha$. Then for any compact subset
$W\subset N\setminus N_V$ there is a constant $C\in\R$ and $a_0\in A$
such that for $x\in W$ and $a\ge a_0$
\begin{eqnarray*}
|s^z(x^aAw)|&=&|s^z((xA^{a^{-1}})^aw)|\\
&=&a^{(2z-d)\alpha}  |s^z(xAw)|\\
&\le& C a^{(2z-d)\alpha} \ .
\end{eqnarray*}

For any $\lambda_0$ we now define the trivialization
$\Phi_{\lambda_0}:\bigcup_{\lambda\in\aca} \cS_{U,k}(\sigma_{\lambda},\vp)\rightarrow  \cS_{U,k}(\sigma_{\lambda_0},\vp) \times \aca$
such that the restriction  of $\Phi_{\lambda_0}$ to the fibre $\cS_{U,k}(\sigma_{\lambda},\vp)$  
is multiplication by $s^{(\lambda_0-\lambda)/\alpha}$.
The transition map $\Phi_{\lambda_0}\circ \Phi_{\lambda_1}^{-1}$
is multiplication by $s^{(\lambda_0-\lambda_1)/\alpha}$
and thus independent of $\lambda$. In particular, it
is a holomorphic family of continuous maps.
\hB

\subsubsection{}

The Schwartz space is the space of smooth rapidly decaying sections on $B_U$. 
\begin{ddd}\label{weih245}
We define the Schwartz space
$S_{U}(\sigma_\lambda,\vp)$ as the intersection of the spaces
$S_{U,k}(\sigma_\lambda,\vp)$ over all $k\in\nat_0$.
\end{ddd}

\subsubsection{}

Since the trivializations $\Phi_\lambda$ are compatible with the
inclusions  $S_{U,k^\prime}(\sigma_\lambda,\vp)\hookrightarrow
 S_{U,k}(\sigma_\lambda,\vp)$, $k^\prime\ge k$, we obtain
the following corollary. 
\begin{kor}
The family $\{S_{U}(\sigma_\lambda,\vp)\}_{\lambda\in \aca}$
is a trivial holomorphic vector bundle of Fr\'echet spaces.
\end{kor}

\subsubsection{}

Note that the inclusions $S_{U,k+1}(\sigma_\lambda,\vp)\hookrightarrow S_{U,k}(\sigma_\lambda,\vp)$ are compact. The proof of this fact is a simple application of the Lemma of Arzela-Ascoli.
The compactness of these inclusions imply the following fact.
\begin{kor}
The Schwartz space
$S_{U}(\sigma_\lambda,\vp)$ is a Montel space. In particular, it is reflexive.
\end{kor}

\subsubsection{}

If the cusp associated to $U\subset P$
has full rank, then we have for all $k\in\nat_0$
$$S_{U,k}(\sigma_\lambda,\vp)=\cS_{U}(\sigma_\lambda,\vp)=C^\infty(B_{U},V_{B_{U}}(\sigma_\lambda,\vp))\ .$$

\subsubsection{}

Next we show that the Schwartz spaces are compatible
with twisting.
Let $\sigma$ be a Weyl invariant representation of $M$ and 
$(\pi_{\sigma,\mu},V_{\pi_{\sigma,\mu}})$ be a finite-dimensional
representation of $G$ as in \ref{weih144}. See \ref{weih146} for the definition of $Z(\lambda)$ and $\Pi(\lambda)$.
\begin{lem}\label{compat1}
\mbox{}\\
\begin{enumerate}
\item The Schwartz space
$\cS_{U}(\sigma_\lambda,\vp)$ coincides with the
closed subspace of ${}^U C^\infty(\partial X,V(\sigma_\lambda,\vp))$
of sections which vanish at $\infty_P$ of infinite order.
\item
We have holomorphic families of continuous maps
\begin{eqnarray*}
\{i_{\sigma,\mu}^U\}&:&S_{U}(\sigma_\lambda,\vp)\rightarrow
S_U(1_{\lambda+\mu},\pi_{\sigma,\mu}\otimes \vp)\\
\{p_{\sigma,\mu}^U\}&:&S_U(1_{\lambda-\mu},\pi_{\sigma,\mu}\otimes \vp)\rightarrow
\cS_{U}(\sigma_\lambda,\vp)\ .
\end{eqnarray*}

\item
$S_U(1_{\lambda+\mu},\pi_{\sigma,\mu}\otimes \vp)$
is a $Z(\gaaa)$-module.
\item
If $\lambda\not\in I_\aaaa$, then
$\{i_{\sigma,\mu}^U\}$ maps $\cS_{U}(\sigma_\lambda,\vp)$
isomorphically onto $$\ker\{Z(\lambda):S_U(1_{\lambda+\mu},\pi_{\sigma,\mu}\otimes \vp)\rightarrow S_U(1_{\lambda+\mu},\pi_{\sigma,\mu}\otimes \vp)\}\ ,$$
and the restriction of 
$\{p_{\sigma,\mu}^U\}$ to  $$\ker\{\tilde Z(\lambda):S_U(1_{\lambda-\mu},\pi_{\sigma,\mu}\otimes \vp)\rightarrow S_U(1_{\lambda-\mu},\pi_{\sigma,\mu}\otimes \vp)\}$$
is an isomorphism onto $\cS_{U}(\sigma_\lambda,\vp)$
(here $\tilde Z(\lambda)$ is the adjoint of
$\pi^{1_{-\lambda+\mu}\otimes \tilde\pi_{\sigma,\mu},\id_{\vp}}(\Pi(\lambda))$).
\item
The composition
$$\{j_{\sigma,\mu}^U\}:S_U(1_{\lambda+\mu},\pi_{\sigma,\mu}\otimes \vp)
\stackrel{1-Z(\lambda)}{\rightarrow} \ker(Z(\lambda))\stackrel{(\{i_{\sigma,\mu}^U\})^{-1}}{\rightarrow}
\cS_{U}(\sigma_\lambda,\vp)$$
which is intitially defined for $\lambda\not\in I_\aaaa$
extends to a meromorphic family of continuous maps.
Similarly, the composition
$$\{q_{\sigma,\mu}^U\}:\cS_{U}(\sigma_\lambda,\vp)\stackrel{(\{p_{\sigma,\mu}^U\})^{-1}}{\rightarrow} \ker(\tilde Z(\lambda)) \rightarrow S_U(1_{\lambda-\mu},\pi_{\sigma,\mu}\otimes \vp)$$
extends to a  meromorphic family of continuous maps.
\end{enumerate}
\end{lem}
\proof 
We leave the proof of 1. at this place to the interested reader.
Alternatively the assertion can be considered as an immediate consequence
of the material of Subsection \ref{trivas}.
\subsubsection{} 

2. follows from 1. and the fact that $\{i_{\sigma,\mu}^U\}$
and $\{p_{\sigma,\mu}^U\}$ are induced by an inclusion
$$i_{\sigma,\mu}:V(\sigma_\lambda,\vp)\rightarrow V(1_{\lambda+\mu},\pi_{\sigma,\mu}\otimes \vp)$$
and a projection $$p_{\sigma,\mu}:V(1_{\lambda-\mu},\pi_{\sigma,\mu}\otimes \vp)\rightarrow V(\sigma_\lambda,\vp)$$
of bundles. Another argument for the assertion about
$\{i_{\sigma,\mu}^U\}$ would be to check that
$\{i_{\sigma,\mu}^U\}$ maps $S_{U,k,d}(\sigma_\lambda,\vp)$ continuously to
$S_{U,k^\prime,d}(1_{\lambda+\mu},\pi_{\sigma,\mu}\otimes \vp)$,
$k^\prime=k+\mu/\alpha$,
for all $k,d\in\nat_0$ by comparing the norms explicitly.
A similar argument works for $\{p_{\sigma,\mu}^U\}$ as well.

\subsubsection{}
In order to see 3. note that the action of $cZ(\gaaa)$ is implemented
by differential operators. Using the fact that the action
$\pi^{1_{\lambda+\mu}\otimes \pi_{\sigma,\mu},\id_\vp}$
of $\cZ(\gaaa)$ commutes with the action
$\pi^{1_{\lambda+\mu},\pi_{\sigma,\mu}\otimes \vp}$
of $P_UA$ which is used in order to characterize
the growth of elements of $S_{U,k}(1_{\lambda+\mu},\pi_{\sigma,\mu}\otimes \vp)$
one could alternatively check that this space is a $\cZ(\gaaa)$-module.

\subsubsection{}

We now prove the first
assertion of 4. and leave the second to the reader since the argument
is similar. The representation  $V_{1_{\lambda+\mu}}\otimes V_{\pi_{\sigma,\mu}}$ of $P$ fits
into an exact sequence
$$0\rightarrow V_{\sigma_\lambda} \stackrel{i_{\sigma,\mu}}{\rightarrow} V_{1_{\lambda+\mu}}\otimes V_{\pi_{\sigma,\mu}}
\rightarrow W\rightarrow 0\ ,$$
which induces a corresponding exact sequence of bundles
\begin{equation}
\label{hj112}0\rightarrow V(\sigma_\lambda,\vp)\stackrel{i_{\sigma,\mu}}{\rightarrow} V(1_{\lambda+\mu},\pi_{\sigma,\mu}\otimes\vp) 
\rightarrow  V(w ,\vp) \rightarrow 0  \ ,\end{equation}
where $w $ denotes the representation of $P$ on $W $. Using 1. we obtain an exact
sequence of sections which vanish of infinite order at $\infty_P$.
Going over to $U$-invariant sections we obtain  
the exact sequence
$$0\rightarrow S_U(\sigma_\lambda,\vp)\stackrel{\{i_{\sigma,\mu}^U\}}{\rightarrow} S_U(1_{\lambda+\mu},\pi_{\sigma,\mu}\otimes\vp)
\stackrel{p}{\rightarrow} S_U(w ,\vp)  \ .$$
The action of $\cZ(\gaaa)$  commutes with $p$.
The operator $\Pi(\lambda)$ is designed such that 
$\pi^{w ,\vp}(\Pi(\lambda))=1$
and $\pi^{\sigma_\lambda,\vp}(\Pi(\lambda))=0$ (note that $\lambda\not\in I_\aaaa$).
Since $\{i_{\sigma,\mu}^U\}$ is a morphism of $Z(\gaaa)$-modules
we see that $\{i_{\sigma,\mu}^U\}$ maps to $\ker(Z(\lambda))$.
In order to show that it is onto take $f\in \ker(Z(\lambda))$.
Then $$p(f)=p((1-Z(\lambda))f)=(1-\pi^{w,\vp}(\Pi(\lambda)))p(f)=0\ ,$$
and $f$ is in the range of $\{i_{\sigma,\mu}^U\}$.

\subsubsection{}

Finally we show the first assertion of 5. and leave the second to the reader
since the argument is similar, again.
We choose a holomorphic family of splits $j :V(1_{\lambda+\mu},\pi_{\sigma,\mu}\otimes\vp) \rightarrow 
V(\sigma_\lambda,\vp)$ of the sequence of bundles (\ref{hj112}) (not necessarily $U$-invariant). The split 
$j $ induces a  holomorphic family of maps $J :C^\infty(\partial X,V(1_{\lambda+\mu},\pi_{\sigma,\mu}\otimes\vp))\rightarrow
C^\infty(\partial X,V(\sigma_\lambda,\vp))$.
The map $J$ is compatible with the subspaces of sections
vanishing of infinite order at $\infty_P$. We define
the meromorphic family of continuous maps
$\tilde j_{\sigma,\mu}^U$ as the restricition of $J \circ (1-Z(\lambda))$
to $S_U(1_{\lambda+\mu},\pi_{\sigma,\mu}\otimes \vp)$.
Since $$\tilde j_{\sigma,\mu}^U\circ i_{\sigma,\mu}^U = J \circ (1-Z(\lambda))\circ i_{\sigma,\mu}^U=\id$$
and the restriction of $\tilde j_{\sigma,\mu}^U$
to $\ker(1-Z(\lambda))$ vanishes we see by 4. that
$\tilde j_{\sigma,\mu}^U$ maps to  $S_U(\sigma_\lambda,\vp)$,
and that it coincides with $\{j_{\sigma,\mu}^U\}$ for $\lambda\not\in I_\aaaa$.
\hB

\subsection{Asymptotics for the trivial group}\label{trivas}

\subsubsection{}\label{weih600}
We consider the decomposition $G\setminus wP=\bar{N}MAN$,
$g=\bar{n}(g)m(g)a(g)n(g)$. There is a unique diffeomorphism
$F:N\setminus \{1\} \rightarrow \bar{N}\setminus \{1\}$
such that $nw\in F(n)MAN$. Indeed, $F(n)=\bar{n}(nw)$. One can check that
\begin{eqnarray}
 a(n^{ma}w)&=&a^2 a(nw)\label{weih141}\\
 m(n^aw)&=&m(nw)\nonumber\\
F(n^{ma})&=& F(n)^{ma}\nonumber\ .
\end{eqnarray}

\subsubsection{}\label{weih213}

Let $(\vp,V_\vp)$ be an admissible twist for $U$.
For simplicity we normalize the restriction of $\vp$
to $A$ such that
the lowest $A$-weights of all its irreducible coponents are zero.
Let $l_\vp\in\aaaa^*$ be the highest weight
of $\vp$.

\subsubsection{}

In the present subsection we describe the space
$C^\infty(\partial X,V(\sigma_\lambda,\vp))=:B_{\{1\}}(\sigma_\lambda,\vp)$
as an extension
$$0\rightarrow \cS_{\{1\}}(\sigma_\lambda,\vp)\rightarrow B_{\{1\}}(\sigma_\lambda,\vp)
\rightarrow \cR_{\{1\}}(\sigma_\lambda,\vp)\rightarrow 0\ .$$
The Schwartz space $\cS_{\{1\}}(\sigma_\lambda,\vp)$
coincides with the subspace of $C^\infty(\partial X,V(\sigma_\lambda,\vp))$ of all
section that vanish at $\infty_P$ of infinite order.
Note that $\bar N$ is diffeomorphic to $\partial X\setminus 0_P$
by $\bar n\mapsto \bar n \infty_P$
(recall that $0_P=w\infty_P$).
Using a theorem
of E. Borel to the effect that each formal power series
can be realized as a Taylor series of a smooth function we obtain
an exact sequence
\begin{equation}\label{borell}0\rightarrow S_{\{1\}}(\sigma_\lambda,\vp)\rightarrow C^\infty(\partial X,V(\sigma_\lambda,\vp))
\stackrel{TS}{\rightarrow} \Hom(\cU(\bar\naaa),V_{\sigma_\lambda}\otimes V_\vp)  \rightarrow 0\ ,\end{equation}
where $TS$ is given by $TS(f)(A):=f(Ae)$.

\subsubsection{}\label{neuj100}

The space $\Hom(\cU(\bar\naaa),V_{\sigma_\lambda}\otimes V_\vp)$ admits
an action of $AM_U$ by $(u.f)(A):=(\sigma_\lambda(u)\otimes \vp(u))f(A^{u^{-1}})$.
With respect to the action of $A$ it
can be decomposed as
$$\Hom(\cU(\bar\naaa),V_{\sigma_\lambda}\otimes V_\vp)=\prod_{n\in\nat_0} \Hom(\cU(\bar\naaa),V_{\sigma_\lambda}\otimes V_\vp)^n\ ,$$
where $A$ acts on the $M_U$-module $\Hom(\cU(\bar\naaa),V_{\sigma_\lambda}\otimes V_\vp)^n$
with weight $\rho-\lambda+\alpha n$. 

\subsubsection{}

By definition a function $p:\bar N\rightarrow V_{\sigma_\lambda}\otimes V_\vp$ is a polynomial
iff $p\circ \exp:\bar\naaa\rightarrow V_{\sigma_\lambda}\otimes V_\vp$ is a polynomial.
Such a function
$p$ is called homogeneous of degree $n$, iff
$$(\sigma_{\lambda}(a)\otimes \vp(a))p(\bar n^{a^{-1}})= a^{\rho-\lambda+n\alpha} p(\bar n)\ .$$
\subsubsection{}\label{neuj101}

The space $\Hom(\cU(\bar\naaa),V_{\sigma_\lambda}\otimes V_\vp)^n$ can be identified
$AM_U$-equivariantly with the space of polynomials $\Pol(\bar N,V_{\sigma_\lambda}\otimes V_\vp)^n$
on $\bar N$ with values in $V_{\sigma_\lambda}\otimes V_\vp$ being homogeneous of degree $n$.
Here a polynomial $p\in \Pol(\bar N,V_{\sigma_\lambda}\otimes V_\vp)^n$
corresponds to the map $\cU(\bar\naaa)\ni A\mapsto
p(Ae)\in V_{\sigma_\lambda}\otimes V_\vp$.

\subsubsection{}\label{weih503}

If a homogeneous polynomial $p\in \Pol(\bar N,V_{\sigma_\lambda}\otimes V_\vp)^n$
is considered as a section $f_p$ of $V(\sigma_\lambda,\vp)$ defined over
$\partial X\setminus \{0_P\}$, then we have for $x\in N$
$$f_p(xw)=f_p(F(x)m(xw)a(xw)n(xw))=\sigma(m(xw))^{-1} a(xw)^{\lambda-\rho} f_p(F(x))\ .$$
In particular, for $x\in N$, $x\not=1$, and $a\in A$ we have 
\begin{eqnarray*}
\vp(a)^{-1} f_p(x^aw)&=&\vp(a)^{-1}a(x^aw)^{\lambda-\rho} \sigma(m(x^aw))^{-1}f_p(F(x^a))\\
&=&
a^{2(\lambda-\rho)-n\alpha} a(xw)^{\lambda-\rho}\sigma(m(xw))^{-1} f_p(F(x))\ .
\end{eqnarray*}

\begin{ddd}\label{weih231}
We define $A_{\{1\}}(\sigma_\lambda,\vp)^n\subset C^\infty(\partial X\setminus
\{\infty_P,0_P\},V(\sigma_\lambda,\vp))$ as the subspace
of sections $f$ satisfying 
\begin{equation}\label{weih217}
\vp(a)^{-1}f(x^aw)=a^{2(\lambda-\rho)-n\alpha} f(xw)
\end{equation}
for all $x\in N\setminus \{1\}$, $a\in A$.
We define the subspace $\cR_{\{1\}}(\sigma_\lambda,\vp)^n\subset A_{\{1\}}(\sigma_\lambda,\vp)^n$
as the subspace spanned by the sections $f_p$ corresponding
to $p\in\Pol(\bar N,V_\vp)^n$.
We further define $\cR_{\{1\},k}(\sigma_\lambda,\vp):=\bigoplus_{n=0}^k \cR_{\{1\}}(\sigma_\lambda,\vp)^n$ and $\cR_{\{1\}}(\sigma_\lambda,\vp):=\prod_{n=0}^\infty\cR_{\{1\}}(\sigma_\lambda,\vp)^n$.
\end{ddd}

\subsubsection{}

The spaces $A_{\{1\}}(\sigma_\lambda,\vp)^n$ for the various $\lambda$
can be identified using multiplication by suitable powers of
the function $xw\mapsto a(xw)^\alpha$ and therefore form a trivial
vector bundle of Fr\'echet spaces over $\aca$. The spaces
$\cR_{\{1\}}(\sigma_\lambda,\vp)^n$ form  trivial finite-dimensional
subbundles trivialized by the identification with the trivial bundles
$\aca\times  \Pol(\bar N,V_{\sigma} \otimes V_\vp)^n\rightarrow \aca$.
Thus $\cR_{\{1\},k}(\sigma_\lambda,\vp)$ has the structure of a trivial
holomorphic vector bundle.

\subsubsection{}\label{neuj202}

We fix a smooth cut-off function $\chi\in C^\infty(A)$ such that
$\chi(a)=0$ in a neighbourhood of $A_-$ and $\chi(a)=1$ if $a^\alpha>2$.
Multiplication by the function 
$$\Omega_P\ni x\,0_P\mapsto \chi(a(xw))$$ 
defines inclusions
$L: \cR_{\{1\}}(\sigma_\lambda,\vp)^n \rightarrow C^\infty(B_{\{1\}},V_{\{1\}}(\sigma_\lambda,\vp))$
for each $n$, and summing these maps up we obtain inclusions
$$L: \cR_{\{1\},k}(\sigma_\lambda,\vp) \rightarrow C^\infty(B_{\{1\}},V_{B_{\{1\}}}(\sigma_\lambda,\vp))\ .$$
In the present paper the symbol $L$ is used for various maps of this kind. It will always be clear from the context which version of $L$ is meant. 

\subsubsection{}

If $f\in \cS_{\{1\},k}(\sigma_\lambda,\vp)$,
then
$$\lim_{a\to\infty} a^{k\alpha-2(\lambda-\rho)}\vp(a)^{-1} f(x^aw) =0$$
uniformly for $x$ in compact subsets of $N\setminus \{1\}$.
Thus $L(\cR_{\{1\},k}(\sigma_\lambda,\vp))$ intersects $\cS_{\{1\},k}(\sigma_\lambda,\vp)$ trivially. We define
$$B_{\{1\},k}(\sigma_\lambda,\vp):=\cS_{\{1\},k}(\sigma_\lambda,\vp)\oplus L(\cR_{\{1\},k}(\sigma_\lambda,\vp))\ .$$
This space fits into the exact sequence
\begin{equation}\label{weih210}
0\rightarrow \cS_{\{1\},k}(\sigma_\lambda,\vp) \rightarrow  B_{\{1\},k}(\sigma_\lambda,\vp)
\stackrel{AS}{\rightarrow}\cR_{\{1\},k}(\sigma_\lambda,\vp)\rightarrow 0\ ,
\end{equation}
where $AS$ takes the finite asymptotic expansion. To $AS$ applies the same remark as as above for $L$. This symbol appears in various versions which will be denoted by the same symbol, and it will be clear from the context which version is meant.

\subsubsection{}

The sequence (\ref{weih210}) is split by $L$
and defines $B_{\{1\},k}(\sigma_\lambda,\vp)$ as a topological vector space
and the family of spaces $\{B_{\{1\},k}(\sigma_\lambda,\vp)\}_{\lambda\in\aca}$
as a trivial bundle of Fr\'echet spaces.
Moreover, we have continuous and compact inclusions
$B_{\{1\},k+1}(\sigma_\lambda,\vp)\hookrightarrow B_{\{1\},k}(\sigma_\lambda,\vp)$.
The intersection of these spaces over all $k\in\nat_0$ is the Fr\'echet
and Montel space $B_{\{1\}}(\sigma_\lambda,\vp)$, which  concides
with $C^\infty(\partial X,V(\sigma_\lambda,\vp))$ because of
exactness of (\ref{borell}).
Note that the
sequence
$$0\rightarrow \cS_{\{1\}}(\sigma_\lambda,\vp)\rightarrow B_{\{1\}}(\sigma_\lambda,\vp)
\rightarrow \cR_{\{1\}}(\sigma_\lambda,\vp)\rightarrow 0$$
does not admit any continuous split.

\subsubsection{}

Though this is clearly possible we will not attempt to
trivialize the family of spaces $B_{\{1\}}(\sigma_\lambda,\vp)$.
But we keep in mind that $B_{\{1\}}(\sigma_\lambda,\vp)$ is an intersection
of trivial holomorphic bundles (see \ref{weih211}).
Below we will obtain a similar description of the spaces
$B_{U,k}(\sigma_\lambda,\vp)$.

\subsection{The spaces $B_{U,k}(\sigma_\lambda,\vp)$}\label{asz}

\subsubsection{}

In the present subsection $\sigma$ is any finite-dimensional representation of $M$, and $\vp$ denotes a twist. Below we define the  spaces $B_{U}(\sigma_\lambda,\vp)$.
In particular, we construct the spaces 
$\cR_{U}(\sigma_\lambda,\vp)$ which describe the asymptotic behaviour
of elements of $B_{U}(\sigma_\lambda,\vp)$.
For simplicity we assume that the twist $\vp$ is normalized
as in \ref{weih213}.

\subsubsection{}\label{weih244}

We assume that $U\subset P$ defines a cusp of smaller rank.
The case of cusps of full rank will be discussed in  
Subsection \ref{maxrank}.
We want to define
a map $\pi^{P_{U}}_*:A_{\{1\}}(\sigma_\lambda,\vp)^n\rightarrow
C^\infty(\Omega_P\setminus N_V \{0_P\},V(\sigma_\lambda,\vp))$
by
$$\pi^{P_{U}}_*(f):=\int_{P_U}
\pi^{\sigma_\lambda,\vp}(u) fdu  \ .$$
If $x\not\in N_V$, then the integrand  $(\pi^{\sigma_\lambda,\vp}(u) f)(xw)$ is well-defined for all $u\in P_{U}$.
We will see below that this integral converges for
$\Ree(\lambda)<\rho^U-l_\vp/2$.
Let $S(N_V):=\{x\in N_V|a(xw)=1\}$ (see \ref{weih213} for a definition of $l_\vp$ and Definition \ref{weih215} for $\rho^U$).

\subsubsection{}

Employing the assumption that $N_V$
is invariant under conjugation by $A$ (regularity of the cusp, Definition \ref{t799})
we can consider polar coordinates $S(N_V)\times A\ni (\xi,a)\mapsto \xi^a\in N_V$
of $N_V\setminus \{1\}$.
\begin{lem}\label{210}
There is a measure $d\xi$ on $S(N_V)$ such that the Haar measure $dx$
of $N_V$ is given by $ a^{2\rho_U} dad\xi $.
\end{lem}
\proof 
There is a family of measures $d\xi(a)$ on $S(N_V)$
such that for any $f\in C_c(N_V)$ we have
$$\int_{N_V} f(x) dx =\int_A \int_{S(N_V)} f(\xi^a) d\xi(a) a^{2\rho_U} da\ .$$
We compute for any $b\in A$
\begin{eqnarray*}
\int_{N_V} f(x^{b^{-1}}) dx &=& b^{2\rho_U}\int_{N_V} f(x) dx\\
&=& b^{2\rho_U} \int_A \int_{S(N_V)} f(\xi^a) d\xi(a) a^{2\rho_U} da\ ,\\
\int_{N_V} f(x^{b^{-1}}) dx &=&\int_A \int_{S(N_V)} f(\xi^{ab^{-1}}) d\xi(a) a^{2\rho_U} da\\
&=&b^{2\rho_U} \int_A \int_{S(N_V)} f(\xi^a) d\xi(ab) a^{2\rho_U} da\ .
\end{eqnarray*}
We conclude that $d\xi(ab)=d\xi(a)$ for all $b\in A$ and hence $d\xi(a)=d\xi(1)=:d\xi$. \hB

\subsubsection{}\label{weih221}

Recall the the Definition \ref{weih1422} of $P_U=N_V M_U$.
We introduce the abbreviation $$f_0:= \int_{M_{U}} \pi^{\sigma_\lambda,\vp}(t) f dt\ .$$
Then we can write for $x\not\in N_V$ and using the homogeneity (\ref{weih217})
\begin{eqnarray}
\pi^{P_{U}}_*(f)(xw)&=&\int_A \int_{S(N_V)}    \vp(\xi^a)^{-1} f_0( \xi^a xw) a^{2\rho_U}    d\xi da\label{hj11}\\&=&\int_A \int_{S(N_V)}   \vp(a) \vp(\xi)^{-1} \vp(a)^{-1}f_0( a\xi a^{-1} xw) a^{2\rho_U}    d\xi da\nonumber\\
&=&\int_A \int_{S(N_V)}   \vp(a) \vp(\xi)^{-1} \vp(a)^{-1}f_0( a\xi x^{a^{-1}}a^{-1} w) a^{2\rho_U}    d\xi da\nonumber\\
&=& \int_A\int_{S(N_V)}   \vp(a) \vp(\xi)^{-1} f_0(\xi x^{a^{-1}}w) a^{2(\lambda-\rho^U)-n\alpha}  d\xi  da\ .\label{hj12}
\end{eqnarray}
We decompose the outer integral into the integrals over $A_+$ and $A_-$
and obtain $\pi^{P_{U} }_+(f)$, $\pi^{P_{U} }_- (f)$ such that $\pi^{P_{U}}_*(f)=\pi^{P_{U} }_+(f) +\pi^{P_{U} }_-(f) $.
It is clear from (\ref{hj11}) that $\pi^{P_{U}  }_-(f)$ converges for all $\lambda\in\aca$.
We conclude from (\ref{hj12}) that $\pi^{P_{U}}_+(f)$ converges for
$\Ree(\lambda)<\rho^U+(n\alpha-l_\vp)/2$.

In this domain of convergence we compute for $b\in A$
$$\vp(b)^{-1}\pi^{P_{U} }_*(f)(x^bw)= b^{2(\lambda-\rho^U)-n\alpha}\pi^{P_{U} }_*(f)(xw)\ .$$

\subsubsection{}

Motivated by this calculation we make the following definition.
 
\begin{ddd}
We define $A_{P_U}(\sigma_\lambda,\vp)^n$ to be the space
of all $P_U$-invariant $f\in C^\infty(\Omega_P\setminus N_V 0_P,V(\sigma_\lambda,\vp))$
satisfying $$\vp(a)^{-1}f(x^aw)=a^{2(\lambda-\rho^U)-n\alpha}
f(xw)$$ for all $a\in A$, $x\in N\setminus N_V$.
\end{ddd}

\subsubsection{}

\begin{lem}
The family of spaces $\{A_{P_U}(\sigma_\lambda,\vp)^n\}_{\lambda\in\aca}$ forms a trivial holomorphic family
of Fr\'echet and Montel spaces.
\end{lem}
\proof
For some $z\in\R$ with $z<-\rho_U/\alpha$
we define the function
$$s_1(x):=\int_{P_{U}} a(y.xw)^{z\alpha} dy\ .$$
We further define  $s$ as the $\alpha/(z\alpha+\rho_U)$'th power of $s_1$.
Now we consider the function $s$ 
as a section $s\in C^\infty(\Omega_P\setminus N_V0_P,V(1_{\rho+\alpha}))$
defining $s(xw):= s(x)$.
Multiplication by $s^\mu$ defines 
an continuous isomorphism of $A_{P_U}(\sigma_\lambda,\vp)^n$
with $A_{P_U}(\sigma_{\lambda+\mu},\vp)^n$. We employ these isomorphisms in order
to define a holomorphic trivialization of the bundle.
\hB

\subsubsection{}

For $\Ree(\lambda)<\rho^U+(n\alpha-l_\vp)/2$ we have defined above a holomorphic family of continuous maps
$$\pi^{P_{U} }_*:A_{\{1\}}(\sigma_\lambda,\vp)^n\rightarrow A_{P_{U}}(\sigma_\lambda,\vp)^n\ .$$
\begin{lem}\label{aaww}
The family $\pi^{P_{U} }_*$ extends meromorphically to all of $\aca$ with at most first order poles in the set $\rho^U+\frac{n\alpha-l_\vp}{2} +\frac{1}{2}\nat_0$. The residues are finite-dimensional.
\end{lem}
\proof 
Recall the decomposition $\pi^{P_{U} }_*=\pi^{P_{U} }_++\pi^{P_{U} }_-$ introduced in \ref{weih221}. We have already seen that $\pi^{P_{U} }_-$ has a holomorphic continuation.
It suffices to show that $\pi^{P_{U} }_+$ has a meromorphic continuation.

Consider $f\in A_{\{1\}}(\sigma_\lambda,\vp)^n$.
Let $$F(x):=\int_{S(N_V)}  \vp(\xi)^{-1} f_0(\xi xw) d\xi$$
(see \ref{weih221} for the definition of $f_0$).
This function is smooth in a small neighbourhood
of $1\in N$ and on $N\setminus N_V$. Using the Taylor formula
for each $r\in\nat_0$ we can write
$F(x^a)=\sum_{q=0}^r F_q(x) a^{q\alpha}+ a^{(r+1)\alpha}R_r(x,a)$,
where $F_q$ is a homogeneous polynomial of degree $q\alpha$ on $N$
and $R_r(x,a)$ is uniformly bounded as $a\to 0$.
We write
$\pi^{P_{U} }_+(f)(xw):=\int_{A_+}  \vp(a) F(x^{a^{-1}}) a^{2(\lambda-\rho^U)-n\alpha}da$
and insert the expansion for $F$ in order to obtain
$$\pi^{P_{U} }_+(f)(xw):=J_r^1(f)(xw)+J_r^2(f)(xw)\ ,$$ where
\begin{eqnarray*}
J_r^1(f)(xw)&:=& \sum_{q=0}^r  \int_{A_+} \vp(a)  a^{2(\lambda-\rho^U)-(q+n)\alpha}da F_q(x)\\
 J_r^2(f)(xw)&:= &\int_{A_+}  \vp(a)R_r(x,a^{-1})   a^{2(\lambda-\rho^U)-(r+1+n)\alpha}da\ .\end{eqnarray*}
The integral $J_r^2$ converges for $\Ree(\lambda)<\rho^U+ (n+r+1)\alpha/2-l_\vp/2$.
In order to evaluate $J_r^1$ we introduce  the operator $B:=(d/da)_{|a=1}\vp(a)\in \End(V_\vp)$.
In order to define this derivative  we embed $A$ into the multiplicative group of $\R$.
Since $\vp$ is assumed to be algebraic as a representation of $A$ the eigenvalues of $B$ are integral. 
We have
$$J_r^1(f)(xw) = - \sum_{q=0}^r   (B+2(\lambda-\rho^U)-(n+q)\alpha)^{-1}  F_q(x)\ .$$
In particular we see that $J_r^1$ has a meromorphic continuation
to all of $\aca$.
Since we can choose $r$ arbitrarily large we obtain a meromorphic continuation
of $\pi^{P_{U} }_+$ to all of $\aca$. \hB

\subsubsection{}

Recall the Definition \ref{weih231} of the subspace of polynomials
$\cR_{\{1\}}(\sigma_\lambda,\vp)^n\subset A_{1}(\sigma_\lambda,\vp)^n$.
We consider the following diagram
$$\xymatrix{\cR_{\{1\}}(\sigma_\lambda,\vp)^n\ar[d]\ar[r]&A_{1}(\sigma_\lambda,\vp)^n\ar[d]^{\pi^{P_{U} }_*}\\?\ar[r]&A_{P_{U}}(\sigma_\lambda,\vp)^n}\ .$$
We can complete the diagram using the concept of an image bundle Definition \ref{imagebundle}.
\begin{ddd}\label{weih254}
We define 
$\cR_{U}(\sigma_\lambda,\vp)^n \subset A_{P_{U}}(\sigma_\lambda,\vp)^n$
as the image bundle 
of the restriction of $\pi^{P_{U} }_*$ to $\cR_{\{1\}}(\sigma_\lambda,\vp)^n$.
We define the meromorphic family of maps
$$[\pi^{U}_*]:\cR_{\{1\}}(\sigma_\lambda,\vp)^n\rightarrow \cR_{U}(\sigma_\lambda,\vp)^n$$
to be  induced by $\pi^{P_{U}}_*$. Furthermore, we define
$\cR_{U,k}(\sigma_\lambda,\vp):=\bigoplus_{n=0}^k  \cR_{U}(\sigma_\lambda,\vp)^n$
and let $$[\pi^{U}_*]:\cR_{\{1\},k}(\sigma_\lambda,\vp)\rightarrow \cR_{U,k}(\sigma_\lambda,\vp)$$
denote the corresponding meromorphic family of maps.
\end{ddd}
Similar to the usage of $AS$ and $L$
the symbol $[\pi^{\Gamma}_*]$ denotes various versions of the push-down
on the level of asymptotic terms.  It will be clear from the context wich version is meant. 

\subsubsection{}\label{weih500}

Using the second assertion of Lemma \ref{bunle}
we choose once and for all meromorphic families of right-inverses $$[Q]: \cR_{U}(\sigma_\lambda,\vp)^n \rightarrow \cR_{\{1\}}(\sigma_\lambda,\vp)^n $$
of $[\pi^{U}_*]$.

\subsubsection{}\label{weih253}

Let $\chi\in C^\infty(B_U)$ be a cut-off function such that
$1-\chi$ has compact support, and which vanishes in a neighbourhood
of $U\backslash (P_U 0_P)$. We can assume that $\chi$
is $P_U$-invariant (otherwise we replace it by
the average $xw\mapsto  \int_{U\backslash P_U} \chi(u^{-1}xw) du$).
Multiplication by $\chi$ defines
an inclusion $L:\cR_{U,k}(\sigma_\lambda,\vp)\rightarrow
C^\infty(B_{U},V_{B_{U}}(\sigma_\lambda,\vp))$. 
As in Subsection \ref{trivas} we see that
$L(\cR_{U,k}(\sigma_\lambda,\vp))\cap \cS_{U,k}(\sigma_\lambda,\vp)=\{0\}$.

\subsubsection{}\label{neuj915}

\begin{ddd}\label{neuj104}
We define the Fr\'echet space
$$B_{U,k}(\sigma_\lambda,\vp):=  \cS_{U,k}(\sigma_\lambda,\vp)\oplus L(\cR_{U,k}(\sigma_\lambda,\vp))\ .$$ Furthermore we
  define the Fr\'echet and Montel space $B_{U}(\sigma_\lambda,\vp)$ to be the intersection
of the spaces $B_{U,k}(\sigma_\lambda,\vp)$ for all $k\in\nat$.
\end{ddd}
Note that $B_{U,k}(\sigma_\lambda,\vp)$ fits into the split exact sequence
\begin{equation}\label{sse0}0\rightarrow  \cS_{U,k}(\sigma_\lambda,\vp)\rightarrow  B_{U,k}(\sigma_\lambda,\vp)
\stackrel{AS}{\rightarrow} \cR_{U,k}(\sigma_\lambda,\vp)\rightarrow 0\ ,\end{equation}
where $AS$ takes the finite asymptotic expansion.
Since the spaces $\cS_{U,k}(\sigma_\lambda,\vp)$ and $\cR_{U,k}(\sigma_\lambda,\vp)$ form trivial holomorphic bundles we can employ the split $L$ of the sequence
in order to equip the family of spaces $\{B_{U,k}(\sigma_\lambda,\vp)\}_{\lambda\in\aca}$
with the structure of a trivial holomorphic bundle over $\aca$.
The space 
$B_{U}(\sigma_\lambda,\vp)$
fits into the exact sequence
\begin{equation}\label{sse1}0\rightarrow S_{U}(\sigma_\lambda,\vp)\rightarrow
B_{U}(\sigma_\lambda,\vp)\rightarrow \cR_{U}(\sigma_\lambda,\vp)\rightarrow 0\end{equation}
which does not admit any continuous split.
The spaces $B_{U}(\sigma_\lambda,\vp)$ form a projective limit
of locally trivial holomorphic bundles in the sense of  \ref{weih211}.

\subsubsection{}

We now show compatibility of the spaces
$B_{U}(\sigma_\lambda,\vp)$ with twisting.
We assume that $\sigma$ is Weyl invariant as in \ref{weih241}. 
Let $(\pi_{\sigma,\mu},V_{\pi_{\sigma,\mu}})$ be a finite-dimensional
representation of $G$ as in \ref{weih144}.
Recall the notation $p_{\sigma,\mu}$ and $i_{\sigma,\mu}$ from \ref{weih142}. 

\begin{lem}\label{compat2}
\begin{enumerate}
\item For each $n\in\nat_0$ we have the following commutative diagrams
\begin{eqnarray}&&
\begin{array}{ccc}
R_{\{1\}}(\sigma_\lambda,\vp)^n&\stackrel{i_{\sigma,\mu}^n}{\rightarrow}&
R_{\{1\}}(1_{\lambda+\mu},\pi_{\sigma,\mu}\otimes \vp)^m\\
\downarrow [\pi^{U}_*]&&\downarrow [\pi^{U}_*]\\
R_{\{U\}}(\sigma_\lambda,\vp)^n&\stackrel{[i_{\sigma,\mu}^{U}]}{\rightarrow}&
R_{\{U\}}(1_{\lambda+\mu},\pi_{\sigma,\mu}\otimes \vp)^m
\end{array}\label{uz76}\\&&
\begin{array}{ccc}
R_{\{1\}}(\sigma_\lambda,\vp)^n&\stackrel{p_{\sigma,\mu}^n}{\leftarrow}&
R_{\{1\}}(1_{\lambda-\mu},\pi_{\sigma,\mu}\otimes \vp)^m\\
\downarrow [\pi^{U}_*]&&\downarrow [\pi^{U}_*]\\
R_{\{U\}}(\sigma_\lambda,\vp)^n&\stackrel{[p_{\sigma,\mu}^{U}]}{\leftarrow}&
R_{\{U\}}(1_{\lambda-\mu},\pi_{\sigma,\mu}\otimes \vp)^m
\end{array}   \ ,\nonumber
\end{eqnarray}
where $m:=n+\mu/\alpha$.
\item
We have holomorphic families of continuous maps
$$i_{\sigma,\mu}^U:B_{U}(\sigma_\lambda,\vp)\rightarrow
B_U(1_{\lambda+\mu},\pi_{\sigma,\mu}\otimes \vp)$$
and
$$p_{\sigma,\mu}^U:B_U(1_{\lambda-\mu},\pi_{\sigma,\mu}\otimes \vp)\rightarrow
B_{U}(\sigma_\lambda,\vp)\ .$$
\item
$B_U(1_{\lambda+\mu},\pi_{\sigma,\mu}\otimes \vp)$
is a $Z(\gaaa)$-module.
\item
If $\lambda\not\in I_\aaaa$, then
$i_{\sigma,\mu}^U$ maps $B_{U}(\sigma_\lambda,\vp)$
isomorphically onto $$\ker\{Z(\lambda):B_U(1_{\lambda+\mu},\pi_{\sigma,\mu}\otimes \vp)\rightarrow B_U(1_{\lambda+\mu},\pi_{\sigma,\mu}\otimes \vp)\}\ ,$$
and the restriction of
$p_{\sigma,\mu}^U$ to  $$\ker\{\tilde Z(\lambda):B_U(1_{\lambda-\mu},\pi_{\sigma,\mu}\otimes \vp)\rightarrow B_U(1_{\lambda-\mu},\pi_{\sigma,\mu}\otimes \vp)\}$$
is an isomorphism onto $B_{U}(\sigma_\lambda,\vp)$ 
(here $\tilde Z(\lambda)$ is the adjoint of
$\pi^{1_{-\lambda+\mu}\otimes\tilde\pi_{\sigma,\mu}, \id_{\tilde\vp}}(\Pi(\lambda))$).
\item
The composition
$$j_{\sigma,\mu}^U:B_U(1_{\lambda+\mu},\pi_{\sigma,\mu}\otimes \vp)
\stackrel{1-Z(\lambda)}{\rightarrow} \ker_U(Z(\lambda))\stackrel{(i_{\sigma,\mu}^U)^{-1}}{\rightarrow}
B_{U}(\sigma_\lambda,\vp)$$
which is initially defined for $\lambda\not\in I_\aaaa$
extends to a meromorphic family of continuous maps.
Similarly, the composition
$$q_{\sigma,\mu}^U:B_{U}(\sigma_\lambda,\vp)\stackrel{(p_{\sigma,\mu}^U)^{-1}}{\rightarrow} \ker(\tilde Z(\lambda)) \rightarrow B_U(1_{\lambda-\mu},\pi_{\sigma,\mu}\otimes \vp)$$
extends to a meromorphic family of continuous maps.
\end{enumerate}
\end{lem}
\proof
\subsubsection{}
We consider the first diagram of 1. and leave the second to the reader
since the argument is similar.
First one checks that the inclusion
$i_{\sigma,\mu}$ induces an inclusion
$$i_{\sigma,\mu}^{P_U,n}:\cA_{P_U}(\sigma_\lambda,\vp)^n\rightarrow
\cA_{P_U}(1_{\lambda+\mu},\pi_{\sigma,\mu}\otimes \vp)^{m}$$
for all $n\in\nat_0$, $m:=n+\mu/\alpha$. Now we obtain diagram
(\ref{uz76}),
where the map $[i_{\sigma,\mu}^{U}]$ is the restriction
of $i_{\sigma,\mu}^{P_U,n}$ which is well-defined by the
naturality of the image bundle construction Lemma \ref{bunle}.

\subsubsection{}
Assertion 2. follows from Lemma \ref{compat1}, 2. and 1.
Indeed, $i^{U}_{\sigma,\mu}$ maps
$B_{U,k}(\sigma_\lambda,\vp)$ to
$$\{i^U_{\sigma,\mu}\}(\{S_{U,k}(\sigma_\lambda,\vp)\})+L([i_{\sigma,\mu}^{U}](R_{U,k}(\sigma_\lambda,\vp)))
\subset B_{U,k^\prime}(1_{\lambda+\mu},\pi_{\sigma,\mu}\otimes \vp) \ ,$$
where $k^\prime=k+\mu/\alpha$. The argument for the second
assertion of 2. is similar.

\subsubsection{}
In order to prove 3. we first show that the spaces
$R_{U}(1_{\lambda+\mu},\pi_{\sigma,\mu}\otimes \vp)^m$ are a $\cZ(\gaaa)$-modules.
Note that the representations $\pi^{1_{\lambda+\mu}\otimes\pi_{\sigma,\mu},\id_\vp}$
of $\cZ(\gaaa)$ and $\pi^{1_{\lambda+\mu},\pi_{\sigma,\mu}\otimes \vp}$
of $AP_U$ commute. Since we employ the latter representation in order to define homogeneities
we see that the spaces $A_{P_U}(1_{\lambda+\mu},\pi_{\sigma,\mu}\otimes \vp)^m$
are $\cZ(\gaaa)$-modules. The space $R_{\{1\}}(1_{\lambda+\mu},\pi_{\sigma,\mu}\otimes \vp)$
is clearly a $Z(\gaaa)$-module since it is the quotient
of $C^\infty(\partial X,V(1_{\lambda+\mu},\pi_{\sigma,\mu}\otimes \vp))$
by the submodule $S_{\{1\}}(1_{\lambda+\mu},\pi_{\sigma,\mu}\otimes \vp)$
(Lemma \ref{compat1}, 3.). Thus the subspaces
$R_{\{1\}}(1_{\lambda+\mu},\pi_{\sigma,\mu}\otimes \vp)^m$ are
$Z(\gaaa)$-modules. We obtain the $Z(\gaaa)$-module structure
on  $R_{U}(1_{\lambda+\mu},\pi_{\sigma,\mu}\otimes \vp)^m$
from Lemma \ref{bunle} since $\pi^{P_U}_*$ is $\cZ(\gaaa)$-equivariant.

Assertion 3. now follows from Lemma \ref{compat1}, 3. and the fact that  for any $A\in \cZ(\gaaa)$ the commutator
$[\pi^{1_{\lambda+\mu}\otimes\pi_{\sigma,\mu},\id_\vp}(A),\chi]$
is a differential operator with compactly supported coefficients on
$B_U$.

\subsubsection{}
We now prove the first assertion of 4. and leave the second to the reader.
The exact sequence of bundles (\ref{hj112}) induces an exact sequence
$$0\rightarrow R_U(\sigma_\lambda,\vp)^n\stackrel{[i^U_{\sigma,\mu}]}{\rightarrow}
R_U(1_{\lambda+\mu},\pi_{\sigma,\mu}\otimes \vp)^m
\stackrel{p}{\rightarrow} R_U(w,\vp)^k\ ,$$
where $m=n+\mu/\alpha$ and we employ an appropriate definition of
homogeneity for the last space. If $\lambda\not\in I_\aaaa$,
then we have $\ker(Z(\lambda):R_U(1_{\lambda+\mu},\pi_{\sigma,\mu}\otimes \vp)^m
\rightarrow R_U(1_{\lambda+\mu},\pi_{\sigma,\mu}\otimes \vp)^m)=\ker(p)$.
For these $\lambda$ the map $[i^U_{\sigma,\mu}]$ identifies
$R_U(\sigma_\lambda,\vp)^n$ with the kernel of $Z(\lambda)$.
If $f\in B_U(1_{\lambda+\mu},\pi_{\sigma,\mu}\otimes \vp)$
satisfies $Z(\lambda)(f)=0$, then $Z(\lambda)AS(f)=0$.
From what was shown above and (\ref{sse1}) it follows
that there exists $g\in B_U(\sigma_\lambda,\vp)$
such that $[i^U_{\sigma,\mu}]\circ AS(g) = AS(f)$.
Then $f-i^U_{\sigma,\mu}(g)\in S_U(1_{\lambda+\mu},\pi_{\sigma,\mu}\otimes \vp)$,
and by Lemma \ref{compat1}, 4. there exists $h\in S_U(\sigma_\lambda,\vp)$
such that $\{i^U_{\sigma,\mu}\}(h)=f-i^U_{\sigma,\mu}(g)$.
Thus $f= i^U_{\sigma,\mu}(h+g)$, and this finishes the proof of 4.

\subsubsection{}
We now prove the first assertion of 5. and leave the second to the reader.
We employ the notation introduced in the proof of Lemma \ref{compat1}, 5.
For each $n\in\nat_0$ there is $m\in\nat_0$ such that
$J:S_{U,m}(1_{\lambda+\mu},\pi_{\sigma,\mu}\otimes \vp) \rightarrow
S_{U,n}(\sigma_\lambda,\vp)$ is a   holomorphic family of maps.
 We define the  meromorphic family of continuous maps
$\tilde j^U_{\sigma,\mu}:B_{U,m}(1_{\lambda+\mu},\pi_{\sigma,\mu}\otimes \vp) \rightarrow B_{U,n}(\sigma_\lambda,\vp)$ by
$$\tilde j^U_{\sigma,\mu}(f):=
L\circ [i^U_{\sigma,\mu}]^{-1}\circ AS\circ (1-Z(\lambda))(f) +
\{j^U_{\sigma,\mu}\}(f- L\circ  AS \circ (1-Z(\lambda))(f))\ .$$
Then $\tilde j^U_{\sigma,\mu}\circ i^U_{\sigma,\mu}(f)=f$
for any $f\in B_{U}(\sigma_\lambda,\vp)$.
Since $\tilde j^U_{\sigma,\mu}$ vanishes on $\ker(1-Z(\lambda))$
we conclude that
the restriction of $\tilde j^U_{\sigma,\mu}$ to
$B_{U}(1_{\lambda+\mu},\pi_{\sigma,\mu}\otimes \vp)$
coincides with $j^U_{\sigma,\mu}$ for $\lambda\not\in I_\aaaa$.
This proves 5. \hB

\subsection{The push-down for Schwartz spaces}\label{psuscg}

\subsubsection{}

In this subsection we show that the push-down induces a map
between Schwartz spaces
$$\{\pi^U_*\}:\cS_{\{1\},k}(1_\lambda,\vp)\rightarrow  \cS_{U,k}(1_\lambda,\vp) \ .$$
Note that we only consider the spherical case $\sigma=1$.
\begin{lem}\label{schwpush}
\begin{enumerate}
\item The push-down  $\{\pi^U_*\}$ converges for
$\Ree(\lambda)<\rho^U+(k\alpha-l_\vp)/2$.
\item For $k_1<k$ it induces a holomorphic
family of maps
$\{\pi^U_*\}:\cS_{\{1\},k}(1_\lambda,\vp)\rightarrow  \cS_{U,k_1}(1_\lambda,\vp)$.
\end{enumerate}
\end{lem}
\proof
\subsubsection{}
Let $f\in C_c^\infty(B_{\{1\}},V_{B_{\{1\}}}(1_\lambda,\vp))$.
Then $\{\pi^{U}_*\}(f)=\sum_{u\in U}\pi^{1_\lambda,\vp}(u)f$
converges and defines an element of
$C_c^\infty(B_{U},V_{B_{U}}(1_\lambda,\vp))$.
\subsubsection{}
Let $h\in\nat_0$, $W\subset N\setminus N_V$ be compact, $|.|$
be a $\vp(M_{U})$-invariant norm on $V_\vp$, and
$D\in \cU(\naaa)^{\le d \alpha}$ (see \ref{weih243} for notation). Then we consider the seminorm
$$q_{W,D,h}(f):=\sup_{x\in W}\sup_{a\in A_+}  a^{h\alpha-2(\lambda-\rho^U)}|\vp(a)^{-1}(f)(x^aDw)|\ .$$
The first assertion of the Lemma will easily follow from estimates of the form
$$q_{W,D,k}(\{\pi^{U}_*\}f)<C\|f\|_{k,d}\ ,$$
and the second assertion will follow from
the convergence
of sums of the form
$$\sum_{u\in U} q_{W,D,k_1}(\pi^{1_\lambda,\vp}(u) f)<C\|f\|_{k,d} \ .$$

\subsubsection{}
If $u\not=1$, then we can write
$u^{-1}=m_u\xi_u^{b_u}$ for suitable $\xi_u\in S(N_V)$, $b_u\in  A$, and $m_u\in M_{U}$ (see \ref{weih244} for the notation $S(N_V)$).
There is a constant $C\in R$ such that for all $f\in C_c^\infty(B_{\{1\}},V_{B_{\{1\}}}(1_\lambda,\vp))$,
$x\in W$, $a\in A$, $1\not=u\in U$ with
$b_u\ge a$ we have the estimate
\begin{eqnarray*}
|\vp(a)^{-1}\pi^{1_\lambda,\vp}(u) f(x^aDw)|&=&
|\vp(a)^{-1}\vp(\xi_u^{b_u})^{-1}
f(m_u\xi_u^{b_u}x^aDw)|\\
&=&
|\vp(a^{-1}b_u)\vp(\xi_u)^{-1} \vp(b_u^{-1})f(m_u(\xi_ux^{ab_u^{-1}})^{b_u}Dw)|\\
&\le&
C\|f\|_{d,k} |\vp(a^{-1}b_u)| b_u^{2(\lambda-\rho)-k\alpha}
\end{eqnarray*}
and for $b_u\le a$
\begin{eqnarray*}
|\vp(a)^{-1}\pi^{1_\lambda,\vp}(u) f(x^aDw)|&=&
|\vp(a)^{-1} \vp(\xi_u^{b_u})^{-1} f(m_u\xi_u^{b_u}x^aDw)|\\
&=&
|\vp(\xi_u^{a^{-1}b_u})^{-1} \vp(a)^{-1}f(m_u(\xi^{b_ua^{-1}}x)^{a}Dw)|\\
&\le&
C\|f\|_{d,k} a^{2(\lambda-\rho)-k\alpha}\ .
\end{eqnarray*}

\subsubsection{}
If $\Ree(\lambda)<\rho^U+(k\alpha-l_\vp)/2$, then
summing up these estimates over $U$ and estimating the sum
by an integral over $P_{U}$ we obtain $C_1\in\R$ such that
all $x\in W$, $a\in A^+$, $f\in C_c^\infty(B_{\{1\}},V_{B_{\{1\}}}(1_\lambda,\vp))$
$$|\vp(a)^{-1}\{\pi^{U^0}_*\}(f)(x^aDw)|< C_1 a^{2(\lambda-\rho^U)-k\alpha}\ .$$
This implies $$q_{W,D,k}(\{\pi^{U}_*\}f)<C_1\|f\|_{k,d}\ .$$

\subsubsection{}
Multiplying the two inequalities above by  $a^{k_1\alpha-2(\lambda-\rho^U)}$
and taking the supremum over $x$ and $a$ we obtain
$q_{W,D,k_1}(\pi^{1_\lambda,\vp}(u) f)<b_u^{-2\rho_U-(k-k_1)\alpha}$.
Summing this over $U$ and estimating the sum over $U$
by the integral over $P_{U}$ we obtain
$$\sum_{u\in U} q_{W,D,k_1}(\pi^{1_\lambda,\vp}(u) f)<C\|f\|_{k,d}\ .$$
\hB

\subsubsection{}

Recall (Definition \ref{weih245}) that $\cS_{\{1\}}(1_\lambda,\vp)$ is the intersection of the spaces
$\cS_{\{1\},k}(1_\lambda,\vp)$ for $k\in \nat_0$. Hence Lemma \ref{schwpush}
implies the following corollary.

\begin{kor}\label{tzu}
We have a holomorphic family of maps
$$\{\pi^U_*\}:\cS_{\{1\}}(1_\lambda,\vp)\rightarrow \cS_U(1_\lambda,\vp)$$
defined on all of $\aca$.
\end{kor}

\subsubsection{}

We now turn to the problem of constructing a right-inverse of the push-down for Schwartz spaces.
Its adjoint plays an important role in the construction of the restriction map. Below we encounter the effect of a loss of regularity ($-2\rho_U$). This is one of the places where the presence of cusps makes the theory
much more complicated in comparison with the convex cocompact case.

\subsubsection{}

\begin{lem} \label{schwspl}
For all $k_1\in\nat_0$ satisfying $k_1 \alpha<k\alpha-2\rho_U$
there exists a  holomorphic family of maps
$$\{Q\}:\cS_{U,k}(1_\lambda,\vp)\rightarrow \cS_{\{1\},k_1}(1_\lambda,\vp)$$
such that $\{\pi^U_*\}\circ \{Q\}$ is the natural inclusion
$\cS_{\{1\},k}(1_\lambda,\vp)\hookrightarrow \cS_{\{1\},k_1}(1_\lambda,\vp)$.
\end{lem}
\proof
We define a holomorphic family $\{Q\}$ of right-inverses of $\{\pi^U_*\}$.
We employ a cut-off function $\chi^U\in C^\infty(\Omega_P)$
with $\sum_{u\in U} u^*\chi^U=1$, and such that
$\chi^U(xD)$ is bounded for each $D\in\cU(\naaa)$. Here we have identified
$\Omega_P$ with $N$. Then $\chi^U$ considered as a function on $N$ can be derived with respect
to the left invariant differential operator $D$. 

To see that such a function exists we equip $\Omega_P\cong N$ with a $MN$-invariant Riemannian metric
such that $U$ acts isometrically.
Note that
the Riemannian manifold $U\backslash N$ admits a lower bound of the
injectivity radius.

For $f\in \cS_{U,k}(1_\lambda,\vp)$
we define $\{Q\}(f)$ to be the lift of $f$ to $\Omega_P$ multiplied by
$\chi^U$. It is then easy to see that
$\{Q\}: \cS_{U,k}(1_\lambda,\vp) \rightarrow  \cS_{\{1\},k_1}(1_\lambda,\vp)$ for all $k_1\in\nat_0$ satisfying $k_1 \alpha<k\alpha-2\rho_U$.
Therefore we obtain a holomorphic family of maps
$\{Q\}:\cS_{U}(1_\lambda,\vp)\rightarrow \cS_{\{1\}}(1_\lambda,\vp)$.
By definition $\{\pi^U_*\}\circ \{Q\}=\id$.
\hB

\subsection{Push down for cusps of smaller rank}\label{weih261}

\subsubsection{}

We assume that the cusp associated to $U\subset P$ does not
have full rank.
In order to construct the push-down it remains to extend
$\{\pi^U_*\}$  (constructed in Lemma \ref{schwpush}) for each $k\in\nat$ from $\cS_{\{1\},k}(1_\lambda,\vp)$ to
$B_{\{1\},k}(1_\lambda,\vp)$ such that the following diagram
is commutative:
$$
\begin{array}{ccccccccc}
0&\rightarrow&\cS_{\{1\},k}(1_\lambda,\vp)&\rightarrow& B_{\{1\},k}(1_\lambda,\vp)&\stackrel{AS}{\rightarrow}& \cR_{\{1\},k}(1_\lambda,\vp)&\rightarrow &0\\
 &           &\downarrow \{\pi^{U}_*\} &           & \downarrow \pi^{U}_* &      & \downarrow [\pi^{U}_*]          &&\\
0&\rightarrow&\cS_{U,k}(1_\lambda,\vp)&\rightarrow& B_{U,k}(1_\lambda,\vp)& \stackrel{AS}{\rightarrow}&\cR_{U,k}(1_\lambda,\vp)&\rightarrow &0\ .\end{array}
$$
The main technical result in this direction is the following Proposition \ref{mainpure}.

\subsubsection{}

Note that $\cR_{\{1\}}(1_\lambda,\vp)^n$ is a finite-dimensional
representation of the torus $M_{U^0}$ (see \ref{weih251} for the explanation of $M_{U^0}$). Let $\cX(M_{U^0})$ denote the
set of characters of $M_{U^0}$.
For any $\theta\in \cX(M_{U^0})$
let $\cR_{\{1\}}(1_\lambda,\vp)^n_\theta$ be the subspace
of all $f\in \cR_{\{1\}}(1_\lambda,\vp)^n$ such that
$$\pi^{1_\lambda,\vp}(t)f=\theta(t)f\ .$$
Then we have a finite decomposition $$\cR_{\{1\}}(1_\lambda,\vp)^n=
\bigoplus_{\theta\in\cX(M_{U^0}) }
\cR_{\{1\}}(1_\lambda,\vp)^n_\theta\ .$$

\subsubsection{}

Recall the construction of $L$ from \ref{weih253}.
Furthermore see Definition \ref{weih254} for  $[\pi^{U^0}_*]$, Definition \ref{weih215} for $\rho^U$ 
and \ref{weih213} for $l_\vp$.

\begin{prop}\label{mainpure}
The composition 
$$\pi^{U^0}_*\circ L:\cR_{\{1\}}(1_\lambda,\vp)^n_\theta\rightarrow
B_{U^0}(1_\lambda,\vp)$$
converges for $\Ree(\lambda)<\rho^U+(n\alpha-l_\vp)/2$
and has a meromorphic continuation to all of $\aca$
such that $AS\circ \pi^{U^0}_*\circ L = [\pi^{U^0}_*]$.
\end{prop}
\proof
\subsubsection{}

We decompose the push-down 
$\pi^{U^0}_*$ into two intermediate maps. The first map
is the push-down
$\pi^{[U^0,U^0]}_*\circ L :
\cR_{\{1\}}(1_\lambda,\vp)^n_\theta\rightarrow  B_{[U^0,U^0]}(1_\lambda,\vp)$
with respect to the commutator group
$[U^0,U^0]$, and the second is a relative push-down
$\pi^{U^0/[U^0,U^0]}$.
Note that $M_{U^0}$ centralizes $[U^0,U^0]$ (see \ref{weih131})
and therefore acts on the sequence
$$0\rightarrow \cS_{[U^0,U^0]}(1_\lambda,\vp)\rightarrow B_{[U^0,U^0]}(1_\lambda,\vp)\rightarrow
\cR_{[U^0,U^0]}(1_\lambda,\vp)\rightarrow 0\ .$$
The representation of $M_{U^0}$ on  $\cR_{[U^0,U^0]}(1_\lambda,\vp)$
is compatible with the decomposition \linebreak[4] 
$\bigoplus_{n=0}^\infty \cR_{[U^0,U^0]}(1_\lambda,\vp)^n$.

\subsubsection{}\label{weih257}

Let $\cR_{[U^0,U^0]}(1_\lambda,\vp)^n_\theta$ be the subspace of all
$f\in \cR_{[U^0,U^0]}(1_\lambda,\vp)^n$ satisfying $$\pi^{1_\lambda,\vp}(t)f=\theta(t)f\ .$$
Then we have a further finite decomposition
$$\cR_{[U^0,U^0]}(1_\lambda,\vp)^n=
\bigoplus_{\theta\in\cX(M_{U^0}) }
\cR_{[U^0,U^0]}(1_\lambda,\vp)^n_\theta\ .$$

\subsubsection{}

We will first show
\begin{lem} \label{lpro1}
The push-down over $[U^0,U^0]$ defines a meromorphic
family of maps
$$\pi^{[U^0,U^0]}:\cR_{\{1\}}(1_\lambda,\vp)^n_\theta\rightarrow
\cS_{[U^0,U^0]}(1_\lambda,\vp) \oplus L(\cR_{[U^0,U^0]}(1_\lambda,\vp)^n_\theta)$$
such that $$AS\circ \pi^{[U^0,U^0]}\circ L=[\pi^{[U^0,U^0]}]\ .$$
\end{lem}
It is clear that the average $\{\pi_*^{U^0/[U^0,U^0]}\}: \cS_{[U^0,U^0]}(1_\lambda,\vp)\rightarrow
\cS_{U^0}(1_\lambda,\vp)$ over $U^0/[U^0,U^0]$
 converges and depends holomorphically
on $\lambda$. A simple way to see this formally is to
write
$$\{\pi_*^{U^0/[U^0,U^0]}\}=\{\pi^{U^0}_*\}\circ\{Q\}\circ  \pi_*^{[U^0,U^0]}$$
using the split $\{Q\}$ constructed in Lemma \ref{schwspl}
and Corollary \ref{tzu}.

Proposition \ref{mainpure} now  immediately follows from
\begin{lem}\label{lpro2}
The composition
$$\pi^{U^0/[U^0,U^0]}_*\circ L: \cR_{[U^0,U^0]}(1_\lambda,\vp)^n_\theta
\rightarrow  B_{U^0}(1_\lambda,\vp)$$
defines a meromorphic family of maps
such that $$AS\circ \pi^{U^0/[U^0,U^0]}_* \circ L \circ [\pi_*^{[U^0,U^0]}]=
[\pi^{U^0}_*]\ .$$
\end{lem}

\subsubsection{}
We now start with the proof of Lemma \ref{lpro1}.
Since $M_{U^0}$ is abelian the homomorphism  $m:N_V\rightarrow M_{U^0}$ vanishes
on $[N_V,N_V]$.  Therefore
$[U^0,U^0]$ is a discrete subgroup of the center of $N$.
If $f\in \cR_{\{1\}}(1_\lambda,\vp)^n$, then for all $d\in\nat_0$ we have
$\|L(f)\|_{n,d}<\infty$.
Therefore the same arguments as in the proof of Lemma \ref{schwpush} show that  the push-down
$\pi^{[U^0,U^0]}_*\circ L$ converges for $\Ree(\lambda)<\rho^{{[U^0,U^0]}}+(n\alpha-l_\vp)/2$.

\subsubsection{}\label{weih256}

We consider the abelian group $Z:=P_{{[U^0,U^0]}}=[N_V,N_V]$. 
If we are given an
unitary character $\vartheta\in \cX(Z)$,
then we consider the push-down $\pi^{Z,\vartheta}_*$,
which associates to $f\in C^\infty(\partial X,V(1_\lambda,\vp))$
the average
$$\pi^{Z,\vartheta}_*(f):=\int_{Z} \vartheta(z)^{-1} \pi^{1_\lambda,\vp}(z)(f_{|\Omega}) dz$$
provided the integral converges. 

\subsubsection{}

Observe that $[U^0,U^0]$
is a cocompact discrete subgroup of $Z$.
We normalize the Haar measure of $Z$ such that $\vol({[U^0,U^0]}\backslash Z)=1$.
Let $\cX(Z,{[U^0,U^0]})$ denote the set of all unitary characters of
$Z$ which are trivial on ${[U^0,U^0]}$. We identify $\cX(Z)$
with $\imath \zaaa^*$ and $\cX(Z,{[U^0,U^0]})$ with a lattice
in $\imath \zaaa^*$, where $\zaaa$ denotes the Lie algebra of $Z$.
Our approach is based on the
Poisson summation formula which states that
$$\pi^{[U^0,U^0]}_*(f)=\sum_{\vartheta\in\cX(Z,{[U^0,U^0]})}
\pi^{Z,\vartheta}_*(f)$$ in the domain of convergence.

\subsubsection{}

We choose an euclidean structure on $\zaaa$.
This structure induces norms on $\zaaa$ and $\zaaa^*$. Furthermore we use
the euclidean structure in order to define a Laplace operator 
$\Delta\in \cU(\zaaa)$ which is normalized such that if
$\vartheta\in\cX(Z)$, then $\vartheta(\Delta^l x)=|\vartheta|^{2l}
\vartheta(x)$.
Let $COP:\cU(\zaaa)\rightarrow \cU(\zaaa)\otimes \cU(\zaaa)$ denote
the coproduct on $\cU(\zaaa)$. There are $A_{l,j},B_{l,j},C_{l,j}\in\cU(\zaaa)$
such that
$$(COP\otimes 1)\circ COP(\Delta^l)=\sum_{j} A_{l,j}\otimes B_{l,j}\otimes C_{l,j}\ .$$
Let $D\mapsto \tilde D$ be the canonical antiautomorphism of $\cU(\zaaa)$.

\subsubsection{}

Fix $l\in\nat$ and let $f\in \cR_{\{1\}}(1_\lambda,\vp)^n_\theta$.
Recall that $L(f)(xw)=\chi(a(xw)) f(xw)$, where $\chi\in C^\infty(A)$
was some cut-off function which vanishes on a neighbourhood of $A_-$
and is equal to one near $\infty$.
By $S(Z)$ we denote the unit-sphere in $Z$.
For $\xi\in S(Z)$ and $a\in A$ we compute
\begin{eqnarray}
&&\pi^{1_\lambda,\vp}(\Delta^l\xi^a) L(f)(.w)\nonumber\\
&=&\pi^{1_\lambda,\vp}(a)\pi^{1_\lambda,\vp}((\Delta^l)^{a^{-1}}\xi) \pi^{1_\lambda,\vp}(a^{-1})\chi (a(.w))f(.w)\nonumber\\
&=&a^{\rho-\lambda-4l\alpha}\pi^{1_\lambda,\vp}(a)\pi^{1_\lambda,\vp}(\Delta^l\xi) \vp(a^{-1}) \chi(a(a.a^{-1}w)) f(a. a^{-1}w)  \nonumber   \\
&=&a^{\lambda-\rho-(n+4l)\alpha}\pi^{1_\lambda,\vp}(a)\pi^{1_\lambda,\vp}(\Delta^l) \vp(\xi) \chi(a^2a(\xi^{-1}.w)) f(\xi^{-1}.w)\nonumber\\
&=&a^{\lambda-\rho-(n+4l)\alpha}\sum_{j}
\pi^{1_\lambda,\vp}(a)  \vp(A_{l,j}\xi)  \chi(a^2 a(\xi^{-1}\tilde B_{l,j}.w)) f(\xi^{-1}\tilde C_{l,j}.w)\nonumber\\
&=&
a^{2(\lambda-\rho)-(n+4l)\alpha}\sum_{j}
\vp(a) \vp(A_{l,j}\xi)  \chi(a^2 a(\xi^{-1}\tilde B_{l,j}a^{-1}.aw)) f(\xi^{-1}\tilde C_{l,j}a^{-1}.aw)\label{zzz}\ .
\end{eqnarray}

\subsubsection{}

If $W\subset N$ is any compact subset and $D\in \cU(\naaa)$, then there
is a constant $C\in\R$ such that
$$|\pi^{1_\lambda,\vp}(\Delta^l\xi^a) L(f)(yDw)|\le C   a^{2(\lambda-\rho)-(n+4l)\alpha+l_\vp}$$
for all $\xi\in S(Z)$ and $a\in A_+$.

If $\Ree(\lambda)$ is sufficiently small, then by partial integration for
$\vartheta\in \cX(Z,{[U^0,U^0]})$, $\vartheta\not=0$
\begin{eqnarray*}
\pi^{Z,\vartheta}_*(L(f))&=&\int_{Z}
\vartheta(x)^{-1}\pi^{1_\lambda,\vp}(x) L(f) dx\\
&=&\frac{1}{|\vartheta|^{2l}}\int_{Z} \vartheta(\Delta^lx)^{-1} \pi^{1_\lambda,\vp}(x) L(f) dx\\
&=&\frac{1}{|\vartheta|^{2l}}\int_{Z} \vartheta(x)^{-1}\pi^{1_\lambda,\vp}(\Delta^l x) L(f) dx\\
&=&\frac{1}{|\vartheta|^{2l}}\int_{S(Z)}  \int_A
\vartheta(\xi^a)^{-1}\pi^{1_\lambda,\vp}(\Delta^l \xi^a) L(f)a^{2\rho_{{[U^0,U^0]}}} da d\xi\ .
\end{eqnarray*}
By the estimate above
the integral converges locally uniformly on
$\{\Ree(\lambda)<\rho^{{[U^0,U^0]}}+(n+4l)\alpha/2 -l_\vp/2\}$,
and for any compact subset of this region,
$D\in\cU(\naaa)$, and compact subset $W\subset N$
there is a constant $C_1$ such that
$$|\pi^{Z,\vartheta}_*(L(f))(yDw)|\le \frac{C_1}{ |\vartheta|^{2l}}$$
for all $y\in W$.

\subsubsection{}

Choosing $2l>\dim(Z)$ we see that the sum
$$\sum_{0\not=\vartheta\in\cX(Z,{[U^0,U^0]})}
\pi^{Z,\vartheta}_*(L(f))$$ converges in $C^\infty(\Omega_P,V(1_\lambda,\vp))$.

\subsubsection{}

Refining the estimates above we now show that this
sum in fact converges in the space of rapidly decreasing functions.
Let $W_1\subset N\setminus Z$ be a compact subset and $D\in\cU(\naaa)$.
It immediately follows from (\ref{zzz}), that there is a constant $C_2$
such that for all $\xi\in S(Z)$, $a,b\in A_+$, with $a \ge b$ and $y\in W_1$
$$|\vp(b)^{-1}(\pi^{1_\lambda,\vp}(\Delta^l\xi^a) L(f))(y^bDw)|\le
C_2 a^{2(\lambda-\rho)-(n+4l)\alpha}| \vp(ab^{-1})|  \ .$$
Let $COP(D)=\sum_{h} D_h^l\otimes D_h^r$
denote the coproduct of $D$, where here
$COP:\cU(\naaa)\rightarrow \cU(\naaa)\otimes\cU(\naaa)$.

We find a constant $C_3\in \R$ such that
for all $\xi\in S(Z)$, $a,b\in A_+$ with $a \le b$ and $y\in W_1$
\begin{eqnarray*}
&&|\vp(b)^{-1}(\pi^{1_\lambda,\vp}(\Delta^l\xi^a) L(f))(y^bDw)|\\
&=&
a^{2(\lambda-\rho)-(n+4l)\alpha}|\sum_{j,h}
\vp(ab^{-1}) \vp(A_{l,j}\xi)  \chi(a^2 a(\xi^{-1}\tilde B_{l,j}a^{-1}y^bD^l_haw)) f(\xi^{-1}\tilde C_{l,j}a^{-1}y^bD^r_haw)|\\
&=&
a^{-4l\alpha}
b^{2(\lambda-\rho)-n\alpha}
|\sum_{j,h} \vp((A_{l,j}\xi)^{ab^{-1}})  \chi(b^2 a((\xi^{-1}\tilde B_{l,j})^{b^{-1}a}y(D_h^l)^{b^{-1}}w)) f((\xi^{-1}\tilde C_{l,j})^{b^{-1}a}y(D_h^r)^{b^{-1}}w)|\\
&\le& C_3 b^{2(\lambda-\rho)-(n+4l)\alpha}\ .
\end{eqnarray*}

Using this estimate we see that
$$\vp(b)^{-1}\sum_{0\not=\vartheta\in\cX(Z,{[U^0,U^0]})}
\pi^{Z,\vartheta}_*(L(f))(y^bDw)$$
can be estimated by $C b^{2(\lambda-\rho^{{[U^0,U^0]}})-(n+4l)\alpha}$,
where $C$ can be choosen uniformly for $y\in W_1$ and $\lambda$ in compact subsets
of $\{\Ree(\lambda)<\rho^{{[U^0,U^0]}}+(n+4l)\alpha/2 -l_\vp/2\}$.
Since we can choose $l$ arbitrary large we obtain
a holomorphic continuation of the sum above to all
of $\aca$. Moreover we see that this sum
is rapidly decreasing  with respect to $b$.

\subsubsection{}

We now consider $\pi^{Z,0}_*(L(f))$.
If $\Ree(\lambda)$ is sufficiently small then we can write
$$\pi^{Z,0}_*(L(f))=\int_A \int_{S(Z)} \pi^{1_\lambda,\vp}(\xi^a) L(f) d\xi a^{2\rho_{[U^0,U^0]}} da\ .$$
We again employ the relation
\begin{eqnarray*}
&&\pi^{1_\lambda,\vp}(\xi^a) L(f)(yw)\\
&=&
a^{2(\lambda-\rho)-n\alpha}
\vp(a) \vp(\xi)  \chi(a^2 a(\xi^{-1}a^{-1}yaw)) f(\xi^{-1}a^{-1}yaw)\ .
\end{eqnarray*}
Note that
$$F(a,y):= \int_{S(Z)}\vp(\xi)  \chi(a^2 a(\xi^{-1}yw)) f(\xi^{-1}yw) d\xi$$
is a smooth function of $y$ near $y=0$ which is independent
of $a$ for large $a$.
Let
$$F(a,y^{a^{-1}})=\sum_{q=0}^r F_q(y)a^{-q\alpha} + a^{-(r+1)\alpha} R_r(a,y)$$
be the asymptotic expansion for large $a\in A$ obtained from the Taylor
series of $F(a,y)$ at $y=0$. The remainder $R_r(a,y)$ remains bounded
as $a\to\infty$. We write
$$\pi^{Z,1}_*(L(f))(yw)=I_+(y)+I_-(y)\ ,$$ where
\begin{eqnarray*}
I_+(y)&:=&\int_{A_+}  a^{2(\lambda-\rho^{[U^0,U^0]})-n\alpha}  \vp(a) F(a,y^{a^{-1}}) da\\
I_-(y)&:=&\int_{A_-}\int_{S(Z)} \pi^{1_\lambda,\vp}(\xi^a) L(f)(yw) d\xi a^{2\rho_{[U^0,U^0]}} da
\ .\end{eqnarray*}
The integral $I_-$ converges for all $\lambda\in\aca$
and defines a holomorphic family of smooth functions.
We write $I_+(y):=J_r^1(y)+J_r^2(y)$, where
\begin{eqnarray*}
J_r^1(y)&:=&  \int_{A_+}  a^{2(\lambda-\rho^{[U^0,U^0]})-n\alpha}  \vp(a) \sum_{q=0}^r F_q(y)a^{-q\alpha} da\\
J_r^2(y)&:=&  \int_{A_+}  a^{2(\lambda-\rho^{[U^0,U^0]})-n\alpha}  \vp(a)
a^{-(r+1)\alpha} R_r(a,y) da\ .
\end{eqnarray*}
The integral $J_r^2$
converges for $\Ree(\lambda)<\rho^{[U^0,U^0]}+(n+r+1)\alpha/2-l_\vp/2$
and defines a smooth function in $y$. The integral $J_r^1$ can be
evaluated:
$$J_r^1(y)=-\sum_{q=0}^r (B+2(\lambda-\rho^{[U^0,U^0]})-(n+q)\alpha)^{-1} F_q(y)\ ,$$
where $B:=\frac{d}{da}_{|a=1}\vp(a)\in \End(V_\vp)$ (compare with the proof of Lemma \ref{aaww}).
It obviously defines a meromorphic family of smooth functions.
Since we can choose $r$ arbitrary large we obtain a meromorphic continuation of
$\pi^{Z,0}_*(L(f))$ to all of $\aca$. Note that if $y\in N\setminus Z$ and
$b\in A$
is sufficiently large, then we have $\pi^{Z,0}_*(L(f))(y^bw)=\pi_*^{P_{[U^0,U^0]}}(f)(y^bw)$,
where
$\pi_*^{P_{[U^0,U^0]}}$
was discussed in Subsection \ref{asz}.
We conclude that
$\pi^{Z,0}_*(L(f))\in S_{{[U^0,U^0]}}(1_\lambda,\vp)\oplus L(\cR_{[U^0,U^0]}(1_\lambda,\vp)^n)$.

We have shown that
$$\pi^{[U^0,U^0]}_*\circ L:\cR_{\{1\}}(1_\lambda,\vp)^n\rightarrow
S_{{[U^0,U^0]}}(1_\lambda,\vp)\oplus L(\cR_{[U^0,U^0]}(1_\lambda,\vp)^n_\theta)\subset
B_{{[U^0,U^0]}}(1_\lambda,\vp)$$
has a meromorphic continuation to all of $\aca$ such that
$AS\circ \pi^{[U^0,U^0]}_*\circ L = [\pi^{[U^0,U^0]}_*]$.
This finishes the proof of the lemma.
\hB

\subsubsection{} \label{weih258}

We now prove Lemma \ref{lpro2} in a similar manner.
Since $M_{U^0}$ centralizes $N_V$ we can form the  direct product $P_1:=P_{U^0}/Z=M_{U^0}T$ where
$T:=N_V/Z$ (the group $Z$ was defined in \ref{weih256}).  
By construction of $U^0$
there is a lattice $V_1 \subset T$ and a
homomorphism $m:T\rightarrow M_{U^0}$ 
such that ${U^0/{[U^0,U^0]}}=\{m(v)v|v\in {V_1}\}\subset P_1$.
For $f\in \cR_{[U^0,U^0]}(1_\lambda,\vp)^n_\theta$ (this space is defined in \ref{weih257}) we have
$$\pi^{1_\lambda,\vp}(m(v)v)L(f)(yw)= \theta^m(v) \vp(\tilde v) f(\tilde v^{-1}yw)\ ,$$
where $\tilde v\in N$ is any lift of $v\in T$, and $\theta^m=\theta\circ m\in \cX(T)$.

\subsubsection{}

We normalize the Haar measure on $T$ such that
$\vol({V_1}\backslash T)=1$.
Given a unitary character
$\vartheta\in \cX(T)$,
we consider the push-down $\pi^{T,\vartheta}_*$,
which associates to $f\in C^\infty(B_Z,V_{B_Z}(1_\lambda,\vp))$
the average
$$\pi^{T,\vartheta}_*(f):=\int_{T} \vartheta(z)^{-1}  \pi^{1_\lambda,\vp}(z) f  dz$$
provided the integral converges.
Note that we define $L:\cR_{[U^0,U^0]}(1_\lambda,\vp)_\theta^n\rightarrow
B_{[U^0,U^0]}(1_\lambda,\vp)$ using a $Z$-invariant
cut-off function (recall that $Z=P_{[U^0,U^0]}$). It follows from the proof of Lemma \ref{lpro1}
that all $f\in \cR_{[U^0,U^0]}(1_\lambda,\vp)_\theta^n$ are $Z$-invariant.
Therefore $L(f)\in C^\infty(B_Z,V_{B_Z}(\theta,\vp))$ and
$\pi^{T,\vartheta}_*(L(f))$ is well-defined (up to convergence).

\subsubsection{}

 We fix an isomorphism of abelian groups $T\cong \R^{\dim(T)}$
preserving the Haar measure
such that we obtain an euclidean structure on $T$ and its Lie algebra $\taaa$.
As before let  $\cX(T,{V_1})$ denote the set of all unitary characters of
$T$ which are trivial on ${V_1}$. We identify $\cX(T)$
with $\imath \taaa^*$ and $\cX(T,{V_1})$ with a lattice
in $\imath \taaa^*$. If $f\in \cR_{[U^0,U^0]}(1_\lambda,\vp)^n_\theta$, then
by the Poisson summation formula
$$\pi^{U^0/{[U^0,U^0]}}_*(L(f))=\sum_{\vartheta\in\cX(T,{V_1})}
\pi^{T,\vartheta-\theta^m}_*(L(f))$$ in the domain of convergence.

\subsubsection{}

Note that if $\theta^m\in \cX(T,{V_1})$, then $\theta=0$
by the construction of $M_{U^0}$.
If $\theta\not=0$, then will show that
$\sum_{\vartheta\in\cX(T,{V_1})}
 \pi^{T,\vartheta-\theta^m}_*(L(f))$
gives rise to a rapidly decreasing section for all $\lambda\in \aca$.
If $\theta=0$, then again $\sum_{0\not=\vartheta\in\cX(T,{V_1})}
 \pi^{T,\vartheta-\theta^m}_*(L(f))$ is rapidly decreasing.
The remaining term $\pi^{T,0}_*(L(f))$ contributes to
$\cS_{U^0}(1_\lambda,\vp)\oplus L(\cR_{U^0}(1_\lambda,\vp)^n_\theta)$
and depends meromorphically on $\lambda$.

\subsubsection{}

Observe that  $A$ normalizes $Z$. This is a consequence  of the assumption that the cusp is regular, and that the Langlands decomposition of $P$ is adapted   (Definition \ref{t799}). It is in fact  the reason for making this assumption. We see that $A$  acts on $T$.
In particular $\taaa^*$ decomposes into two
eigenspaces $\taaa=\taaa_1^*\oplus\taaa_2^*$
(we write $\vartheta=\vartheta_1\oplus \vartheta_2$ for the corresponding decomposition of $\vartheta\in\taaa^*$)
with respect to $A$ such that $A$ acts on $\taaa_i^*$ by $a^{i\alpha}$, $i=1,2$.
We define the $A$-homogeneous "norm" on $\taaa^*$ by
 $|\vartheta|:=(|\vartheta_1|^4+|\vartheta_2|^2)^{1/2}$.
Furthermore we use
the euclidean structure on $\taaa_i$
in order to define  Laplace operators
$\Delta_i\in \cU(\taaa_i)$ and the $A$-homogeneous operator
$\Delta:=\Delta_1^2+\Delta_2$. We fix the normalizations such that
for $\vartheta\in\cX(T)$ we have
$\vartheta(\Delta^l x)=|\vartheta|^{2l}
\vartheta(x)$.
Let $\rho_{U^0/{[U^0,U^0]}}\in\aaaa^*$ be such that $A$ acts on
$\Lambda^{max}\taaa$ by the character $a^{2\rho_{U^0/{[U^0,U^0]}}}$.
Note that $\rho_{U^0/{[U^0,U^0]}}+\rho_{[U^0,U^0]}=\rho_{U}$.

\subsubsection{}

Let $COP:\cU(\taaa)\rightarrow \cU(\taaa)\otimes \cU(\taaa)$ denote
the coproduct on $\cU(\taaa)$. There are $A_{l,j},B_{l,j},C_{l,j}\in\cU(\taaa)$
such that
$$(COP\otimes 1)\circ COP(\Delta^l)=\sum_{j} A_{l,j}\otimes B_{l,j}\otimes C_{l,j}\ .$$
Let $D\mapsto \tilde D$ be the canonical antiautomorphism of $\cU(\taaa)$.

Recall that $L(f)(xw)=\chi(xw) f(xw)$, where $\chi$ is $Z$-invariant.
By $S(T)$ we denote the unit-sphere in $T$.
For $\xi\in S(T)$ and $a\in A$ we compute
\begin{eqnarray}
&&\pi^{1_\lambda,\vp}(\Delta^l\xi^a) L(f)(.w)\nonumber\\
&=&\pi^{1_\lambda,\vp}(a)\pi^{1_\lambda,\vp}((\Delta^l)^{a^{-1}}\xi) \pi^{1_\lambda,\vp}(a^{-1})\chi f(.w)\nonumber\\
&=&a^{\rho-\lambda-4l\alpha}\pi^{1_\lambda,\vp}(a)\pi^{1_\lambda,\vp}(\Delta^l\xi) \vp(a^{-1}) \chi(a.a^{-1}w) f(a. a^{-1}w)  \nonumber   \\
&=&a^{\lambda+\rho-2\rho^{[U^0,U^0]}-(n+4l)\alpha}\pi^{1_\lambda,\vp}(a)\pi^{1_\lambda,\vp}(\Delta^l) \vp(\xi) \chi(a\xi^{-1}.a^{-1}w) f(\xi^{-1}.w)\nonumber\\
&=&a^{\lambda+\rho-2\rho^{[U^0,U^0]}-(n+4l)\alpha}\sum_{j}
\pi^{1_\lambda,\vp}(a)  \vp(A_{l,j}\xi)  \chi(a\xi^{-1}\tilde B_{l,j}.a^{-1}w) f(\xi^{-1}\tilde C_{l,j}.w)\nonumber\\
&=&
a^{2(\lambda-\rho^{[U^0,U^0]})-(n+4l)\alpha}\sum_{j}
\vp(a) \vp(A_{l,j}\xi)  \chi(a\xi^{-1}\tilde B_{l,j}a^{-1}.w) f(\xi^{-1}\tilde C_{l,j}a^{-1}.aw)\label{zzz1}\ .
\end{eqnarray}
If $W\subset N$ is any compact subset and $D\in \cU(\naaa)$, then there
is a constant $C\in\R$ such that
$$|\pi^{1_\lambda,\vp}(\Delta^l\xi^a) L(f)(yDw)|\le C   a^{2(\lambda-\rho^{[U^0,U^0]})-(n+4l)\alpha+l_\vp}$$
for all $\xi\in S(T)$ and $a\in A_+$.

\subsubsection{}

If $\Ree(\lambda)$ is sufficiently small, then by partial integration for
$\vartheta\in \cX(T,{V_1})$, $\vartheta-\theta^m\not=0$
\begin{eqnarray*}
\pi^{Z,\vartheta-\theta^m}_*(L(f))&=&\int_{Z}
(\vartheta-\theta^m)^{-1}(x)\pi^{1_\lambda,\vp}(x) L(f) dx\\
&=&\frac{1}{|\vartheta-\theta^m|^{2l}}\int_{T} (\vartheta-\theta^m)^{-1}(\Delta^lx) \pi^{1_\lambda,\vp}(x) L(f) dx\\
&=&\frac{1}{|\vartheta-\theta^m|^{2l}}\int_{T} (\vartheta-\theta^m)^{-1}(x)\pi^{1_\lambda,\vp}(\Delta^l x) L(f) dx\\
&=&\frac{1}{|\vartheta-\theta^m|^{2l}}\int_{S(T)}  \int_A
(\vartheta-\theta^m)^{-1}(\xi^a)\pi^{1_\lambda,\vp}(\Delta^l \xi^a) L(f)a^{2\rho_{{U^0/{[U^0,U^0]}}}} da d\xi\ .
\end{eqnarray*}
By the estimate above
the integral converges locally uniformly on
$\{\Ree(\lambda)<\rho^{U}+(n+4l)\alpha/2 -l_\vp/2\}$,
and for any compact subset of this region,
$D\in\cU(\naaa)$, and compact subset $W\subset N$
there is a constant $C_1$ such that
$$|\pi^{T,\vartheta-\theta^m}_*(L(f))(yDw)|\le \frac{C_1}{ |\vartheta-\theta^m|^{2l}}$$
for all $y\in W$.

\subsubsection{}

Choosing $2l>\dim(T)$  we see that the sum
$$\sum_{\vartheta\in\cX(T,{V_1}), \vartheta-\theta^m\not=0}
\pi^{T,\vartheta-\theta^m}_*(L(f))$$ converges in $C^\infty(\Omega,V(1_\lambda,\vp))$.

\subsubsection{}

Refining the estimates above we now show that this
sum in fact converges in the space of rapidly decreasing functions.
Let $W_1\subset N\setminus N_V$ be a compact subset and $D\in\cU(\naaa)$.
It immediately follows from (\ref{zzz1}), that there is a constant $C_2$
such that for all $\xi\in S(T)$, $a,b\in A_+$, with $a \ge b$ and $y\in W_1$
$$|\vp(b)^{-1}(\pi^{1_\lambda,\vp}(\Delta^l\xi^a) L(f))(y^bDw)|\le
 C_2 a^{2(\lambda-\rho^{[U^0,U^0]})-(n+4l)\alpha}| \vp(ab^{-1})|  \ .$$

Let $COP(D)=\sum_{h} D_h^l\otimes D_h^r$
denote the coproduct of $D$, where here
$COP:\cU(\naaa)\rightarrow \cU(\naaa)\otimes\cU(\naaa)$.

We find a constant $C_3$ such that
for all $\xi\in S(T)$, $a,b\in A_+$, with $a \le b$ and $y\in W_1$
\begin{eqnarray*}
&&|\vp(b)^{-1}(\pi^{1_\lambda,\vp}(\Delta^l\xi^a) L(f))(y^bDw)|\\
&=&
a^{2(\lambda-\rho^{[U^0,U^0]})-(n+4l)\alpha}|\sum_{j,h}
\vp(ab^{-1}) \vp(A_{l,j}\xi)  \chi(a\xi^{-1}\tilde B_{l,j}a^{-1}y^bD^l_hw) f(\xi^{-1}\tilde C_{l,j}a^{-1}y^bD^r_haw)\\
&=&
a^{-4l\alpha}
b^{2(\lambda-\rho^{[U^0,U^0]})-n\alpha}
|\sum_{j,h} \vp((A_{l,j}\xi)^{ab^{-1}})  \chi(a\xi^{-1}\tilde B_{l,j}a^{-1}y^bD^l_hw) f((\xi^{-1}\tilde C_{l,j})^{b^{-1}a}y(D^r_h)^{b^{-1}}w)\\
&\le& C_3 b^{2(\lambda-\rho^{[U^0,U^0]})-(n+4l)\alpha}\ .
\end{eqnarray*}

Using this estimate we see that
$$\vp(b)^{-1}\sum_{\vartheta\in\cX(T,{V_1}),\vartheta-\theta^m\not=0}
\pi^{T,\vartheta-\theta^m}_*(L(f))(y^bDw)$$
can be estimated by $C b^{2(\lambda-\rho^{U})-(n+4l)\alpha}$,
where $C$ can be choosen uniformly for $y\in W_1$ and $\lambda$ in compact subsets
of $\{\Ree(\lambda)<\rho^{U}+(n+4l)\alpha/2 -l_\vp/2\}$.
Since we can choose $l$ arbitrary large we obtain
a holomorphic continuation of the sum above to all
of $\aca$ and that it is rapidly decreasing with respect to $b$.

\subsubsection{}

Let now $\theta=0$ and consider
$\pi^{T,0}_*(L(f))$.
We write for $\Ree(\lambda)$ sufficiently small
$$\pi^{T,0}_*(L(f))=\int_A \int_{S(T)} \pi^{1_\lambda,\vp}(\xi^a) L(f) d\xi a^{2\rho_{U^0/{[U^0,U^0]}}} da\ .$$
We again employ
\begin{eqnarray*}
&&\pi^{1_\lambda,\vp}(\xi^a) L(f)(yw)\\
&=&
a^{2(\lambda-\rho^{[U^0,U^0]})-n\alpha}
\vp(a) \vp(\xi)  \chi(a\xi^{-1}y^{a^{-1}}a^{-1}w) f(\xi^{-1}a^{-1}yaw)\ .
\end{eqnarray*}
Note that
$$F(a,y):= \int_{S(T)}\vp(\xi)  \chi(a\xi^{-1}ya^{-1}w) f(\xi^{-1}yw) d\xi$$
is a smooth function of $y$ near $y=0$ which is independent of $a$ for large $a$. 
Let
$$F(a,y^{a^{-1}})=\sum_{q=0}^r F_q(y)a^{-q\alpha} + a^{-(r+1)\alpha} R_r(a,y)$$
be the asymptotic expansion for large $a\in A$ obtained from the Taylor
series of $F(a,y)$ at $y=0$.
The remainder $R_r(a,y)$ remains bounded
as $a\to\infty$. We write
$$\pi^{T,0}_*(L(f))(yw)=I_+(y)+I_-(y)\ ,$$ where
\begin{eqnarray*}
I_+(y)&=&\int_{A_+}  a^{2(\lambda-\rho^{U})-n\alpha}  \vp(a) F(a,y^{a^{-1}}) da\\
I_-(y)&=&\int_{A_-}\int_{S(T)} \pi^{1_\lambda,\vp}(\xi^a) L(f)(yw) d\xi a^{2\rho_{U^0/{[U^0,U^0]}}} da\ .
\end{eqnarray*}
The integral $I_-$ converges for all $\lambda\in\aca$
and defines a holomorphic family of smooth functions.
We write $I_+(y):=J_r^1(y)+J_r^2(y)$, where
\begin{eqnarray*}
J_r^1(y)&:=&  \int_{A_+}  a^{2(\lambda-\rho^U)-n\alpha}  \vp(a) \sum_{q=0}^r F_q(y)a^{-q\alpha} da\\
J_r^2(y)&:=&  \int_{A_+}  a^{2(\lambda-\rho^U)-n\alpha}  \vp(a)
a^{-(r+1)\alpha} R_r(a,y) da\ .
\end{eqnarray*}
The integral $J_r^2$
converges for $\Ree(\lambda)<\rho^U+(n+r+1)\alpha/2-l_\vp/2$
and defines a smooth function in $y$. The integral $J_r^1$ can be
evaluated:
$$J_r^1(y)=\sum_{q=0}^r (B+2(\lambda-\rho^U)-(n+q)\alpha)^{-1} F_q(y)\ .$$
It obviously defines a meromorphic family of smooth functions.
Since we can choose $r$ arbitrary large we obtain a meromorphic continuation of
$\pi^{U,1}_*(L(f))$ to all of $\aca$.

\subsubsection{}

Let now $f\in \cR_{\{1\}}(1_\lambda,\vp)_0^n$.
If $y\in N\setminus T$ and
$b\in A$
is sufficiently large, then we have
\begin{eqnarray*}
\pi^{T,0}_*\circ L\circ [\pi^{U^0/{[U^0,U^0]}}_*](f)(y^b)&=&
\pi^{T,0}_*\circ \pi^{Z,0}_*(f)(y^b)\\
&=&\pi^{P_{U^0}}_*(f)(y^b)\\
&=&[\pi^{U^0}_*](f)(y^b)\ .
\end{eqnarray*}
We now have shown that if $f\in \cR_{[U^0,U^0]}(1_\lambda,\vp)^n_\theta$, then
$\pi^{U^0/[U^0,U^0]}_*(f)\in \cS_{U^0}(1_\lambda,\vp)\oplus L(\cR_{U^0}(1_\lambda,\vp)^n)$,
and that $\pi^{U^0/[U^0,U^0]}_*$ depends meromorphically on $\lambda$. \hB

\subsubsection{}

The finite group $U/U^0$ acts on $\cS_{U^0,k}(1_\lambda,\vp)$
and $\cR_{U^0}(1_\lambda,\vp)^n$. We can define
$\{\pi^{U/U^0}\}:\cS_{U^0,k}(1_\lambda,\vp)\rightarrow
\cS_{U,k}(1_\lambda,\vp)$, $\pi^{U/U^0}:B_{U^0}(1_\lambda,\vp)\rightarrow
B_U(1_\lambda,\vp)$, and $[\pi^{U/U^0}_*]:\cR_{U^0}(1_\lambda,\vp)^n\rightarrow
\cR_U(1_\lambda,\vp)^n$.
Using  $\pi^{U}_*:=\pi^{U/U^0}_*\circ \pi^{U^0}_*$
and Proposition \ref{mainpure} (and its proof) we obtain
 \begin{prop}\label{muncor}
The composition
$$\pi^{U}_*\circ L:\cR_{\{1\}}(1_\lambda,\vp)^n_\theta\rightarrow
B_{U}(1_\lambda,\vp)$$
converges for $\Ree(\lambda)<\rho^U+(n\alpha-l_\vp)/2$
and has a meromorphic continuation to all of $\aca$
such that $AS\circ \pi^{U}_*\circ L = [\pi^{U}_*]$.
If $\theta\not=0$, then it is in fact
 a holomorphic family of maps
$\pi^{U}_*\circ L:\cR_{\{1\}}(1_\lambda,\vp)^n_\theta\rightarrow
\cS_{U}(1_\lambda,\vp)$. If $\theta=0$,
then $\pi^{U}_*\circ L$ has at most first order poles in the set $\rho^U+\frac{n\alpha-l_\vp}{2} +\frac{1}{2}\nat_0$.  
\end{prop}

\subsubsection{}

Combining Proposition \ref{muncor} with Lemma \ref{schwpush}
we obtain
\begin{kor}\label{u76}
Assume that the cusp associated to $U\subset P$ does not have full rank.
Then for any $k_1< k$ the push-down
$\pi_*^U: B_{\{1\},k}(1_\lambda,\vp) \rightarrow B_{U,k_1}(1_\lambda,\vp)$
forms a meromorphic family of continuous maps
with finite-dimensional singularities
defined on $\{\Ree(\lambda)<\rho^U+(k\alpha-l_\vp)/2\}$.
It fits into the following commutative diagram
$$
\begin{array}{ccccccccc}
0&\rightarrow&\cS_{\{1\},k}(1_\lambda,\vp)&\rightarrow& B_{\{1\},k}(1_\lambda,\vp)&\stackrel{AS}{\rightarrow}& \cR_{\{1\},k}(1_\lambda,\vp)&\rightarrow &0\\
 &           &\downarrow \{\pi^{U}_*\} &           & \downarrow \pi^{U}_* &      & \downarrow [\pi^{U}_*]          &&\\
0&\rightarrow&\cS_{U,k}(1_\lambda,\vp)&\rightarrow& B_{U,k}(1_\lambda,\vp)& \stackrel{AS}{\rightarrow}&\cR_{U,k}(1_\lambda,\vp)&\rightarrow &0\ .\end{array}
$$
Moreover we have a  meromorphic family of maps
$\pi_*^U:B_{\{1\}}(1_\lambda,\vp)   \rightarrow    B_{U}(1_\lambda,\vp)$
with finite-dimensional singularities and defined on all of $\aca$
such that
$$
\begin{array}{ccccccccc}
0&\rightarrow&\cS_{\{1\}}(1_\lambda,\vp)&\rightarrow& B_{\{1\}}(1_\lambda,\vp)&\stackrel{AS}{\rightarrow}& \cR_{\{1\}}(1_\lambda,\vp)&\rightarrow &0\\
 &           &\downarrow \{\pi^{U}_*\} &           & \downarrow \pi^{U}_* &      & \downarrow [\pi^{U}_*]          &&\\
0&\rightarrow&\cS_{U}(1_\lambda,\vp)&\rightarrow& B_{U}(1_\lambda,\vp)& \stackrel{AS}{\rightarrow}&\cR_{U}(1_\lambda,\vp)&\rightarrow &0\ .\end{array}
$$
is commutative.
\end{kor}

\subsection{Push-down for cusps of full rank and for general $\sigma$} \label{maxrank}

\subsubsection{}

In this subsection $(\sigma,V_\sigma)$ denotes a Weyl-invariant
representation of $M$ as in \ref{weih241}.

Let  $C^{-\infty}(\infty_P,V(\sigma_\lambda,\vp))\subset C^{-\infty}(\partial X,V(\sigma_\lambda,\vp))$
be the space of distributions which are supported on $\infty_P$.
This space can be identified $P_U$-equivariantly with the tensor product of 
a generalized Verma module by $V_\vp$
$$(\cU(\gaaa)\otimes_{\cU(\paaa)} V_{\sigma_{\lambda+2\rho}})\otimes V_\vp$$
such that $(X\otimes s\otimes v)\in (\cU(\gaaa)\otimes_{\cU(\paaa)} V_{\sigma_{\lambda+2\rho}})\otimes V_\vp$  maps $f\in C^\infty(\partial X,V(\tilde \sigma_{-\lambda},\tilde\vp))$ to $(s\otimes v)(f(Xe))$.
We can further identify
$C^{-\infty}(\infty_P,V(\sigma_\lambda,\vp))$ with $R_{\{1\}}(\tilde\sigma_{-\lambda},\tilde\vp)^*$.

\subsubsection{}\label{weih121}

We now assume that the cusp associated to $U\subset P$
has full rank. In this case $N$ is the Zariski closure of $U^0$.  
The space $(\cU(\gaaa)\otimes_{\cU(\paaa)} V_{\sigma_{\lambda+2\rho}})\otimes V_\vp$ carries an algebraic representation
of $N$. Therefore
\begin{eqnarray*}
{}^U[(\cU(\gaaa)\otimes_{\cU(\paaa)} V_{\sigma_{\lambda+2\rho}})\otimes V_\vp]&\subset&
{}^{U^0}[(\cU(\gaaa)\otimes_{\cU(\paaa)} V_{\sigma_{\lambda+2\rho}})\otimes V_\vp]\\
&=&{}^N[(\cU(\gaaa)\otimes_{\cU(\paaa)} V_{\sigma_{\lambda+2\rho}})\otimes V_\vp]
 \end{eqnarray*}
is finite-dimensional since the space of highest weight vectors of 
the $\gaaa$-module
$\cU(\gaaa)\otimes_{\cU(\paaa)} V_{\sigma_{\lambda+2\rho}}$ is finite-dimensional \cite{MR552943}.

\subsubsection{}\label{weih300}
Let $\cE_{\infty_P}(\sigma,\vp)$ denote the sheaf of holomorphic families
$f_\nu\in{}^UC^{-\infty}(\infty_P,V(\sigma_\nu,\vp))$.
Since $\cE_{\infty_P}(\sigma,\vp)$ is torsion-free it is the space of sections
of a unique holomorphic vector bundle $E_{\infty_P}(\sigma,\vp)$ over $\aca$.
By $E_{\infty_P}(\sigma_\lambda,\vp)$ we denote the fibre of $E_{\infty_P}(\sigma,\vp)$ at $\lambda\in\aca$.
We will discuss this bundle in detail in Lemma \ref{vermaolbrich}.

\subsubsection{}

We now define the function space $B_U(\sigma_\lambda,\vp)$
in the case of a cusp of full rank.

\begin{ddd}
We define $R_U(\sigma_\lambda,\vp):= E_{\infty_P}(\tilde\sigma_{-\lambda},\tilde\vp)^*$.
Furthermore we set $$B_U(\sigma_\lambda,\vp):=\cS_U(\sigma_\lambda,\vp)\oplus \cR_U(\sigma_\lambda,\vp)\ .$$
Let $AS:B_U(\sigma_\lambda,\vp) \rightarrow \cR_U(\sigma_\lambda,\vp)$
be the projection and $L:\cR_U(\sigma_\lambda,\vp)\rightarrow B_U(\sigma_\lambda,\vp)$
be the inclusion. 
\end{ddd}
These families of spaces form trivial holomorphic bundles of Fr\'echet and Montel spaces over $\aca$.

\subsubsection{}

\begin{ddd}\label{rolf}
We define $$[ext^U]:E_{\infty_P}(\tilde\sigma_{-\lambda},\tilde\vp)\hookrightarrow
\cR_{\{1\}}(\tilde\sigma_{-\lambda},\tilde\vp)^*$$ as the natural inclusion.
We define the push-down
$$[\pi^U_*]:\cR_{\{1\}}(\sigma_{\lambda},\vp)\rightarrow 
R_U(\sigma_\lambda,\vp)$$ to be the adjoint of $[ext^U]$.
Furthermore we define
$$(\pi^U_*)_1:C^\infty(\partial X,V(\sigma_\lambda,\vp))\rightarrow  
\cS_U(\sigma_\lambda,\vp)$$ by $(\pi^U_*)_1(f):=\sum_{u\in U}\pi^{\sigma_\lambda,\vp}(u)(f_{|\Omega})$
(convergence provided). Finally we set
$$\pi^U_*:=(\pi^U_*)_1\oplus [\pi^U_*]:C^\infty(\partial X,V(\sigma_\lambda,\vp))\rightarrow B_U(\sigma_\lambda,\vp) \ .$$
\end{ddd}
Note that
$[\pi^U_*]$ is a holomorphic family of surjective maps.
Once and for all we fix a right-inverse $[Q]$ of $[\pi^U_*]$.
We leave it to the interested reader to show that the
statements of Lemma \ref{compat2} hold true for cusps of full rank as well.

\subsubsection{}

\begin{lem}
The push-down 
$(\pi^U_*)_1:C^\infty(\partial X,V(1_\lambda,\vp))\rightarrow
\cS_U(1_\lambda,\vp)$ converges for $\Ree(\lambda)<-l_\vp/2$
and has a meromorphic continuation to all of $\aca$
with finite-dimensional singularities.
\end{lem}
\proof
One checks that the corresponding parts of the proofs
of Lemma \ref{schwpush}  and Proposition \ref{mainpure} apply to the case
of cusps of full rank as well. \hB

\subsubsection{}

Note that in Subsection \ref{weih261} we have considered the push-down in the spherical case $\sigma=1$. We now deal with the general case using the concept of twisting.
We  assume that $U\subset P$ defines a cusp of smaller rank and consider a Weyl invariant  $\sigma$ (see \ref{weih241}).
We show the existence of a meromorphic family of push-down maps
$$\pi^{U}_*:C^\infty(\partial X,V(\sigma_\lambda,\vp))\rightarrow
B_U(\sigma_\lambda,\vp)$$ using twisting.
We employ a finite-dimensional representation $(\pi_{\sigma,\mu},V_{\pi_{\sigma,\mu}})$ of $G$ as in \ref{weih144}.
Note that $$\pi^U_*:C^\infty(\partial X,V(1_{\lambda+\mu},\pi_{\sigma,\mu}\otimes\vp))
\rightarrow B_U(1_{\lambda+\mu},\pi_{\sigma,\mu}\otimes\vp)$$
is $\cZ(\gaaa)$-equivariant. Therefore we can make the following definition:
\begin{ddd}\label{cv1}
If $\lambda\not\in I_\aaaa$, then using Lemma \ref{compat2}, 4.,
we define the push-down
$$\pi^U_*:C^\infty(\partial X,V(\sigma_\lambda,\vp))\rightarrow 
B_U(\sigma_\lambda,\vp)$$
by the following coummutative diagram:
$$\begin{array}{ccccccc}
0&\rightarrow&C^\infty(\partial X,V(\sigma_\lambda,\vp))&\stackrel{i_{\sigma,\mu}}{\rightarrow}&C^\infty(\partial X,V(1_{\lambda+\mu},\pi_{\sigma,\mu}\otimes\vp))&\stackrel{Z(\lambda)}{\rightarrow}&
C^\infty(\partial X,V(1_{\lambda+\mu},\pi_{\sigma,\mu}\otimes\vp))\\
&&\downarrow \pi^U_*&&\downarrow \pi^U_*&&\downarrow \pi^U_*\\
0&\rightarrow&B_U(\sigma_\lambda,\vp)&\stackrel{i^U_{\sigma,\mu}}{\rightarrow}&
\cB_U(1_{\lambda+\mu},\pi_{\sigma,\mu}\otimes\vp)&\stackrel{Z(\lambda)}{\rightarrow}&
\cB_U(1_{\lambda+\mu},\pi_{\sigma,\mu}\otimes\vp)\end{array}\ .$$
\end{ddd}

\subsubsection{}

It is clear that in the domain of convergence $\Ree(\lambda)\ll 0$
this definition coincides with Definition \ref{pushdowndef}.
In order to see that
$\pi^U_*$ extends to a meromorphic family note that
$i^U_{\sigma,\mu}$ admits a meromorphic family of
left inverses $j^U_{\sigma,\mu}$ (Lemma \ref{compat2}, 5.), and that we can express
the push-down for $\sigma$ through the spherical push-down
by $j^U_{\sigma,\mu}\circ \pi^U_* \circ i^U_{\sigma,\mu}$.

\subsubsection{}

We claim that $\pi^U_*$ has finite-dimensional
singularities. Note that given $\nu\in \aaaa^*$
there exists $k\in\nat_0$ such that the restriction of
$\pi^U_*$ to $S_{\{1\},k}(\sigma_\lambda,\vp)\cap B_{\{1\}}(\sigma_\lambda,\vp)$
converges for all $\Ree(\lambda)<\nu$.
The rank of the singularities of
$\pi^U_*$ for $\Ree(\lambda)<\nu$ is bounded
by the codimension of
$S_{\{1\},k}(\sigma_\lambda,\vp)\cap B_{\{1\}}(\sigma_\lambda,\vp)$
in $B_{\{1\}}(\sigma_\lambda,\vp)$ which is finite. This proves the claim.

\subsubsection{}

The following corollary is a consequence of the discussion
above and Corollary \ref{u76}.
\begin{kor}\label{neuj106}
Let $U\subset P$ define a regular cusp and $\vp$ be an admissible twist.
Then the push-down is a meromorphic family of
continuous maps $\pi_*^U:B_{\{1\}}(\sigma_\lambda,\vp)   \rightarrow    B_{U}(\sigma_\lambda,\vp)$
with finite-dimensional singularities such that
$$
\begin{array}{ccccccccc}
0&\rightarrow&\cS_{\{1\}}(\sigma_\lambda,\vp)&\rightarrow& B_{\{1\}}(\sigma_\lambda,\vp)&\stackrel{AS}{\rightarrow}& \cR_{\{1\}}(\sigma_\lambda,\vp)&\rightarrow &0\\
 &           &\downarrow \{\pi^{U}_*\} &           & \downarrow \pi^{U}_* &      & \downarrow [\pi^{U}_*]          &&\\
0&\rightarrow&\cS_{U}(\sigma_\lambda,\vp)&\rightarrow& B_{U}(\sigma_\lambda,\vp)& \stackrel{AS}{\rightarrow}&\cR_{U}(\sigma_\lambda,\vp)&\rightarrow &0\ .\end{array}
$$
is commutative.
\end{kor}

\subsection{Compatibility with embedding}\label{weih116}

\subsubsection{}

In this subsection we assume that $G^n$ belongs to the list
$$\{Spin(1,n), SO(1,n)_0, SU(1,n), Sp(1,n)\}\ .$$
First assume that $U\subset P$ defines a cusp of smaller rank.
Then we have a commutative diagram
$$\begin{array}{ccc} A_{\{1\}}(1^{n+1}_{\lambda},\vp)^m&\stackrel{i^*}{\rightarrow}
&A_{\{1\}}(1^{n}_{\lambda-\zeta},\vp)^m\\
\downarrow \pi^{P^{n+1}_U}_*&&\downarrow \pi^{P^n_U}_*\\
A_{P_U}(1^{n+1}_{\lambda},\vp)^m&\stackrel{i_U^*}{\rightarrow}
&A_{P_U}(1^{n}_{\lambda-\zeta},\vp)^m
\end{array}\ .$$
In fact, commutativity is obvious in the domain of convergence.

\subsubsection{}

By the definition of the spaces 
$R_{\{1\}}(1^{n+1}_{\lambda},\vp)$, $R_{\{1\}}(1^{n}_{\lambda-\zeta},\vp)$
as the spaces of asymptotics of smooth sections and the fact that
that $i^*$ maps smooth sections to smooth sections we obtain $U$-equivariant maps $R_{\{1\}}(1^{n+1}_{\lambda},\vp)^m\stackrel{i^*}{\rightarrow}
R_{\{1\}}(1^{n}_{\lambda-\zeta},\vp)^m$, $m\in\nat_0$.

It follows from the naturality of the image bundle Lemma \ref{bunle}
that we have a commutative diagram
$$\begin{array}{ccc} R_{\{1\}}(1^{n+1}_{\lambda},\vp)^m&\stackrel{i^*}{\rightarrow}
&R_{\{1\}}(1^{n}_{\lambda-\zeta},\vp)^m\\
\downarrow [\pi^{U,n+1}_*]&&\downarrow[\pi^{U,n}_*]\\
R_U(1^{n+1}_{\lambda},\vp)^m&\stackrel{i_U^*}{\rightarrow}
&R_{U}(1^{n}_{\lambda-\zeta},\vp)^m
\end{array}\ .$$

\subsubsection{}

On the level of Schwartz spaces we have  for all $k\in\nat_0$
$$\begin{array}{ccc} S_{\{1\},k}(1^{n+1}_{\lambda},\vp)&\stackrel{i^*}{\rightarrow}
&S_{\{1\},k}(1^{n}_{\lambda-\zeta},\vp)\\
\downarrow \{\pi^{U,n+1}_*\}&&\downarrow\{\pi^{U,n}_*\}\\
S_{U,k}(1^{n+1}_{\lambda},\vp)&\stackrel{i_U^*}{\rightarrow}
&S_{U,k}(1^{n}_{\lambda-\zeta},\vp)
\end{array}\ .$$
In order to see that $i^*$ and $i^*_U$ (initially defined on spaces of smooth sections) 
induce maps between Schwartz spaces one checks that these maps are bounded with respect to the norms
$\|.\|_{k,d}$ introduced in \ref{weih259}. Commutativity of the diagram is clear.

\subsubsection{}
We now easily obtain the diagram 
$$\begin{array}{ccc} B_{\{1\},k}(1^{n+1}_{\lambda},\vp)&\stackrel{i^*}{\rightarrow}
&B_{\{1\},k}(1^{n}_{\lambda-\zeta},\vp)\\
\downarrow \pi^{U,n+1}_*&&\downarrow\pi^{U,n}_* \\
B_{U,k}(1^{n+1}_{\lambda},\vp)&\stackrel{i_U^*}{\rightarrow}
&B_{U,k}(1^{n}_{\lambda-\zeta},\vp)
\end{array}\ .$$
It gives Proposition \ref{rrttee} in the case of a pure cusp of lower rank.
 
\subsubsection{} 
 
We now assume that $U\subset P^n$ defines a cusp of full rank.
We are going to construct a map
$i_U^*:B_U(1^{n+1}_{\lambda},\vp)\rightarrow
B_U(1^n_{\lambda-\zeta},\vp)$ as the direct sum of 
$$(i_U^*)_1:B_U(1^{n+1}_{\lambda },\vp)\rightarrow
S_U(1^{n}_{\lambda },\vp)$$ and
$$(i_U^*)_2:B_U(1^{n+1}_{\lambda },\vp)\rightarrow
R_U(1^{n}_{\lambda-\zeta},\vp)\ .$$
While $(i_U^*)_1$ is just the restriction of sections
the definition of $(i_U^*)_2$ is more complicated.
In order to verify its properties we need some results of Subsection \ref{samel}.

\subsubsection{}

In the following we discuss  $(i_U^*)_2$.
This map comes as a meromorphic family and will essentially be fixed
by the condition 
\begin{equation}\label{tobsucht}\pi^{U,n}_*\circ i^*= i_U^*\circ \pi^{U,n+1}_*\ .\end{equation}
The details are as follows.

\subsubsection{} 

Since $E_{\infty_P}(1^{n}_{-\lambda+\zeta},\tilde\vp)$ (see \ref{weih300} for notation) is finite-dimensional it consists of distributions of uniformly bounded order. Further, since these distributions are supported in $\infty_P$
we can choose
 $k\in \nat$ be such that
$E_{\infty_P}(1^{n}_{-\mu+\zeta},\tilde\vp)$ pairs trivially with the space $\cS_{\{1\},k-2\rho_U}(1^n_{\mu-\zeta}, \vp)$ for all $\mu$ in some compact neighbourhood of $\lambda$. Thus  we have an inclusion
 $E_{\infty_P}(1^{n}_{-\mu+\zeta},\tilde\vp)\subset\cR_{\{1\},k}(1^n_{\mu-\zeta}, \vp)^*$
(see Definition \ref{weih231} for notation and use the dual of (\ref{weih210})).

Note that the target of $(i_U^*)_2$ is the dual of $E_{\infty_P}(1^{n}_{-\lambda+\zeta},\tilde\vp)$.
We define $$(i^*_U)_2: B_U(1^{n+1}_{\lambda }, \vp) \rightarrow
\cR_U(1^n_{ \lambda-\zeta},\vp)$$
by the condition that
\begin{equation}\label{weih404}
\langle\phi,(i^*_U)_2(f)\rangle=\langle i_*(\phi), L\circ [Q]\circ AS(f)\rangle
\end{equation}
for all $\phi\in E_{\infty_P}(1^{n}_{-\lambda+\zeta},\tilde\vp)$, where
$i_*:C^{-\infty}(\partial X^n,V(1^n_{-\lambda+\zeta},\tilde\vp))\rightarrow C^{-\infty}(\partial X^{n+1},V(1^{n+1}_{-\lambda },\tilde\vp))$
is the natural inclusion adjoint to $i^*$ (see \ref{weih301}),
 $L$ is the split of  (\ref{weih210}), $[Q]$ is defined in \ref{weih500}, and 
$AS(f)\in  \cR_{\{1\},k}(1^n_{\lambda}, \vp)$.
This formula defines $(i^*_U)_2(f)$ for generic $\lambda\in \aca$ where $[Q]$ is regular.
The meromorphic family is given by 
$$(i^*_U)_2:=(i_*)_{|E_{\infty_P}(1^{n}_{-\lambda+\zeta},\tilde\vp)}^*\circ L\circ  [Q]\circ AS\ .$$

\subsubsection{}

\begin{lem}\label{sucht}
The definition of $(i^*_U)_2$ is independent of the choice of $k$, the split $L$ and $[Q]$.
\end{lem}
\proof
Let $L^\prime\circ [Q^\prime]\circ AS^\prime$ be defined with different choices (assume that $k^\prime\ge k$).
Let $f_\mu\in B_U(1^{n+1}_\mu,\vp)$ be the germ of a holomorphic family near $\lambda$, and consider a family  $\phi_\mu\in E_{\infty_P}(1^{n}_{-\mu+\zeta},\tilde\vp)$.
Then we have
\begin{eqnarray*}
\lefteqn{AS\circ \pi_*^{U,n+1}(L^\prime\circ [Q^\prime]\circ AS^\prime(f)- L\circ [Q]\circ AS(f))}&&\\&\stackrel{Cor. \ref{u76}}{=}&
[\pi_*^{U,n+1}]\circ AS (L^\prime\circ [Q^\prime]\circ AS^\prime(f)-L\circ [Q]\circ AS(f))\\
&=&[\pi_*^{U,n+1}]  ([Q^\prime]\circ AS^\prime(f)-[Q]\circ L\circ AS(f))\\
&=&AS^\prime(f)-AS(f)\ .
\end{eqnarray*}
 We define (see Lemma \ref{schwspl} for $\{Q\}$)
\begin{eqnarray*}
\Delta&:=&L^\prime\circ [Q^\prime]\circ AS^\prime(f)+\{Q\}\circ L(AS^\prime(f)-AS(f))\\&-&  L\circ [Q]\circ AS(f)\ .
\end{eqnarray*}
Then by construction we have
\begin{equation}\label{weih402}
\pi_*^{U,n+1}(\Delta)=0\ .
\end{equation}
Using results which we will prove (independently of the present stuff) in Subsection \ref{samel}
we show that $\langle i_*(\phi),\Delta\rangle=0$. Proposition \ref{maertins} states an equality
of two spaces defined in \ref{weih400} and \ref{weih401}:
$$\Ext_U(1^{n+1}_{-\mu},\tilde\vp)=E_U(1^{n+1}_{-\mu},\tilde\vp)\ .$$
On the one hand,  $$E_U(1^{n+1}_{-\mu},\tilde\vp)\subset {}^UC^{-\infty}(\partial X^{n+1},V(1^{n+1}_{-\mu},\tilde \vp))$$ is the space of evaluations at $\mu$ of holomorphic families of 
$U$-invariant distributions. In particular we have $i_*(\phi_\mu)\in E_U(1^{n+1}_{-\mu},\tilde\vp)$.
On the other hand, for generic $\mu$ the space $\Ext_U(1^{n+1}_{-\mu},\tilde\vp)$ is contained in the annihilator of the kernel of $\pi^{U,n+1}_*$. Using (\ref{weih402}) we obtain
$\langle i_*(\phi_\mu),\Delta_\mu\rangle=0$.

By our choice of $k$ we have
$\langle i_*(\phi_\mu),g\rangle=0$ for every $g\in \cS_{\{1\},k-2\rho_U}(1^{n+1}_{\mu}, \vp)$ and $\mu$ near $\lambda$
since $g$ vanishes at $\infty_P$ with order larger than $k$.
Note that
$L(AS^\prime(f_\mu)-AS(f_\mu))\in \cS_{U,k}(1^{n+1}_{\mu}, \vp)$.
Since  $\{Q\}\circ L(AS^\prime(f_\mu)-AS(f_\mu))\in \cS_{\{1\},k-2\rho_U}(1^{n+1}_{\mu}, \vp)$
 it follows that
$$\langle i_*(\phi_\mu), \{Q\}\circ L(AS^\prime(f_\mu)-AS(f_\mu))\rangle=0$$
for generic $\mu$ near $\lambda$.
Finally we conclude for these $¸\mu$ that 
\begin{eqnarray*}
\lefteqn{
\langle i_*(\phi_\mu), L^\prime\circ [Q]^\prime\circ AS^\prime(f_\mu)\rangle}&&\\&=&
\langle i_*(\phi_\mu), L^\prime\circ [Q^\prime]\circ AS^\prime(f_\mu)+\{Q\}\circ L(AS^\prime(f_\mu)-AS(f_\mu))\rangle\\
&=&\langle i_*(\phi_\mu), L\circ [Q]\circ AS(f_\mu)\rangle\ .\end{eqnarray*}
\hB

\subsubsection{}

\begin{lem}
The equality (\ref{tobsucht}) holds true.
\end{lem}
 \proof 
Let $f\in  B_{\{1\}}(1^{n+1}_{\lambda }, \vp)$.
Then we have
\begin{eqnarray}\lefteqn{
\pi_*^{U,n+1}(f-L\circ [Q]\circ AS\circ \pi^{U,n+1}_*(f)}&&\nonumber\\&&-\{Q\}(\pi^{U,n+1}_*(f)- \pi^{U,n+1}_*\circ L\circ [Q]\circ AS\circ \pi^{U,n+1}_*(f))\, =\,0\ .\label{tob}
\end{eqnarray}
We have seen in the proof of Lemma \ref{sucht} that $i_*(\phi)$ annihilates the kernel of $\pi^{U,n+1}_*$
(for generic $\lambda$, where all the maps are regular). We combine this fact with (\ref{tob}) in order to derive Equality (\ref{entscheident}) below.
For all  $\phi\in E_{\infty_P}(1^{n}_{-\lambda+\zeta},\tilde\vp)$ we have  
 \begin{eqnarray}
\langle\phi,    AS\circ i^*_U\circ \pi^{U,n+1}_*(f)\rangle&\stackrel{def}{=}&
\langle\phi,    (i^*_U)_2 \circ \pi_*^{U,n+1}(f)\rangle\nonumber\\
&\stackrel{(\ref{weih404})}{=}&\langle i_*(\phi),L\circ [Q]\circ AS\circ \pi^{U,n+1}_*(f)\rangle  \nonumber\\
&=&\langle i_*(\phi),L\circ [Q]\circ AS\circ \pi^{U,n+1}_*(f)\rangle  \nonumber\\
&&+\{Q\}(\pi^{U,n+1}_*(f)-\pi^{U,n+1}_*\circ L\circ [Q]\circ AS\circ \pi^{U,n+1}_*(f))\rangle\nonumber\\
&=&\langle i_*(\phi),f \rangle   \label{entscheident}\\
&=&\langle \phi,i^*(f) \rangle  \nonumber\\
&=&\langle \phi,[\pi^{U,n}_*]\circ AS\circ i^*(f)  \rangle\nonumber\ .
\end{eqnarray}
Since clearly $\{\pi^{U,n}_*\}\circ i^*=(i_U^*)_1\circ  \{\pi^{U,n+1}_*\}$
we conclude the required indentity (\ref{tobsucht}). \hB

\subsubsection{}

Let us combine the results of the present subsection into one statement.
Let now $U\subset P^n$ define a cusp of arbitrary rank.
Note that $(i_U^*)_2$ and therefore $i^*_U$ may have poles in the case of a 
cusp of full rank. The following Proposition settles Proposition \ref{rrttee} in the case of pure cusps.

\begin{prop}\label{klopp}
We have  the following commutative diagram (to be understood as an identity of meromorphic families if $i_U^*$ has poles)
\begin{equation}\label{weih502}
\begin{array}{ccc} 
C^\infty(\partial X,V(1^{n+1}_{\lambda },\vp)) &\stackrel{i^*}{\rightarrow}&C^\infty(\partial X^n, V(1^n_{\lambda-\zeta},\vp))\\
 \downarrow \pi^{U,n+1}_*&&\downarrow \pi^{U,n}_*\\
 B_U(1^{n+1}_{\lambda },\vp)&\stackrel{i_U^*}{\rightarrow }& B_U(1^n_{\lambda-\zeta},\vp)
\end{array}\ .
\end{equation}
 \end{prop}

\subsection{Extension and restriction}\label{samel}

\subsubsection{}

In the present subsection we assume for simplicity that the twist
$\vp$ is normalized such that all its highest $A$-weights are zero.
This differs from the convention adopted in \ref{weih213}. Our present convention has the effect that
the dual $\tilde \vp$ has the normalization adopted in \ref{weih213}.

 Let $l_\vp$ be the highest weight of $\tilde{\vp}$,
i.e. $l_\vp:=l_{\tilde{\vp}}$.  

\subsubsection{}
The space $C^{-\infty}(\infty_P,V(\sigma_\lambda,\vp))$
carries an action $\pi^{\sigma_\lambda,\vp}$ of $P_UA$.
Using the isomorphism of $P_UA$-modules $$C^{-\infty}(\infty_P,V(\sigma_\lambda,\vp))\cong
(\cU(\gaaa)\otimes_{\cU(\paaa)} V_{\sigma_{\lambda+2\rho}})\otimes V_\vp$$
and of $M_UA$-modules (compare \ref{neuj100})
$$(\cU(\gaaa)\otimes_{\cU(\paaa)} V_{\sigma_{\lambda+2\rho}})\otimes V_\vp
\cong \cU(\bar\naaa)\otimes V_{\sigma_{\lambda+2\rho}}\otimes V_\vp$$
given by the PBW-theorem we see that $A$ acts semisimply.

Let $C^{-\infty}(\infty_P,V(\sigma,\vp))^n$ denote
the subspace on which $A$ acts with weight $-\lambda-\rho-n\alpha$.
Then we have $$C^{-\infty}(\infty_P,V(\sigma_\lambda,\vp))^n=
(R_{\{1\}}(\tilde\sigma_{-\lambda},\tilde\vp)^n)^*\ .$$

\subsubsection{}

Since the action of $P_U$ is algebraic
and $U\subset P_U$ is Zariski dense we have
$${}^U C^{-\infty}(\infty_P,V(\sigma_\lambda,\vp))=
{}^{P_U} C^{-\infty}(\infty_P,V(\sigma_\lambda,\vp))\ .$$
Since $A$ normalizes $P_U$ we
conclude that ${}^U C^{-\infty}(\infty_P,V(\sigma_\lambda,\vp))$
is an $A$-invariant subspace.
In particular,  we obtain a decomposition
$${}^U C^{-\infty}(\infty_P,V(\sigma_\lambda,\vp))=\bigoplus_{n=0}^\infty
{}^U C^{-\infty}(\infty_P,V(\sigma_\lambda,\vp))^n= \bigoplus_{n=0}^\infty{}^U (R_{\{1\}}(\tilde\sigma_{-\lambda},\tilde\vp)^n)^*\ .$$

\subsubsection{}\label{neuj703}

Let $\cE_{\infty_P}(\sigma,\vp)^n$ be the sheaf
of holomorphic families $f_\nu\in {}^U (R_{\{1\}}(\tilde\sigma_{-\nu},\tilde\vp)^n)^*$.
This sheaf is torsion-free,  and it is therefore the sheaf of holomorphic
sections of a holomorphic vector bundle $E_{\infty_P}(\sigma,\vp)^n$.
By $E_{\infty_P}(\sigma_\lambda,\vp)^n$ we denote its fibre at $\lambda$.
For each $n\in\nat_0$ 
we define the space  
$\bar Q_(\sigma_\lambda,\vp)^n$ by the exact sequence
$$0\rightarrow E_{\infty_P}(\sigma_\lambda,\vp)^n\rightarrow 
{}^U(R_{\{1\}}(\tilde\sigma_{-\lambda},\tilde\vp)^n)^*\rightarrow
\bar Q_{\infty_P}(\sigma_\lambda,\vp)^n\rightarrow 0\ .$$
Furthermore let $E_{\infty_P}(\sigma_\lambda,\vp):=\bigoplus_{n=0}^\infty E_{\infty_P}(\sigma_\lambda,\vp)^n$
and $\bar Q_{\infty_P}(\sigma_\lambda,\vp):=\bigoplus_{n=0}^\infty \bar Q_{\infty_P}(\sigma_\lambda,\vp)^n$.

\begin{ddd}\label{neuj107}
We call the elements of $E_{\infty_P}(\sigma_\lambda,\vp)^n$ deformable. An element of\linebreak[4] 
${}^U(R_{\{1\}}(\tilde\sigma_{-\lambda},\tilde\vp)^n)^*$ is called undeformable, if it represents a nontrivial class in $\bar Q_{\infty_P}(\sigma_\lambda,\vp)^n$.
\end{ddd}

\subsubsection{}

Let $\taaa\subset \maaa$ be a Cartan subalgebra of $\maaa$. Then
$\haaa:=\taaa \oplus \aaaa$ is a Cartan subalgebra of $\gaaa$.  
By $\Delta^+(\gaaa,\haaa)$ we denote a positive root system
which is compatible with the orientation of $\aaaa$. 
By $\Delta^+(\maaa,\haaa)\subset \Delta^+(\gaaa,\haaa)$
we denote the subsystem of roots of $\maaa$.
For each $\sigma\in \hat{M}$ we define
\begin{equation}\label{neuj103}
\aaaa^*\ni d(\sigma):=-\rho+\max\{
\frac{\langle \mu_\sigma,\varepsilon\rangle}{\langle \alpha,\varepsilon\rangle}
\:|\:\varepsilon\in \Delta^+(\gaaa,\haaa)\setminus \Delta^+(\maaa,\haaa)\} \alpha \ ,
\end{equation}
where $\mu_\sigma$ is the highest weight of $\sigma$, and
$\langle.,.\rangle$ is any Weyl-invariant scalar product on
$\haaa$.
There is a natural action of $P_UA$ on the space
$\Pol(N,V_{\sigma^w_{-\lambda}}\otimes V_\vp)$ of polynomials on $N$ with values in $V_{\sigma^w_{-\lambda}}\otimes V_\vp$ (compare  \ref{neuj101})
given by
$$(man.f)(x) = (\sigma^w_{-\lambda}\otimes \vp)(man) f((n^{-1}x)^{m^{-1}a^{-1}}),\quad
m\in M_U,a\in A, n\in N_V\ .$$

 \subsubsection{}
\begin{lem}\label{vermaolbrich}
\begin{enumerate}
\item 
There is a holomorphic family of $A$-equivariant 
maps 
$$j_\lambda : E_{\infty_P}(\sigma_\lambda,\vp){\rightarrow}\  {}^{P_U}\Pol(N,V_{\sigma^w_{-\lambda}}\otimes V_\vp) \ .$$
\item
If $\Ree(\lambda)> d(\sigma)$ or $\lambda\not\in I_\aaaa$,
then $j_\lambda$ is an isomorphism and $\bar Q_{\infty_P}(\sigma_\lambda,\vp)=0$.
\end{enumerate}
\end{lem}
\proof
\subsubsection{}\label{gaa100}
Let 
$$\hat{J}^w_\lambda:C^{-\infty}(\partial X,V(\sigma_\lambda))
\rightarrow C^{-\infty}(\partial X,V(\sigma^w_{-\lambda}))$$
be the unnormalized Knapp-Stein intertwining operator
(compare \cite{MR1749869}, Sec. 5, (15)).
In order to fix the conventions we recall its definition.
The restriction of  $\hat{J}^w_\lambda$ to smooth sections is given for
$\Ree(\lambda)<0$ by
$$(\hat J^w_{\lambda})f(g)=\int_{\bar N} f(gw\bar n) d\bar n\ .$$
For the rest of parameters it is defined by meromorphic continuation, and it extends by continuity
to distributions.

\subsubsection{}

By  $$j_\lambda:C^{-\infty}(\infty_P,V(\sigma_\lambda))
\rightarrow  \Pol(N,V_{\sigma^w_{-\lambda}})$$
we denote the off-diagonal part of the Knapp-Stein intertwining operator.
Here we identify 
$\Pol(N,V_{\sigma^w_{-\lambda}})$ with a subspace
of $C^\infty(\Omega_P,V(\sigma^w_{-\lambda}))$
such that $p\in \Pol(N,V_{\sigma^w_{-\lambda}})$
corresponds to $f_p\in C^\infty(\Omega_P,V(\sigma^w_{-\lambda}))$
with $f_p(xw)=p(x)$. 

The off-diagonal part of $\hat{J}^w_\lambda$ maps to polynomials since it is $P$-equivariant and the elements of $C^{-\infty}(\infty_P,V(\sigma_{\lambda}))$ are $P$-finite.
Alternatively, using the identification
$$C^{-\infty}(\infty_P,V(\sigma_{\lambda}))\cong\linebreak[4] \cU(\gaaa)\otimes_{\cU(\paaa)} V_{\sigma_{\lambda+2\rho}}\ ,$$
we can write 
$$j_\lambda(X\otimes v)(n)=\pi^{\sigma^w_{-\lambda}}(X) f_{1_v}(nw)\ ,$$
where $1_v\in \Pol(N,V_{\sigma^w_{-\lambda}})$ is the constant
polynomial with value  $v\in V_{\sigma^w_{-\lambda}}$.

\subsubsection{}

The map $j_\lambda$ is in fact $\gaaa$-equivariant, where
the action $\pi^{\sigma^w_{-\lambda}}$ of $\gaaa$
on $\Pol(N,V_{\sigma^w_{-\lambda}})$ is induced by the
embedding $\Pol(N,V_{\sigma^w_{-\lambda}})\subset C^\infty(\Omega_P,V(\sigma^w_{-\lambda}))$.
Therefore $\ker(j_\lambda)$ is a $\gaaa$-submodule
of the Verma module $\cU(\gaaa)\otimes_{\cU(\paaa)} V_{\sigma_{\lambda+2\rho}}$.
By \cite{MR552943}, Satz 1.17,  this Verma module is irreducible
for $\Ree(\lambda)> d(\sigma)$. If $\lambda\not\in I_\aaaa$,
then $j_\lambda$ is injective by \cite{MR1749869}, Lemma 6.7.

\subsubsection{}

Since $$\dim \Pol(N,V_{\sigma^w_{-\lambda}})^n= \dim C^{-\infty}(\infty_P,V(\sigma_\lambda))^n = \dim (\cU(\bar\naaa)\otimes V_{\sigma_\lambda})^n\ ,$$ we conclude that
$j_\lambda$ is in fact an isomorphism.
After tensoring with $V_\vp$  and taking $U$-invariants 
we obtain an isomorphism 
$$j_\lambda:{}^U(R_{\{1\}}(\tilde\sigma_{-\lambda},\tilde\vp)^n)^* \cong {}^{P_U}\Pol(N,V_{\sigma^w_{-\lambda}}\otimes V_\vp)^n\ .$$
Since $j_\lambda$ and its inverse depend holomorphically on $\lambda$
we have $${}^{P_U}\Pol(N,V_{\sigma^w_{-\lambda}}\otimes V_\vp)^n\cong E_{\infty_P}(\sigma_\lambda,\vp)^n$$ for $\lambda\in\aca$
with $\Ree(\lambda) >  d(\sigma)$ or $\lambda\not\in I_\aaaa$.
This proves the lemma. \hB
 
\subsubsection{}

 Recall the definition of the function spaces  $B_{U,k}(\tilde\sigma_{-\lambda},\tilde\vp)$  given in \ref{neuj104}.  

\begin{ddd}\label{neuj1002}
We define
$D_{U,k}(\sigma_\lambda,\vp)$ to be the dual space to
$B_{U,k}(\tilde\sigma_{-\lambda},\tilde\vp)$. Furthermore let
$D_{U}(\sigma_\lambda,\vp):=\bigcup_{k\in\nat_0} D_{U,k}(\sigma_\lambda,\vp)$.
\end{ddd}
As a consequence of the corresponding properties of the family of function spaces
the family of spaces $D_{U,k}(\sigma_\lambda,\vp)$, $\lambda\in\aca$,  forms local trivial holomorphic bundles of dual Fr\'echet spaces. Furthermore, the  spaces
$D_{U}(\sigma_\lambda,\vp)$ are Montel and form a direct limit
of locally trivial holomorphic bundles.

\subsubsection{}

Recall the definition of the push-down Definition \ref{pushdowndef}.
Its meromorphic continuation was finally established in  Corollary \ref{neuj106}.

\begin{ddd}
We define the extension map
$$ext^U:D_{U}(\sigma_\lambda,\vp)\rightarrow D_{\{1\}}(\sigma_\lambda,\vp)$$
as the adjoint of the push-down $$\pi^U_*:B_{\{1\}}(\tilde\sigma_{-\lambda},\tilde \vp)
\rightarrow B_{U}(\tilde\sigma_{-\lambda},\tilde \vp)\ .$$
\end{ddd}
It follows from the corresponding properties of the push-down that
the extension maps form a meromorphic family of
continuous maps with finite-dimensional singularities.

\subsubsection{}

In the remainder of the present subsection we discuss the case
$\sigma=1$.
Assume that the cusp associated to $U\subset P$ has lower rank. 
For $k\in \nat_0$ and $\Ree(\lambda)>-\rho^U+(l_\vp-k\alpha)/2$ such that $$\pi^U_*:B_{\{1\},k}(1_{-\lambda},\tilde \vp)
\rightarrow B_{U,k}(1_{-\lambda},\tilde \vp)$$
is regular (compare Corollary \ref{u76}) we also have a map
 $$ext^U:D_{U,k}(1_\lambda,\vp)\rightarrow D_{\{1\},k}(1_\lambda,\vp)$$
defined as the adjoint $\pi^U_*$. If $k_1 > k$,
then  
$ext^U:D_{U,k}(1_\lambda,\vp)\rightarrow D_{\{1\},k_1}(1_\lambda,\vp)$
is a meromorphic family
of continuous maps with finite-dimensional singularities
defined for
$\Ree(\lambda)>-\rho^U+(l_\vp-k\alpha)/2$.

\subsubsection{}\label{weih400}

\begin{ddd}\label{neuj108}
We define $$Ext_{U}(1_\lambda,\vp)\subset D_{\{1\}}(1_\lambda,\vp)$$ to
be the subspace of all $f\in D_{\{1\}}(1_\lambda,\vp)$ of the form
$ext^U(h)_\lambda$, where $h_\mu\in D_U(1_\mu,\vp)$ is
a meromorphic family defined near $\lambda$ such that $\mu\mapsto ext^U(h)_\mu$ is
regular at $\mu=\lambda$. In a similar manner we define
$Ext_{U,k}(1_\lambda,\vp)$ to by requiring in addition that $h_\mu\in D_{U,k}(1_\mu,\vp)$.
\end{ddd}
The subspaces $Ext_{U,k}(1_\lambda,\vp)\subset Ext_{U}(1_\lambda,\vp)$
are defined for $\Ree(\lambda)>-\rho^U+(l_\vp-k\alpha)/2$.
The space $Ext_{U}(1_\lambda,\vp)$ plays the role of the range of the extension.

\subsubsection{}\label{neuj471}

It is clear that $$Ext_{U}(1_\lambda,\vp)\subset  {}^U D_{\{1\}}(1_\lambda,\vp)\ .$$
In order to describe to which extent the space ${}^U D_{\{1\}}(1_\lambda,\vp)$ is exhausted by $Ext_{U}(1_\lambda,\vp)$ 
we define $Q_{U}(1_\lambda,\vp)$ to be the following quotient:
$$0\rightarrow Ext_{U}(1_\lambda,\vp) \rightarrow {}^U D_{\{1\}}(1_\lambda,\vp)
\rightarrow Q_{U}(1_\lambda,\vp)\rightarrow 0\ .$$
The elements in the space $Q_{U}(1_\lambda,\vp)$ turn out to be somewhat uncontrollable.
We therefore take much effort to show that this space is trivial under certain
conditions.

\subsubsection{}\label{weih401}

We now define the space of deformable $U$-invariant distributions (compare
Def.~\ref{neuj107} for a similar definition with an additional support condition).

\begin{ddd}\label{neuj1088}
We define the subspace $$E_U(1_\lambda,\vp)\subset {}^UD_{\{1\}}(1_\lambda,\vp)$$
as the space of evaluations of germs at $\lambda$ of holomorphic families $f_\nu\in {}^UD_{\{1\}}(1_\nu,\vp)$. In a similar manner we define
$E_{U,k}(1_\lambda,\vp)\subset {}^UD_{\{1\},k}(1_\lambda,\vp)$ as the subspace of evaluations of families with the additional property that $f_\nu\in{}^UD_{\{1\},k}(1_\nu,\vp)$.
\end{ddd}

\subsubsection{}

We define the space $\bar Q_{U}(1_\lambda,\vp)$ as the quotient
$$0\rightarrow E_{U}(1_\lambda,\vp) \rightarrow {}^U D_{\{1\}}(1_\lambda,\vp)
\rightarrow \bar Q_{U}(1_\lambda,\vp)\rightarrow 0\ .$$
\begin{ddd}
An element of ${}^U D_{\{1\}}(1_\lambda,\vp)$
which represents a non-trivial class in $\bar Q_{U}(1_\lambda,\vp)$
is called undeformable.
\end{ddd}
(compare Definition \ref{neuj107}).

\subsubsection{}

It follows immediately from the definitions that
$$Ext_{U}(1_\lambda,\vp) \subset E_U(1_\lambda,\vp)\ .$$
Hence we have a surjection
$$Q_U(1_\lambda,\vp)\to \bar Q_U(1_\lambda,\vp)\ .$$
One of the goals of the present subsection is to show
that
$$Ext_U(1_\lambda,\vp)=E_U(1_\lambda,\vp)\ .$$
Furthermore, we want to show that for many (generic) $\lambda\in \aca$
every element of ${}^U D_{\{1\}}(1_\lambda,\vp)$
is deformable, i.e. $Q_{U}(1_\lambda,\vp)\cong 0$.
The results are stated in Proposition \ref{maertins}.

\subsubsection{}

Now we come to the definition of a left-inverse of $ext^U$: the restriction map $res^U$.
We fix $k\in \nat_0$.
Recall the construction of the meromorphic family of right-inverses of $[\pi_*^U]$
$$[Q]:R_{U,k}(1_\lambda,\tilde{\vp})\to R_{\{1\},k}(1_\lambda,\tilde \vp)$$
from \ref{weih500}. Let $k_1\in\nat_0$ be such that $k_1 \alpha<k\alpha-2\rho_U$, and let
$$\{Q\}:\cS_{U,k}(1_\lambda,\vp)\rightarrow \cS_{\{1\},k_1}(1_\lambda,\vp)$$
be as in Lemma \ref{schwspl}.
We define a meromorphic family of maps
$$Q:B_{U,k}(1_\lambda,\tilde{\vp})\rightarrow B_{\{1\},k_1}(1_\lambda,\tilde{\vp})$$
by
$$Q(h):=  \{Q\}\left(h - \pi^U_*\circ L\circ [Q]\circ AS(h)\right) + L\circ [Q]\circ AS(h)$$
Then one easily checks that $\pi^U_*\circ Q$ is just the inclusion
$B_{U,k}(1_\lambda,\tilde{\vp})\rightarrow B_{U,k_1}(1_\lambda,\tilde{\vp})$.

\begin{ddd}\label{neuj707}
We define the meromorphic family of restriction maps
$$res^U: D_{\{1\},k_1}(1_\lambda,\vp)\rightarrow  D_{U,k}(1_\lambda,\vp)$$
as the adjoint of $Q$.
\end{ddd}
 
Note that $res^U$ depends on choices (we do not indicate these choices in the notation for the
restriction map).  
 In particular, these maps
are not compatible if we change $k$ and $k_1$.

\subsubsection{}\label{neuj2001}

 Note that
the composition $res^U\circ ext^U$ coincides with the inclusion
$$D_{U,k_1}(1_\lambda,\vp) \hookrightarrow      D_{U,k}(1_\lambda,\vp)\ .$$
Recall the Definition \ref{neuj108} of $Ext_{U,k_1}(1_\lambda,\vp)$.
\begin{lem}\label{weedef}
Let $k,k_1\in\nat_0$ and $\lambda\in\aca$ be such that
$res^U:D_{\{1\},k_1}(1_\lambda,\vp)\rightarrow  D_{U,k}(1_\lambda,\vp)$ is defined.
Furthermore we assume that
$\pi^U_*:B_{\{1\}}(1_{-\lambda},\tilde\vp)\rightarrow
B_U(1_{-\lambda},\tilde\vp)$ is regular. Then the restriction
of $res^U$ to $Ext_{U,k_1}(1_\lambda,\vp)$ is independent of choices.
\end{lem}
\proof This follows from $res^U\circ ext^U(f)=f$ and the fact that
$ext^U$ is regular at $\lambda$.
\hB

\subsubsection{}\label{neuj300}


Assume that $\lambda\in\aca$ is such that
$$\pi^U_*:B(1_{-\lambda},\tilde\vp)\rightarrow
B_U(1_{-\lambda},\tilde\vp)$$ and $$[Q]:R_{U}(1_{-\lambda},\tilde{\vp})^n\to R_{\{1\}}(1_{-\lambda},\tilde \vp)^n$$ are regular for all $n\in \nat_0$.
We call $\lambda\in\aca$ satisfying these conditions admissible.

\subsubsection{}\label{neuj704}
Note that  $\cR_{\{1\}}(1_{-\lambda},\tilde\vp)^*$
is the space of those distribution sections of $V(1_\lambda,\vp)$
that are supported at the point $\infty_P$.
We define $$Ext_{\infty_P}(1_\lambda,\vp):=\cR_{\{1\}}(1_{-\lambda},\tilde\vp)^*\cap Ext_U(1_\lambda,\vp)\ .$$
Recall the definition \ref{neuj471} of the space
$Q_{U}(1_\lambda,\vp)$. The point of the following lemma is
that every element of $Q_{U}(1_\lambda,\vp)$ can be represented by 
an invariant distribution supported in $\infty_P$.
Note that $$[ext^{U}]:\cR_{U}(1_{-\lambda},\tilde\vp)^*\rightarrow {}^U \cR_{\{1\}}(1_{-\lambda},\tilde\vp)^*$$
generates a subspace of $Ext_{\infty_P}(1_\lambda,\vp)$, where $[ext^{U}]$ is the restriction of $ext^{U}$
to $\cR_{U}(1_{-\lambda},\tilde\vp)^*$, or equivalently, the adjoint
of $[\pi^{U}_*]$ (see Def.~\ref{weih254}).
\begin{lem}\label{homjk}
\begin{enumerate}
\item 
The family of maps  $[ext^{U}]$ generates all of $Ext_{\infty_P}(1_\lambda,\vp)$.
\item
If $\lambda\in\aca$ is admissible in the sense of \ref{neuj300}, then
there is an exact sequence
$$0\rightarrow
Ext_{\infty_P}(1_\lambda,\vp)\rightarrow
{}^{U}\cR_{\{1\}}(1_{-\lambda},\tilde\vp)^* \rightarrow
Q_{U}(1_\lambda,\vp)\rightarrow 0$$
of semisimple $A$-modules.\end{enumerate}
\end{lem}
\proof
Recall that $${}^{U}\cR_{\{1\}}(1_{-\lambda},\tilde\vp)^*=
\bigoplus_{n\in\nat_0} {}^U(\cR_{\{1\}}(1_{-\lambda},\tilde\vp)^n)^*$$
is a weight-decomposition of the $A$-module.
Moreover, we have
$$[ext^{U}]:(\cR_{U}(1_{-\lambda},\tilde\vp)^n)^*\rightarrow (\cR_{\{1\}}(1_{-\lambda},\tilde\vp)^n)^*\ .$$  

If $f\in Ext_{\infty_P}(1_\lambda,\vp)$
is represented by $ext^{U}(h)$ for some meromorphic family
$h_\mu\in D_{U}(1_\mu,\vp)$, then
$\{res^{U}\} \{ext^{U}(h)\} =\{h\}$ vanishes at $\lambda$,
where $\{h\}$ denotes the restriction
of $h$ to $\cS_{\{1\}}(1_{-\lambda},\tilde\vp)$, and $\{res^U\}:=\{Q\}^*$.
Therefore all non-positive Laurent-coefficients of
the expansion of $h$ at $\lambda$
belong to $\cR_{U}(1_{-\lambda},\tilde\vp)^*$.
Thus we can choose the family $h$ such that
$h_\mu\in\cR_{U}(1_{-\mu},\tilde\vp)^*$, and such that $[ext^{U}](h)_\lambda=f$.
This shows that $[ext^{U}]$ generates $Ext_{\infty_P}(1_\lambda,\vp)$,
and that $A$ acts semisimply on 
$Ext_{\infty_P}(1_\lambda,\vp)$.

\subsubsection{}

It remains to show that any element of $Q_U(1_\lambda,\vp)$
can be represented by some element of ${}^{U}\cR_{\{1\}}(1_{-\lambda},\tilde\vp)^*$. 
Let $f\in {}^{U}D_{\{1\}}(1_\lambda,\vp)$.
Then there is $k\in\nat_0$ such that $f\in {}^{U}D_{\{1\},k}(1_\lambda,\vp)$.
We choose $k_1$ such that $$\{res^{U}\}:\cS_{\{1\},k}(1_{-\lambda},\tilde\vp)^*
\rightarrow  \cS_{U,k_1}(1_{-\lambda},\tilde\vp)^*$$ is defined.
We then put $h:=\{res^{U}\}\{f\}$. Let $$T:\cS_{U,k_1}(1_{-\lambda},\tilde\vp)^*
\rightarrow D_{U}(1_\lambda,\vp)$$ be the split
induced by the dual split $L$. We
form $f-ext^{U}\circ T(h)$. This difference represents the same
element in $Q_{U}(1_\lambda,\vp)$ as $f$, but its restriction to
$\cS_{\{1\}}(1_{-\lambda},\tilde\vp)$ vanishes. Indeed, for $g$ in $\cS_{\{1\}}(1_{-\lambda},\tilde\vp)$
we have 
$$ \langle ext^{U}\circ T(h),g\rangle= \langle f, \{Q\}\circ\{\pi_*^U\}(g)\rangle=\langle f, g\rangle
\ .$$ \hB

\newcommand{\pr}{{\mathrm{pr}}}
\newcommand{\F}{{\mathbb{F}}}


\subsubsection{}

\begin{prop}\label{argumentprop}
Let $n\ge 0$. If $\lambda\in\aaaa^*$ is sufficiently large, then the inclusion
$$Ext_{\infty_P}(1_\lambda,\vp)^n\hookrightarrow {}^U(R_{\{1\}}(1_{-\lambda},\tilde\vp)^n)^*$$ is an isomorphism.
\end{prop}
\proof
If $U$ defines a cusp of full rank, then the proposition is a direct consequence of Lemma \ref{vermaolbrich}
and the definition of $[ext^U]$ (Definition \ref{rolf}). Thus we may and will assume in the following
that $U$ defines a cusp of smaller rank.
Note that both spaces appearing in the proposition are finite-dimensional.
The proposition is an immediate consequence of the following lemma.
\begin{lem}\label{argumentlem1}
If $\lambda\in\aaaa^*$ is sufficiently large, then we have the inequality
$$\dim Ext_{\infty_P}(1_\lambda,\vp)^n\ge \dim  {}^U(R_{\{1\}}(1_{-\lambda},\tilde\vp)^n)^*\ .$$
\end{lem}
\proof
By Lemma \ref{vermaolbrich},2, if $\lambda\in\aaaa^*$ is sufficiently large we have
$\dim  {}^U(R_{\{1\}}(1_{-\lambda},\tilde\vp)^n)^*=\dim {}^{P_U} \Pol(N,V_{1_{-\lambda}}\otimes V_\vp)^n$, where
the action of $AP_U$ on a polynomial is given by 
$$(p \: f)(n):=\vp(p)f(n^{p^{-1}})\ .$$
The degree-$n$ subspace of the polynomial maps is characterized by
\begin{equation}\label{argumentdegree1}a\:f=a^{-n\alpha} f\ ,\quad  a\in A\ .\end{equation}
We now consider
$$[\pi^U_*]:R_1(1_{-\lambda},\tilde \vp)^n\to R_U(1_{-\lambda},\tilde \vp)^n$$
which is regular for $\Ree(\lambda)$ sufficiently large. Its adjoint is
$$[ext^U]:(R_U(1_{-\lambda},\tilde \vp)^*)^n\to Ext_{\infty_P}(1_\lambda,\vp)^n\ ,$$
and we have 
$\dim \im [\pi^U_*] =\dim \im [ext^U]$.

\subsubsection{}

Therefore,
 Lemma \ref{argumentlem1}  follows directly  from the following lemma.
\begin{lem}\label{argumentlem2}
If $\lambda\in\aaaa^*$ is sufficiently large, then we have
$$\dim \im [\pi^U_*] \ge \dim {}^{P_U} \Pol(N,V_{1_{-\lambda}}\otimes V_\vp)^n\ .$$
\end{lem}
\proof
Let $\langle.,.\rangle$ be some non-degenerate invariant bilinear form on $\gaaa$.
We define the linear subspace
\begin{equation}\label{argumenteq33}\bar \naaa^U:=\{Y\in \bar \naaa\:|\: \langle [Y,H],X\rangle=0\:\:\forall X\in \naaa_V\}\end{equation}
and the submanifold $\bar N^U:=\exp(\bar \naaa^U)\subset \bar N$.
We furthermore choose a Cartan involution $\theta$ of $\gaaa$ compatible with $\aaaa$. 
\begin{lem}\label{argumentlem6}
\begin{enumerate}
\item The submanifold $\bar N^U$ is $AM_U$-invariant.
\item The multiplication map
$\bar N_V\times \bar N^U\to \bar N$ is a diffeomorphism, where $\bar N_V:=N_V^\theta$.
\item
The composition
$$\bar N\stackrel{\sim}{\to}\bar N_V\times \bar N^U\stackrel{\pr}{\to} \bar N^U$$
is an $AM_U$-equivariant polynomial map. 
\end{enumerate}
\end{lem}
\proof
The subspace $\bar \naaa^U$ is $AM_U$-invariant since $\naaa_V$ is $AM_U$-invariant.
The first assertion now follows from the $AM_U$-equivariance of the exponential map

By the $A$-invariance of $\bar N^U$, $\bar N_V$ and the equivariance of the multiplication, it suffices to show that the multiplication map is a diffeomorphism near $(1,1)$. Infinitesimally it is given by the map $\bar\naaa_V\times \bar \naaa^U\to \bar \naaa$, $(Y_V,Y^U)\mapsto \bar Y_V+Y^U$. This map is an isomorphism. In fact, the dimensions of the domain and the target coincide, and $\bar \naaa_V\cap \bar \naaa^U=\{0\}$.

The last assertion follows from the diagram
$$\xymatrix{\bar \naaa\ar[d]^{\exp}&\bar \naaa_V\times\bar \naaa^U\ar[l]_{q\quad}\ar[d]^{(\exp,\exp)}\ar[r]^{\quad pr}&\bar \naaa^U\ar[d]^{\exp}\\\bar N&\bar N_V\times \bar N^U\ar[l]_{mult\quad}\ar[r]^{\quad\pr}&\bar N^U}$$
and the fact that $q$ is an $AM_U$-equivariant polynomial map with a polynomial inverse.
\hB

\subsubsection{}

Note that $\bar N_VAM_U$ acts on $\Pol(\bar N,V_{1_{-\lambda}}\otimes V_{\tilde\vp})$
by
$$\bar n_Vam\: f(\bar n)=\tilde\vp(am) f((\bar n_V^{-1}\bar n)^{(am)^{-1}})\ .$$
We now consider the space
$$I(\lambda)^n:={}^{\bar N_VM_U}\Pol(\bar N,V_{1_{-\lambda}}\otimes V_{ \tilde \vp})^n\subset\Pol(\bar N,V_{1_{-\lambda}}\otimes V_{\tilde \vp})^n\ ,$$
where the degree $n$-subspace is distinguished by the condition
\begin{equation}\label{argumentdegree2}a\: f=a^{n\alpha} f\ .\end{equation}
We have a degree-preserving inclusion
$$I(\lambda)^n\hookrightarrow R_1(1_{-\lambda},\tilde \vp)^n\ .$$
Furthermore, it follows from Lemma \ref{argumentlem6} that the restriction to $\bar N^U$ induces an isomorphism
$$I(\lambda)^n\stackrel{\sim}{\to} {}^{M_U}\Pol(\bar N^U,1_{-\lambda}\otimes V_{\tilde \vp})^n\ .$$ 
Lemma \ref{argumentlem2} now follows from the following two assertions.

\begin{lem}\label{argumentlem3}
$\dim I(\lambda)^n= \dim {}^{P_U} \Pol(N,V_{1_{-\lambda}}\otimes V_\vp)^n$
\end{lem}
\begin{lem}\label{argumentlem4}
If $\lambda\in\aaaa^*$ is sufficiently large, then the restriction of
$[\pi^U_*]$ to $I(\lambda)^n$ is injective.
\end{lem}

\subsubsection{}
We first show Lemma \ref{argumentlem3}.
We define $N^U:=(\bar N^U)^\theta\subset N$. Then we have an $A$-equivariant diffeomorphism
$N_V\times N^U\stackrel{mult}{\to} N$, and the projection $N\stackrel{\sim}{\to} N_V\times N^U\stackrel{\pr}{\to} N^U$
is an $AM_U$-equivariant polynomial map. These facts follow from Lemma \ref{argumentlem6} by an application of the Cartan involution $\theta$. We conclude that the restriction to $N^U$ induces an isomorphism
$${}^{P_U} \Pol(N,V_{1_{-\lambda}}\otimes V_\vp)^n\stackrel{\sim}{\to} {}^{M_U}\Pol(N^U, V_{1_{-\lambda}}\otimes V_\vp)^n\ .$$
Note that there is a canonical $AM_U$-equivariant isomorphism between $S((\naaa^U)^*)\otimes V_\vp$ and \linebreak[4]
$\Pol(N^U, V_{1_{-\lambda}}\otimes V_\vp)$. Here
$S(.)$ stands for the symmetric algebra. Similarly, we have 
$$S((\bar\naaa^U)^*)\otimes V_{\tilde\vp}\cong \Pol(\bar N^U,V_{1_{-\lambda}}\otimes V_{\tilde \vp})\ .$$
The natural pairing between $\bar\naaa^U$ and $\naaa^U:= (\bar\naaa^U)^\theta$ via the $G$-invariant form $\langle.,.\rangle$ induces a nondegenerate pairing between the two symmetric algebras above. We conclude that the spaces
${}^{M_U}\Pol(N^U, V_{1_{-\lambda}}\otimes V_\vp)^n$
and
$
{}^{M_U}\Pol(\bar N^U,V_{1_{-\lambda}}\otimes V_{\tilde \vp})^n\cong I(\lambda)^n$
are each others duals
(the pairing is degree-preserving in view of the characterizations (\ref{argumentdegree1}) and (\ref{argumentdegree2})).
Lemma \ref{argumentlem3} now follows.

\subsubsection{}

We now start with the proof of Lemma \ref{argumentlem4}.
First of all note that the push-down (see Def.~\ref{weih254})
$$[\pi^U_*]:I(\lambda)^n\to R_U(1_{-\lambda},\tilde \vp)^n\subset A_U(1_{-\lambda},\tilde \vp)^n $$ is given by a convergent integral
$$[\pi^U_*](f)(\bar n)=\int_{N_V}\vp(n_V)^{-1} f(\bar n(n_V\bar n))a(n_V\bar n)^{-\lambda-\rho} dn_V\ ,\bar n\in \bar N\setminus\{1\}\ .$$
Here we have employed again the Bruhat decomposition $g=\bar n(g) m(g) a(g) n(g)$.
The vector space $V_{\tilde \vp}$  has a filtration induced by the action of  $A$ that is preserved  by $M_UN_V$. The induced action of $N_V$ on the associated graded vector space
$\Gr(V_{\tilde \vp})$ is trivial. The filtration of $V_{\tilde \vp}$ induces filtrations on $I(\lambda)^n$ and $A_U(1_{-\lambda},\tilde \vp)^n$. From the integral representation of $[\pi^U_*]$ we see that that this map preserves the filtrations and induces a 
map $\Gr[\pi_U^*]:\Gr I(\lambda)^n\to \Gr A_U(1_{-\lambda},\tilde \vp)^n$.
Since injectivity of the associated graded map implies injectivity of a filtration preserving map, Lemma 
\ref{argumentlem4} is a consequence of 
\begin{lem}\label{argumentlem44}
If $\lambda\in\aaaa^*$ is sufficiently large, then 
$$\Gr [\pi^U_*]:\Gr I(\lambda)^n\to  \Gr A_U(1_{-\lambda},\tilde \vp)^n$$ is injective.
\end{lem}

\subsubsection{}
Note that
$$\Gr [\pi^U_*](f)(\bar n)=\int_{N_V} f(\bar n(n_V\bar n))a(n_V\bar n)^{-\lambda-\rho}dn_V\ ,\bar n\in \bar N\setminus\{1\}\ .$$
We define
$$c_\lambda(\bar n):=\Gr [\pi^U_*](1)(\bar n)= \int_{N_V} a(n_V\bar n)^{-\lambda-\rho}dn_V\ ,\bar n\in \bar N\setminus\{1\}\ .$$
We will show the following lemma.
\begin{lem}\label{argumentlem8}
For $\bar n^U\in \bar N^U\setminus\exp(\bar\naaa_{-2\alpha})$ we have
$$\frac{1}{c_\lambda(\bar n^U)} \Gr [\pi^U_*](f)(\bar n^U)= f(\bar n^U)+O(|\lambda|^{-1})\ .$$
\end{lem}

\subsubsection{}
Let us first show that Lemma \ref{argumentlem8} implies \ref{argumentlem44}.
The natural identification $V_{1_{-\lambda}}\cong \C$ induces identifications
$\Pol(\bar N^U,1_{-\lambda}\otimes V_{\tilde \vp})^n\cong  \Pol(\bar N^U,V_{\tilde \vp})^n$ for all $\lambda$.

Note that $\bar \naaa^U\cap\bar\naaa_{-\alpha}$ is non-trivial. Otherwise we would have 
$\naaa_V\cap\naaa_{\alpha}=\naaa_{\alpha}$ and therefore $\naaa_V=\naaa$, i.e., $U$ would define a cusp of full rank.
Hence we can choose a finite sequence of base points $\bar n^U_i\in \bar N^U\setminus \exp(\bar\naaa_{-2\alpha})$ and vectors
$v_i\in V_\vp$, $i=1,\dots,r:=\dim {}^{M_U}\Pol(\bar N^U,V_{\tilde \vp})^n$ such that
the following map is an isomorphism:
$$\Phi:{}^{M_U}\Pol(\bar N^U,V_{\tilde \vp})^n\to\C^r\ , f\mapsto (\langle v_1,f(n^U_1)\rangle,\dots,\langle v_r,f(n^U_r)\rangle)\ .$$

We now consider the composition
$$A(\lambda):=\Phi\circ  \frac{1}{c_\lambda} \Gr [\pi^U_*]\circ \Phi^{-1}:\C^n\to \C^n\ .$$
Lemma \ref{argumentlem8} implies that
$$A(\lambda)=1+O(|\lambda|^{-1})\ .$$ 
In particular, $A$ is injective, if $\lambda\in\aaaa$ is sufficiently large.
This implies the assertion of Lemma \ref{argumentlem44}.

\subsubsection{}

We now show Lemma \ref{argumentlem8}. 
Fix $\bar n^U\in \bar N^U\setminus \exp(\bar\naaa_{-2\alpha})$. We set 
$$\Psi(n_V):=\log a(n_V\bar n^U) \ ,\quad g(n_V):=f(\bar n(n_V\bar n^U))$$
Then we can write
$$\int_{N_V} f(\bar n(n_V\bar n^U))a(n_V\bar n^U)^{-\lambda-\rho}dn_V=\int_{N_V} \ee^{(-\lambda-\rho)\Psi( n_V)}g(n_V)dn_V\ .$$  
\begin{lem}\label{argumentlem9}
The function $\Psi(n_V)$ has a unique non-degenerate absolute minimum $\Psi(1)=0$ at
$n_V=1$. Furthermore, there exists a compact neighbourhood $K\subset N_V$ of $1$ such that
$\alpha(\Psi(n_V))\ge 1$ for $n_V\not\in K$.
\end{lem}
We first show that Lemma \ref{argumentlem9} implies \ref{argumentlem8}. 
We split the integral
as $\int_{N_V}=\int_{K}+\int_{N_V\setminus K}$. 
Note that $\bar n(N_V\bar n^U)\subset \bar N$ is a pre-compact subset. It follows that
$g$ is smooth and uniformly bounded. We approximate
the first summand by a Gaussian integral at the minimum of $\Psi$ and get
$$\int_K \ee^{(-\lambda-\rho)(\Psi(n_V))}g(n_V)dn_V= \int_K \ee^{(-\lambda-\rho)(\Psi(n_V))}dn_V  \left(g(0)+O(\lambda^{-1})\right)\ .$$ 
Furthermore, 
$\int_K \ee^{(-\lambda-\rho)(\Psi(n_V))}dn_V$ decreases at most as $|\lambda|^{-\dim N_V/2}$.
We claim that the contributions $\int_{N_V\setminus K}$ decrease exponentially so that these parts of the integrals
can only contribute exponentially small error terms. This claim implies Lemma \ref{argumentlem8}. 

In order to see the claim
we write
$$\int_{N_V\setminus K}\ee^{(-\lambda-\rho)(\Psi(n_V))}dn_V =\min_{n_V\in N_V\setminus K}\ee^{(-\lambda+\alpha)(\Psi(n_V))}\int_{N_V\setminus K}\ee^{(-\alpha-\rho)(\Psi(n_V))}dn_V\ .$$
The integral on the right-hand side converges, and
$$\min_{n_V\in N_V\setminus K}\ee^{(-\lambda+\alpha)(\Psi(n_V))}\le \ee^{-c|\lambda|}$$ for a suitable constant $c>0$.
This finishes the proof of Lemma \ref{argumentlem8} under the assumption of Lemma \ref{argumentlem9}.

\subsubsection{}

We now show Lemma \ref{argumentlem9}. In order to compute $a(n_V\bar n^U)$ we may assume
that $G$ is the subgroup of $GL(n+1,\F)$ that preserves the $\F$-valued Hermitian scalar product
\begin{equation}\label{argumenteq45}\bar v_0 w_n+ \bar v_n w_0+\bar v_1 w_1+ \dots + \bar v_{n-1} w_{n-1}\ \end{equation}
on the right $\F$-vector space $\F^{n+1}$.
We choose 
$$A:=\left\{\left(\begin{array}{ccc}a&0&0\\0&1_{n-1\times n-1}&0\\0&0&a^{-1}\end{array}\right)\:|\: a\in \R^+\right\}$$
and get
$$N:=\left\{\left(\begin{array}{ccc}1&w&p-\frac{\|w\|^2}{2}\\0&1&-\bar w^t\\0&0&1\end{array}\right)\:|\: w\in  \F^{n-1}\ ,p\in \im\, \F\right\}\ .$$
The parametrization of $N$ given here is via the exponential map, if we identify
$\naaa=\naaa_\alpha\oplus\naaa_{2\alpha}\cong \F^{n-1}\oplus \Imm\,\F$.
Furthermore,
$$\bar N:=\left\{\left(\begin{array}{ccc}1&0&0\\-\bar v^t&1&0\\q-\frac{\|v\|^2}{2}&v&1\end{array}\right)\:|\: v\in  \F^{n-1}\ ,q\in \Imm\, \F\right\}\ .$$
Since they are $A$-invariant, the subspaces $\naaa_V\subset \naaa$ and $\bar \naaa^U\subset \bar \naaa$ are given in this identification as
$W\oplus P\subset \naaa$ and $V\oplus Q\subset \bar \naaa$ for real subspaces $W,V\subset \F^{n-1}$ and $P,Q\subset \Imm\, \F$. 
We consider the invariant form $\langle A,B\rangle:=\Ree\,\Tr AB$ on $\gaaa$. Given
$W\oplus P$, the space $V\oplus Q$ is characterized by
$\langle\ad(V\oplus W)(H),W\oplus H\rangle=0$ (see (\ref{argumenteq33})).
Explicitly,
$$H=\left(\begin{array}{ccc}1&0&0\\0&0&0\\0&0&-1\end{array}\right)\ ,\langle\ad(v+q)(H),w+p\rangle=-2\Ree (\bar w v^t+qp)\ .$$
Therefore we have
\begin{equation}\label{argumenteq88}V=W^\perp\ , Q=P^\perp\end{equation} with respect to the natural euclidean pairings on $\F^{n-1}$ and $\Imm\, \F$.
The vector $e_0:=(1,\dots,0)\in \F^{n+1}$ is the highest weight vector of the standard representation of $G$ with weight $\alpha$. Similarly, $e_n:=(0,\dots,0,1)$ is the lowest weight vector with weight $-\alpha$.
Let $g=\bar n man\in \bar NMAN$. If $\langle.,.\rangle$ denotes the $\F$-valued scalar product (\ref{argumenteq45}), then
$$\langle ge_0,e_n\rangle=\langle\bar ma n e_0,e_n\rangle=\langle ma n e_0,\bar n^{-1}e_n\rangle=a^\alpha\langle m e_0,e_n\rangle\ .$$ 
The subspace $e_0\F\cong \F$ is invariant under the group $M$. In particular, we have a homomorphism
$\vartheta:M\to \F^*$ such that $m e_0=e_0 \vartheta(m)$. It now follows that
$\langle m e_0,e_n\rangle=\bar \vartheta(m)$. Since $M$ is compact, we have $\|\vartheta(m)\|=1$ and therefore
$a^\alpha=\|\langle ge_0,e_n\rangle\|$.
If we parametrize $(w,p)=n_V\in N_V$ and $(v,q)=\bar n^U\in \bar N^U$ as above, then we get
\begin{equation}\label{argumenteq100} a(n_Vn^U)^{2\alpha}=\|1-w\bar v^t+(p-\frac{\|w\|^2}{2})(q-\frac{\|v\|^2}{2})\|^2\ .\end{equation} 
We fix $(v,q)\in V\oplus Q$ with $v\not=0$.
We must show  that the right-hand side has a unique absolute minimum $1$ at $(w,p)=(0,0)$, and that this minimum is non-degenerate. Using (\ref{argumenteq88}) we get
\begin{eqnarray*}
\Ree(1-w\bar v^t+(p-\frac{\|w\|^2}{2})(q-\frac{\|v\|^2}{2}))&=&1+\frac{\|v\|^2\|w\|^2}{4}\\
\Imm(1-w\bar v^t+(p-\frac{\|w\|^2}{2})(q-\frac{\|v\|^2}{2}))&=&pq-w\bar v^t-\frac{\|w\|^2}{2}q-\frac{\|v\|^2}{2}p\ .
\end{eqnarray*}
First of all, $a(n_Vn^U)^{2\alpha}\ge 1$ and $a(1)^{2\alpha}=1$. Moreover, if
$a(n_V\bar n^U)^{2\alpha}=1$, then $w=0$.
Since $\Ree(\overline{pq}p)=\Ree \bar q\|p\|^2)=0$ we have
$pq\perp \frac{\|v\|^2}{2}p$. Hence, the equality  $a(n_V\bar n^U)=1$ implies in addition to $w=0$ that also $p=0$.
We thus have shown that  $a(n_V\bar n^U)^{2\alpha}$ takes its unique absolute minimum at $1$.
Next we show that it is non-degenerate. The Hessian $h(p,w)$ is the part of the polynomial (\ref{argumenteq100}) which is quadratic in $(w,p)$. It can be written as
$$h(w,p):=\|w\|^2\frac{\|v\|^2}{2}+\|p(q-\frac{\|v\|^2}{2})-w\bar v^t\|^2\ .$$
If $h(w,p)=0$, then from the first summand and $v\not=0$ we get $w=0$ and $p(q-\frac{\|v\|^2}{2})=0$. Since $(q-\frac{\|v\|^2}{2})\not=0$ we conclude that $p=0$.
The last assertion of Lemma \ref{argumentlem9} follows from
$\lim_{n_V\to \infty}a(n_V\bar n^U)^{2\alpha}=\infty$, which is easy to check using the explicit formula  (\ref{argumenteq100}). This finishes the proof of Lemma \ref{argumentlem9}.

We now have also finished the proof of Proposition \ref{argumentprop}.\hB

\subsubsection{}

\begin{prop}\label{maertins}
\begin{enumerate}
\item We have $Ext_{U}(1_\lambda,\vp)=E_{U}(1_\lambda,\vp)$
for all $\lambda\in\aca$.
\item If $\lambda\not\in I_\aaaa$ or $\Ree(\lambda)>-\rho$, then
$Q_U(1_\lambda,\vp)=0$.\end{enumerate}
\end{prop}
\proof
\subsubsection{}
As in\ref{neuj703} let $\cE_{\infty_P}(1,\vp)^n$ be the torsion-free
coherent sheaf on $\aca$ of holomorphic families of
$U$-invariant distributions supported on $\infty_P$,
on which $A$ acts by multiplication by the function
$\lambda\mapsto a^{\lambda-\rho-n\alpha}$.
By $\cE xt_{\infty_P}(1,\vp)^n$ we denote the subsheaf generated by the
restriction of $[ext^{U}]$ to the homogeneous part $(\cR_{U}(1_{-.},\tilde\vp)^n)^*$.

The space  $Ext_{\infty_P}(1_\lambda,\vp)^n$ is the geometric fibre of $\cE xt_{\infty_P}(1,\vp)^n$ for the generic set of $\lambda\in\aca$, where $[ext^U]$ is regular.

\subsubsection{}\label{neuj705}
The sheaf $\cE_{\infty_P}(1,\vp)^n$ is the sheaf of sections of a finite-dimensional trivial holomorphic vector bundle $E_{\infty_P}(1,\vp)^n\to \aca$.
 The torsion-free subsheaf $\cE xt_{\infty_P}(1,\vp)^n$ corresponds to a
bundle $\Ext_{\infty_P}(1,\vp)^n$.

This discussion shows the following: if the inclusion $Ext_{\infty_P}(1_\lambda,\vp)^n\hookrightarrow E_{\infty_P}(1_\lambda,\vp)^n$ is surjective
at one point $\lambda\in \aca$, then it is surjective for generic $\lambda$, i.e. outside a discrete set.
However, we know from Proposition \ref{argumentprop}  that this inclusion is surjective for many $\lambda$, hence it is so generically.

\subsubsection{}\label{miau}

In view of Lemma \ref{homjk}, 2.
 we have
\begin{equation}\label{neuj708}
E_{\infty_P}(1_\lambda,\vp)^n/E xt_{\infty_P}(1_\lambda,\vp)^n\cong \ker\left( Q_U(1_\lambda,\vp)^n\to
\bar Q_{\infty_P}(1_\lambda,\vp)^n\right)\ .
\end{equation}
By \ref{neuj705} the quotient on the left hand side is trivial generically. The same is true
for $\bar Q_{\infty_P}(1_\lambda,\vp)^n$ by Lemma \ref{vermaolbrich}. We conclude that $Q_U(1_\lambda,\vp)^n$ is trivial outside a discrete set.

\subsubsection{}
We now prove the first assertion of Prop.~\ref{maertins}.
We know that $Ext_{U}(1_\lambda,\vp)\subset E_{U}(1_\lambda,\vp)$. 
Let now $f\in E_{U}(1_\lambda,\vp)$ be given 
 as the value at $\lambda$
of a meromorphic family $f_\mu\in {}^{U} D_{\{1\},k}(1_\mu,\vp)$ for some sufficiently large $k$.
Let $res^U$ be the meromorphic family of continuous maps
$res^U:{}^{U} D_{\{1\},k}(1_\mu,\vp):\rightarrow  D_{U,k_1}(1_\mu,\vp)$
for suitable $k_1\in\nat$ (see Definition \ref{neuj707}). By \ref{miau} for generic
$\mu$ we can write $f_\mu=ext^{U} g_\mu$ for some $g_\mu\in
D_{U,k}(1_\mu,\vp)$.
We conclude $ext^{U}\circ res^{U}(f_\mu)= ext^{U}\circ res^{U} \circ ext^U(g_\mu)= ext^U(g_\mu)=f_\mu$  by Lemma \ref{weedef}.
Thus 
\begin{equation}\label{wau}
ext^{U}\circ res^{U}(f_\mu)=f_\mu \quad\mbox{ for all }\mu\ .
\end{equation}
We conclude that $f\in Ext_{U}(1_\lambda,\vp)$.

\subsubsection{}

We now turn to the second assertion. The first assertion implies that the left hand side
of (\ref{neuj708}) is trivial for all $\lambda$. Therefore the map $Q_U(1_\lambda,\vp)^n\to
\bar Q_{\infty_P}(1_\lambda,\vp)^n$ is always injective. We now apply Lemma \ref{vermaolbrich}.
\hB

\subsubsection{}
The argument that proofs Equation (\ref{wau}) also shows the following.
\begin{lem}\label{wauwau}
The composition  $ext^U\circ res^U$
is the identity on $Ext_{U}(1_{\lambda},\vp)$.
\end{lem}

\subsection{The scattering matrix}

\subsubsection{}


Recall that the family of intertwining operators (see \ref{gaa100} for the unnormalized version and \cite{MR1749869} for normalizations) forms a meromorphic family of operators.
It therefore maps (holomorphic) families of invariant sections to (meromorphic) families. 
It follows that 
$J_\lambda$ maps $E_U(1_\lambda,\vp)$ to $E_U(1_{-\lambda},\vp)$
if $\lambda\in \aca$ is such that $J_\lambda$ is regular (e.g. $\lambda\not\in I_\aaaa$, see \ref{gaa101}). By Proposition \ref{maertins} we get a mapping
\begin{equation}\label{uhu}
J_\lambda:Ext_{U}(1_\lambda,\vp)\rightarrow Ext_{U}(1_{-\lambda},\vp)\ .
\end{equation}
Moreover, 
for given $k\in\nat_0$ 
we have
$$J_\lambda:Ext_{U,k}(1_\lambda,\vp)\rightarrow Ext_{U,k_1}(1_{-\lambda},\vp)$$
if $k_1\in\nat_0$ is sufficiently large.

\subsubsection{}\label{neuj2004}

We can now define the scattering matrix
$$S^U_\lambda:D_U(1_\lambda,\vp)\rightarrow D_U(1_{-\lambda},\vp)\ .$$
Fix $k\in\nat_0$ and a compact subset $W\subset \aaaa^*$.
Then we choose $k_0>k$ and $k_1,k_2$ such that
$$J_\lambda: Ext_{U,k_0}(1_\lambda,\vp)\rightarrow Ext_{U,k_1}(1_{-\lambda},\vp)$$
for generic $\lambda$ (e.g.  non-integral)  with $\Ree(\lambda)\in W$ 
and
$$res^U:D_{\{1\},k_1}(1_\lambda,\vp)\rightarrow D_{U,k_2}(1_\lambda,\vp)$$
is defined  (Definition \ref{neuj707}) as a meromorphic family on $W$.
Then we consider the composition
$$S^U_\lambda:=res^U\circ J_\lambda\circ ext^U
:D_{U,k}(1_\lambda,\vp)\rightarrow D_{U,k_2}(1_{-\lambda},\vp)\ .$$
By definition $S^U_\lambda$ is a meromorphic family of continuous maps.
Since
$$J_\lambda\circ ext^U:D_U(1_\lambda,\vp)\to Ext_U(1_{-\lambda},\vp)$$
Lemma \ref{weedef} now implies that  $S^U_\lambda$
is well-defined independently of the choices made for $res^U$.

 \subsubsection{}

Letting $k$ tend to infinity and
$W$ run over a sequence of compact subsets exhausting $\aaaa^*$
we are arrive at the following definition.
\begin{ddd}\label{neuj2003}
We define the scattering matrix as the meromorphic family of continuous maps
$$S^U_\lambda:D_U(1_\lambda,\vp)\rightarrow D_U(1_{-\lambda},\vp)\ ,$$
which is given by  the composition 
$$S^U_\lambda(f):=res^U\circ J_\lambda\circ ext^U(f)\ $$
whenever the constituents are regular. 
\end{ddd}
 
\subsubsection{}
 
If $-\lambda$ is admissible in the sense of \ref{neuj300} (with $\vp$ replaced by $\tilde \vp$),
then the push-down
$$\pi^U_*:B_{\{1\}}(1_\lambda,\tilde{\vp})\to B_U(1_\lambda,\tilde{\vp})$$ is regular and admits a right-inverse. Hence it induces an isomorphism 
$$B_U(1_\lambda,\tilde{\vp})\cong B_{\{1\}}(1_\lambda,\tilde{\vp})/\ker\pi^U_*\ .$$
If $\lambda\not\in I_\aaaa$, then $J_\lambda$ is regular. We claim that it maps $\ker\pi^U_*$
to the kernel of $\pi^*_U$ at $-\lambda$. Let $f\in\ker\pi^U_*$ and $\phi\in D_U(1_{\lambda},\vp)$.
By (\ref{uhu}) there is an $\psi\in\ D_U(1_{-\lambda},\vp)$ such that 
$J_\lambda\circ ext^U(\phi)=ext^U(\psi)$. We find
$$
\langle \pi^U_*\circ J_\lambda (f),\phi\rangle=
\langle f,J_\lambda\circ ext^U(\phi)=
\langle f,ext^U(\psi)\rangle=
\langle \pi^U_*(f),\psi\rangle=0\ .
$$
The claim follows.
Therefore, if $\pm\lambda$ is admissible and non-integral, then  the operator $J_\lambda$
descends to a map 
\begin{equation}\label{neuj801}
\tilde S^U_\lambda:
B_U(1_\lambda,\tilde{\vp})\rightarrow B_U(1_{-\lambda},\tilde{\vp})\ .
\end{equation}
On the other hand, we have a meromorphic family of maps given by the adjoint of the scattering matrix 
$${}^tS^U_\lambda:B_U(1_\lambda,\tilde{\vp})\rightarrow B_U(1_{-\lambda},\tilde{\vp})\ .$$

\begin{lem} \label{gen}
If $\pm\lambda$ is non-integral and admissible, then we have
$\tilde S^U_\lambda={}^tS^U_\lambda$.
\end{lem}
\proof
Recall from Lemma \ref{wauwau} that  $ext^U\circ res^U$
is the identity on $Ext_{U}(1_{-\lambda},\vp)$. For non-integral $\lambda$ we have $Ext_{U}(1_{-\lambda},\vp)={}^U C^{-\infty}(\partial X,V(1_{-\lambda},\vp))$ by Proposition \ref{maertins},~2.

The assumptions on $\lambda$ imply that 
$\tilde S^U_\lambda$ is defined, that 
$\pi^U_*:B_{\{1\}}(1_\lambda,\tilde{\vp})\rightarrow B_U(1_\lambda,\tilde{\vp})$, $ext^U: D_U(1_{\lambda},\vp)\rightarrow
D_{\{1\}}(1_{\lambda},\vp)$, and $J_\lambda$ are regular.
Let $f\in  B_U(1_\lambda,\tilde{\vp})$ be given by $\pi^U_*(F)$, $F\in B_{\{1\}}(1_{\lambda},\vp)$.
Furthermore let $\phi\in D_U(1_{\lambda},\vp)$.
Then we compute
\begin{eqnarray*}
\langle \tilde S^U_\lambda(f),\phi\rangle&=&
\langle \pi^U_*\circ J_\lambda(F),\phi\rangle\\
&=&\langle F,J_\lambda\circ ext^U(\phi)\rangle\\
&=&\langle F,ext^U\circ res^U\circ J_\lambda\circ ext^U(\phi)\rangle\\
&=&\langle f,S^U_\lambda(\phi)\rangle\\
&=&\langle {}^tS^U_\lambda(f),\phi\rangle\ .
\end{eqnarray*}
\hB

\subsubsection{}

Note that the normalization of $\vp$ in the present subsection differs
from that in \ref{weih213}, since we want $\tilde \vp$ to be normalized as required there.
As a consquence 
the space $R_{U}(1_\lambda,\vp)$ contains summands with negative index.
If we put $k_\vp:=l_\vp/\alpha$, then we
have $$R_{U}(1_\lambda,\vp)=\prod^\infty_{n=-k_\vp} R_{U}(1_\lambda,\vp)^n\ ,\quad 
B_U(1_\lambda,\vp)=\bigcap_{k\ge -k_\vp}B_{U,k}(1_\lambda,\vp)\ .$$

\subsubsection{}

\begin{lem}\label{fundd}
Assume that $U$ defines a cusp of smaller rank.
If $l_\vp<2\rho^U$, then there is a natural non-degenerate pairing
between
$B_{U,-k_\vp}(1_\lambda,\vp)$ and $B_{U,0}(1_{-\lambda},\tilde\vp)$
given by integration over $B_U$. We obtain an inclusion
$$B_{U}(1_\lambda,\vp)\subset B_{U,-k_\vp}(1_\lambda,\vp)\hookrightarrow  D_U(1_\lambda,\vp)\ .$$
\end{lem}
\proof
\subsubsection{}
Let $f\in B_{U,-k_\vp}(1_\lambda,\vp)$ and $\phi\in B_{U,0}(1_{-\lambda},\tilde\vp)$.
Let $W\subset N\setminus N_V$ be any compact subset.
Then there exists a constant $C\in \R$ such that for all $\xi\in W$ and 
$a\in A_+$ we have 
\begin{eqnarray}
|\langle f(\xi^aw),\phi(\xi^aw) \rangle | 
&=& |\langle\vp(a)\vp(a)^{-1} f(\xi^aw),\phi(\xi^aw) \rangle|\nonumber\\
&=&|\langle \vp(a)^{-1} f(\xi^aw),\tilde\vp(a)^{-1}\phi(\xi^aw)\rangle|\nonumber\\
&\le& C \|f\|_{-k_\vp,0} \|\phi\|_{0,0} a^{l_\vp-4\rho^U}\label{stegh}
\end{eqnarray}
(see \ref{weih136} for the definition of the norms).

\subsubsection{}
Let $F\subset N_V$ be any compact fundamental domain for the lattice $V\subset N_V$.
We can now write
\begin{eqnarray*}
&&\int_{U\backslash N}|\langle f(xw),\phi(xw) \rangle | dx
\\&=&\frac{1}{[U:U^0]}\int_{N_V\backslash N}  
\int_{F} |\langle f(yvw),\phi(yvw) \rangle |dv dy\\
&=&\frac{1}{[U:U^0]}\int_A \int_{S(N_V\backslash N)} 
\int_{F} |\langle f(\xi^a vw),\phi(\xi^avw) \rangle |dv d\xi a^{2\rho^U}  da\\
&=&\frac{1}{[U:U^0]}\int_A \int_{S(N_V\backslash N)}  
\int_{F} |\langle f((\xi v^{a^{-1}})^aw),\phi((\xi v^{a^{-1}})^aw) \rangle |dv d\xi a^{2\rho^U} da\\
&=& I_+ + I_-\ ,
\end{eqnarray*}
where $I_\pm$ are the integrals over $A_\pm$, respectively.
Inserting (\ref{stegh}) we obtain 
$$
\int_{S(N_V\backslash N)}  
\int_{F} |\langle f((\xi v^{a^{-1}})^aw),\phi((\xi v^{a^{-1}})^aw) \rangle |dv d\xi \le  C_1 \|f\|_{-k_\vp,0} \|\phi\|_{0,0}a^{l_\vp-4\rho^U}$$
and therefore
$$I_+ \le C_2 \|f\|_{-k_\vp,0} \|\phi\|_{0,0} \int_{ A_+}a^{l_\vp-2\rho^U}da \ .$$
The integral on the right-hand side converges since $l_\vp-2\rho^U<0$
by assumption. Since
$I_-$ can clearly be estimated by $C_3 \|f\|_{-k_\vp,0} \|\phi\|_{0,0}$,
we have shown the lemma. \hB

\subsubsection{}\label{neuj2005}


Using the inclusion (Lemma \ref{fundd})
$$B_{U}(1_\lambda,\vp)\hookrightarrow  D_U(1_\lambda,\vp)$$
we can consider the restriction of the scattering matrix
$$(S^U_\lambda)_{|B_{U}(1_\lambda,\vp)}:B_{U}(1_\lambda,\vp)\to  D_{U}(1_{-\lambda},\vp)\ .$$
If $f\in B_{U}(1_\lambda,\vp)$, then
the restriction of the  distribution $ext^U(f)$ to $\Omega_P$ is smooth. 
Since $J_\lambda$ is pseudo-local it
it follows that
$J_\lambda\circ ext(f)$ is smooth on $\Omega_P$. Hence the restriction of the distribution
$S^U_\lambda(f)$   to the Schwartz space
$\cS_U(1_\lambda,\tilde\vp)$ is given by integration against a smooth section.
The following lemma asserts that $S^U_\lambda(f)={}^tS_\lambda^U(f)\in B_{U}(1_{\lambda},\vp)$ considered as an element of $D_{U}(1_{-\lambda},\vp)$ via Lemma \ref{fundd}.

\subsubsection{}

\begin{lem}\label{inside}
Assume that $U$ defines a cusp of smaller rank and that $l_\vp <2\rho^U$.
Then we have an equality of meromorphic families 
\begin{equation}\label{neuj900}
(S^U_\lambda)_{|B_{U}(1_\lambda,\vp)}={}^tS^U_\lambda\ . 
\end{equation}
\end{lem}
\proof
By meromorphic continuation it suffices to show
the equality  (\ref{neuj900}) for non-integral $\lambda$ in the
open subset
$\{|\Ree(\lambda)|<\rho^U-l_\vp/2\}\subset \aca$ such that $-\lambda$ is admissible in the sense of \ref{neuj300}.
Then the push-down converges at $\pm\lambda$.
Let $f\in B_{\{1\}}(1_\lambda,\tilde\vp)$.
We consider an element $\phi\in \cS_{\{1\}}(1_\lambda, \vp)$.

We view $\tilde S^U_\lambda\circ \pi^U_*(f)$ (see (\ref{neuj801})) as an element
of $D_U(1_{-\lambda},\tilde\vp)$ and compute
\begin{eqnarray*}
\langle \tilde S^U_\lambda\circ \pi^U_*(f),\pi^U_*(\phi)\rangle_{B_U} &=&
\langle \pi^U_*\circ J_\lambda (f),\pi^U_*(\phi)\rangle_{B_U}\\
&=&\langle \pi^U_*\circ J_\lambda (f),\phi\rangle_{\partial X}\\
&=&\sum_{u\in U}  \langle \pi^{1_{-\lambda},\tilde\vp}(u) J_\lambda(f),\phi\rangle_{\partial X}\\
&=&\sum_{u\in U}  \langle f,\pi^{1_{-\lambda},\vp}(u) J_\lambda(\phi)\rangle_{\partial X}\\
&=&\langle f,\pi^U_* \circ J_\lambda (\phi)\rangle_{\partial X}\\
&=&\langle \pi_*^U(f),\pi^U_* \circ J_\lambda (\phi)\rangle_{B_U}\\
&=&\langle ext^U\circ \pi_*^U(f),J_\lambda (\phi)\rangle_{\partial X}\\
&=&\langle J_\lambda\circ ext^U\circ \pi^U_*(f),\phi\rangle_{\partial X}\\
&=&\langle ext^U\circ res^U\circ J_\lambda\circ ext^U\circ \pi^U_*(f),\phi\rangle_{\partial X}\\
&=&\langle S^U_\lambda\circ \pi^U_*(f),\pi^U_*(\phi)\rangle_{B_U}\ .
\end{eqnarray*}
Since the $\pi^U_*(\phi)$, $\phi\in  \cS_{\{1\}}(1_\lambda,\vp)$ (resp. $\pi_*^U(f)\in B_{U}(1_\lambda,\tilde\vp)$, $f\in B_{\{1\}}(1_\lambda,\tilde\vp)$),
exhaust $\cS_U(1_\lambda,\vp)$ (resp. $B_{U}(1_\lambda,\tilde\vp)$) we conclude that
$$\tilde S^U_\lambda(h) =S^U_\lambda(h)  $$ as smooth sections for all $h\in B_{U}(1_\lambda,\tilde\vp)$
 and hence $${}^tS^U_\lambda = (S^U_\lambda)_{|B_{U}(1_\lambda,\vp)} $$
by Lemma \ref{gen}. \hB

\subsubsection{}
Next we show that the scattering matrix is off-diagonally
smoothing. Let $\chi,\tilde\chi$ be smooth cut-off functions
on $B_U$ with compact support such that
$\tilde\chi \chi = \chi$.

\begin{lem}\label{offside}
Assume that $U$ defines a cusp of smaller rank and that $l_\vp <2\rho^U$.
Then we have a meromorphic family
of continuous maps
$$(1-\tilde\chi)\circ  S^U_\lambda \circ \chi:D_U(1_\lambda,\vp)\rightarrow B_U(1_{-\lambda},\vp)\ .$$
\end{lem}
\proof
We are going employ the fact that $J_\lambda$ is off-diagonally
smoothing.
Let $\chi^U\in C^\infty(\Omega_P)$ be a cut-off function such that
$\sum_{u\in U} u^* \chi^U\equiv 1$, and such that
the restriction
of the projection $\Omega_P\rightarrow B_U$ to $\supp(\chi^U)$
is proper. Then there exists a cut-off function
$\bar\chi$ on $\partial X$
such that $\bar\chi \chi\chi^U=0$ and
$\bar\chi(1-\tilde\chi)=1-\tilde\chi$.
Finally we choose a compactly supported cut-off function $\hat\chi$
on $B_U$ such that $\chi^U \bar \chi \hat \chi=0$ and $\hat \chi \chi=\chi$.

Let $\phi\in \cS_U(1_{\lambda},\tilde\vp)$ be any test function
and  $f\in D_U(1_\lambda,\vp)$. We approximate
$\chi f$ by a sequence of smooth functions
$f_\alpha$ with compact support.
We can in addition assume that $\supp f_\alpha\subset \{\hat \chi=1\}$
for all $\alpha$.
Then we compute
\begin{eqnarray*}
\langle (1-\tilde\chi)S^U_\lambda(\chi f),\phi\rangle_{B_U}&=&
\lim_{\alpha}  \langle (1-\tilde\chi)S^U_\lambda f_\alpha,\phi\rangle_{B_U}\\
&\stackrel{\mathrm{Lemma}\,\ref{inside}}{=}&\lim_{\alpha}  \langle (1-\tilde\chi){}^tS^U_\lambda f_\alpha,\phi\rangle_{B_U}\\
&=&\lim_{\alpha}  \langle f_\alpha, S^U_\lambda( (1-\tilde\chi) \phi)\rangle_{B_U}\\
&=&\lim_{\alpha}  \langle f_\alpha, \{res^U\} \circ J_\lambda \circ ext^U((1-\tilde\chi) \phi)\rangle_{B_U}\\
&=&\lim_{\alpha} \langle \chi^U f_\alpha ,   J_\lambda \circ ext^U((1-\tilde\chi) \phi)\rangle_{\partial X}\\&=&\lim_{\alpha} \langle \chi^U\hat \chi f_\alpha ,   J_\lambda \circ \bar \chi ext^U((1-\tilde\chi) \phi)\rangle_{\partial X}\\
&=&\lim_{\alpha}\langle \pi^U_* (\bar\chi J_\lambda (\chi^U \hat \chi f_\alpha)),(1-\tilde\chi) \phi      \rangle_{B_U}\\
&\stackrel{(*)}{=}&\langle \pi^U_* \left(\bar\chi J_\lambda ( \chi^U \chi f)-\tilde\chi\pi^U_*  (\bar\chi J_\lambda (\chi \chi^U f))\right) , \phi      \rangle_{B_U}\ .
\end{eqnarray*}
In order to see the equality marked by $(*)$ note that
that $\bar\chi J_\lambda \hat\chi\chi^U$
is a continuous map from $D_U(1_\lambda,\vp)$ to $B_{\{1\}}(1_{-\lambda},\vp)$
(since $J_\lambda$ is off-diagonal smoothing) and that
$\lim_{\alpha}\chi^U \hat \chi f_\alpha=\chi^U\chi f$.  The assertion of the Lemma  follows from
the identity proved above:
$$(1-\tilde\chi)\circ S^U_\lambda\circ \chi (.) = \pi^U_* \left(\bar\chi J_\lambda (\chi \chi^U .)-\tilde\chi \pi^U_*  (\bar\chi J_\lambda (\chi \chi^U .))\right)\ .$$
\hB

\section{The general case} 

\subsection{The space $B_\Gamma(\sigma_\lambda,\vp)$}\label{samel1}

\subsubsection{}

Let $\Gamma\subset G$ be a torsion-free geometrically finite  discrete subgroup (see Definition \ref{t67}) such that all
its cusps are regular (Definition \ref{t799}). Furthermore let $(\vp,V_\vp)$ be an admissible
twist (Definition \ref{weih134}). We are going to employ the notation introduced in Subsection \ref{geomf}.

\subsubsection{}

On $\bar Y$ we choose a partition of unity $\{\chi_p\}_{p\in \cP_\Gamma\cup \{0\}}$ such that
$\chi_p\in C^\infty(\bar Y)$ and $\supp(\chi_p)\subset \bar{Y}_p$ for all $p$.
Restriction of these functions to $B_\Gamma$ gives a partition
of unity $\{\chi_p\}_{p\in\cP_\Gamma^<\cup\{0\}}$ on  $B_\Gamma$ such that
$\supp(\chi_p)\subset B_p$. Here we denote the restriction
of $\chi_p$ to the boundary by the same symbol $\chi_p$.
In a similar manner for $P\in p\in \cP^<_\Gamma$
we let $e_P:B_p\rightarrow
B_{U_P}$ denote the map defined as restriction of the map
$e_P:\bar Y_p\rightarrow \bar Y_{U_P}$.
Let $\chi_P\in C^\infty(\bar{Y}_{U_P})$ be the cut-off
function  which is supported on
the range of $e_P$ and satisfies $\chi_p=e_P^*\chi_P$.

\subsubsection{}\label{neuj3005}

Using cut-off with $\chi_p$ and the map $e_P:B_p\rightarrow
B_{U_P}$ for  $P\in p\in\cP^<_\Gamma$ we define maps
\begin{eqnarray*}
T_P&:&C^\infty(B_\Gamma,V_{B_\Gamma}(\sigma_\lambda,\vp))\rightarrow C^\infty(B_{U_P},V_{B_{U_P}}(\sigma_\lambda,\vp))\\
T^P&:&C^\infty(B_{U_P},V_{B_{U_P}}(\sigma_\lambda,\vp))\rightarrow C^\infty(B_\Gamma,V_{B_\Gamma}(\sigma_\lambda,\vp))\ .
\end{eqnarray*}

\subsubsection{}

We define the Schwartz space for $\Gamma$ as the space of smooth sections
which belong to the  Schwartz space \ref{weih245} near all cusps.

\begin{ddd}
We define $\cS_{\Gamma,k}(\sigma_\lambda,\vp)$ to be the subspace
of all $f\in C^\infty(B_\Gamma,V_{B_\Gamma}(\sigma_\lambda,\vp))$
such that $T_P(f)\in \cS_{U_P,k}(\sigma_\lambda,\vp)$ for all $P\in p\in \cP_\Gamma^<$.
The inclusion $$\cS_{\Gamma,k}(\sigma_\lambda,\vp)\hookrightarrow C^\infty(B_\Gamma,V_{B_\Gamma}(\sigma_\lambda,\vp))$$
and the maps $$T_P:\cS_{\Gamma,k}(\sigma_\lambda,\vp)\rightarrow \cS_{U_P,k}(\sigma_\lambda,\vp)$$
induce on $\cS_{\Gamma,k}(\sigma_\lambda,\vp)$ the structure of a Fr\'echet space. Furthermore we define the Fr\'echet and Montel space
$$\cS_{\Gamma}(\sigma_\lambda,\vp):=\bigcap_{k\in\nat_0}\cS_{\Gamma,k}(\sigma_\lambda,\vp)\ .$$
\end{ddd}

\subsubsection{}

\begin{lem}\label{loctrivg}
The families  $\{\cS_{\Gamma,k}(\sigma_\lambda,\vp)\}_{\lambda\in\aca}$
and $\{\cS_{\Gamma}(\sigma_\lambda,\vp)\}_{\lambda\in\aca}$
form locally trivial holomorphic bundles.
\end{lem}
\proof
Let $\Phi_{P,\lambda_0}$ denote the
trivialization constructed in Lemma \ref{ytra} which is
given by multiplication by $s_P^{(\lambda_0-\lambda)/\alpha}$
(here we add the index $P$ to the notation for the section $s$ constructed  in Lemma \ref{ingf} in order to indicate its dependence on the parabolic subgroup).
Let $s_p:=T^P(s_P)$, where $P\in\tilde \cP^<$ represents $p$.
Let $s_0\in C^\infty(B_\Gamma,V_{B_\Gamma}(1_{\rho+\alpha},\vp))$
be any positive  section (compare \ref{neuj911}).
Then we form $s_\Gamma:=\sum_{p\in \cP_\Gamma} s_p + \chi_0 s_0$.
Multiplication by $s_\Gamma^{(\lambda_0-\lambda)/\alpha}$  defines   isomorphisms
$$\Phi_{\Gamma,\lambda_0}: \cS_{\Gamma,k}(\sigma_\lambda,\vp)\rightarrow
\cS_{\Gamma,k}(\sigma_{\lambda_0},\vp)$$ and
$$\Phi_{\Gamma,\lambda_0}: \cS_{\Gamma}(\sigma_\lambda,\vp)\rightarrow
\cS_{\Gamma}(\sigma_{\lambda_0},\vp)\ .$$
We employ the maps $\Phi_{\Gamma,\lambda_0}$ in order to obtain the
required local trivializations and to implement the holomorphic structures.
\hB

\subsubsection{}

Similarly to the  case of the Schwartz space we define
the space $B_{\Gamma}(\sigma_\lambda,\vp)$
as the space of smooth sections  with the  same asymptotic expansions
near the cusps as the spaces as the space defined in \ref{neuj104}.
Recall the notation $\tilde \cP^{max}$ from \ref{holle}.

\begin{ddd}
We define $$B_{\Gamma,k}(\sigma_\lambda,\vp):=B_{\Gamma,k}(\sigma_\lambda,\vp)_1
\oplus \bigoplus_{P\in \tilde\cP^{max}} \cR_{U_P}(\sigma_\lambda,\vp)\ ,$$
where $$B_{\Gamma,k}(\sigma_\lambda,\vp)_1\subset  C^\infty(B_\Gamma,V_{B_\Gamma}(\sigma_\lambda,\vp))$$ is
 the subspace of all $f$ such that
$T_P(f)\in B_{U_P,k}(\sigma_\lambda,\vp)$ for all $P\in p\in \cP_\Gamma^<$.  
The map $$B_{\Gamma,k}(\sigma_\lambda,\vp)\hookrightarrow C^\infty(B_\Gamma,V_{B_\Gamma}(\sigma_\lambda,\vp))\ ,$$
the maps $$T_P:B_{\Gamma,k}(\sigma_\lambda,\vp)\rightarrow B_{U_P,k}(\sigma_\lambda,\vp)\ ,$$ $P\in \tilde \cP^<$,  
and the natural projections $$AS_P:B_{\Gamma,k}(\sigma_\lambda,\vp)\rightarrow \cR_{U_P}(\sigma_\lambda,\vp)\ ,$$ $P\in \tilde\cP^{max}$,
equip $B_{\Gamma,k}(\sigma_\lambda,\vp)$ with the structure of a Fr\'echet
space.
We further define the Fr\'echet and Montel space
 $$B_{\Gamma}(\sigma_\lambda,\vp):=
\bigcap_{k\in\nat_0} B_{\Gamma,k}(\sigma_\lambda,\vp)\ .$$
\end{ddd}

\subsubsection{}

We define $$\cR_{\Gamma,k}(\sigma_\lambda,\vp):=\bigoplus_{P\in \tilde\cP} \cR_{U_P,k}(\sigma_\lambda,\vp) \ ,\quad \cR_{\Gamma}(\sigma_\lambda,\vp) :=\bigoplus_{P\in \tilde\cP} \cR_{U_P}(\sigma_\lambda,\vp)\ .$$
We have asymptotic term maps $$AS_P: B_{\Gamma,k}(\sigma_\lambda,\vp)\rightarrow \cR_{U_P,k}(\sigma_\lambda,\vp)\ ,$$
$AS_P:=AS\circ T_P$, where $AS$ was defined in  \ref{neuj915}
for pure cusps. The asymptotic term maps admit right-inverses
   $$L_P:\cR_{U_P,k}(\sigma_\lambda,\vp)\rightarrow B_{\Gamma,k}(\sigma_\lambda,\vp)$$
given by the natural inclusion for $P\in\tilde\cP^{max}$,  and by $L_P:=T^P\circ L$
if $P\in \tilde\cP^<$, where $L$ was defined in  \ref{weih253}
for pure cusps.

 Let $$AS:B_{\Gamma,k}(\sigma_\lambda,\vp)\rightarrow
\cR_{\Gamma,k}(\sigma_\lambda,\vp)$$ be induced by the maps $AS_P$ and $$L: \cR_{\Gamma,k}(\sigma_\lambda,\vp)
\rightarrow B_{\Gamma,k}(\sigma_\lambda,\vp)$$ be given by the sum of the maps
$L_P$ for the various $P\in\tilde\cP$.

\subsubsection{}

\begin{lem}\label{spst}
The family $\{B_{\Gamma,k}(\sigma_\lambda,\vp)\}_{\lambda\in\aca}$
forms a trivial holomorphic bundle of Fr\'echet
spaces. We have a split (by $L$) exact sequence
\begin{equation}\label{trtra}0\rightarrow
\cS_{\Gamma,k}(\sigma_\lambda,\vp)\rightarrow B_{\Gamma,k}(\sigma_\lambda,\vp)
\stackrel{AS}{\rightarrow} \cR_{\Gamma,k}(\sigma_\lambda,\vp)\rightarrow 0\ .\end{equation}
Furthermore, the  spaces $\{B_{\Gamma}(\sigma_\lambda,\vp)\}_{\lambda\in\aca}$ form a limit of locally trivial holomorpic bundles in the sense
of Subsection \ref{llim} and fit into the exact sequence (which does not admit any continuous split)
$$0\rightarrow \cS_{\Gamma}(\sigma_\lambda,\vp) \rightarrow B_{\Gamma}(\sigma_\lambda,\vp) \rightarrow
\cR_{\Gamma}(\sigma_\lambda,\vp)\rightarrow 0\ .$$
\end{lem}
\proof
The  assertions follow from Lemma \ref{loctrivg} and the fact that the spaces $\cR_{\Gamma,k}(\sigma_\lambda,\vp)$ form locally trivial holomorphic vector bundles. \hB

\subsubsection{}

We now show that the spaces $B_\Gamma(\sigma_\lambda,\vp)$
are compatible with twisting.
Let $(\pi_{\sigma,\mu},V_{\pi_{\sigma,\mu}})$ be a finite-dimensional
representation of $G$ as in  \ref{weih144}.

\begin{lem}\label{compat3}
\begin{enumerate}
\item
We have holomorphic families of continuous maps (see \ref{neuj1000} for notation)
$$i_{\sigma,\mu}^\Gamma:B_{\Gamma}(\sigma_\lambda,\vp)\rightarrow
B_\Gamma(1_{\lambda+\mu},\pi_{\sigma,\mu}\otimes \vp)$$
and
$$p_{\sigma,\mu}^\Gamma:B_\Gamma(1_{\lambda-\mu},\pi_{\sigma,\mu}\otimes \vp)\rightarrow
B_{\Gamma}(\sigma_\lambda,\vp)\ .$$
\item
$B_\Gamma(1_{\lambda+\mu},\pi_{\sigma,\mu}\otimes \vp)$
is a $\cZ(\gaaa)$-module.
\item
If $\lambda\not\in I_\aaaa$, then
$i_{\sigma,\mu}^\Gamma$ maps $B_{\Gamma}(\sigma_\lambda,\vp)$
isomorphically onto $$\ker_\Gamma(Z(\lambda)):=\ker\{Z(\lambda):
B_\Gamma(1_{\lambda+\mu},\pi_{\sigma,\mu}\otimes \vp)\rightarrow B_\Gamma(1_{\lambda+\mu},\pi_{\sigma,\mu}\otimes \vp)\}$$
(see \ref{weih146} for notation),
and the restriction of
$p_{\sigma,\mu}^\Gamma$ to
$$\ker\{\tilde Z(\lambda):B_\Gamma(1_{\lambda-\mu},\pi_{\sigma,\mu}\otimes \vp)\rightarrow B_\Gamma(1_{\lambda-\mu},\pi_{\sigma,\mu}\otimes \vp)\}$$
is an isomorphism onto $B_{\Gamma}(\sigma_\lambda,\vp)$
(here $\tilde Z(\lambda)$ is the adjoint of
$\pi^{1_{-\lambda+\mu}\otimes \tilde\pi_{\sigma,\mu},  \id_{\tilde\vp}}(\Pi(\lambda))$).
\item
The composition
$$j_{\sigma,\mu}^\Gamma:B_\Gamma(1_{\lambda+\mu},\pi_{\sigma,\mu}\otimes \vp)
\stackrel{1-Z(\lambda)}{\rightarrow} \ker_\Gamma(Z(\lambda))\stackrel{(i_{\sigma,\mu}^\Gamma)^{-1}}{\rightarrow}
B_{\Gamma}(\sigma_\lambda,\vp)$$
which is initially defined for $\lambda\not\in I_\aaaa$
extends to a meromorphic family of continuous maps.
Similarly, the composition
$$q_{\sigma,\mu}^\Gamma:B_{\Gamma}(\sigma_\lambda,\vp)
\stackrel{(p_{\sigma,\mu}^\Gamma)^{-1}}{\rightarrow} \ker_\Gamma(\tilde Z(\lambda)) \rightarrow B_\Gamma(1_{\lambda-\mu},\pi_{\sigma,\mu}\otimes \vp)$$
extends to a tame meromorphic family of continuous maps.
\end{enumerate}
\end{lem}
\proof
\subsubsection{}
In order to see 1. we employ Lemma \ref{compat2} 2., and
the identities $$i_{\sigma,\mu}^{U_P}\circ T_P=T_P\circ i_{\sigma,\mu}^\Gamma\ ,\quad 
p_{\sigma,\mu}^{U_P}\circ T_P=T_P\circ p_{\sigma,\mu}^\Gamma\ .$$

\subsubsection{}
In order to see 2. we use Lemma \ref{compat2}, 3., and the fact
that $[\pi^{1_{\lambda+\mu}\otimes \pi_{\sigma,\mu},\id_\vp}(A),\chi_p]$
is a differential operator with compactly supported coefficients on $B_\Gamma$
for any $A\in\cZ(\gaaa)$.

\subsubsection{}

We show the first assertion of 3. and leave the second to the reader,
since the argument is similar.
Since $i^\Gamma_{\sigma,\mu}$ is $\cZ(\gaaa)$-equivariant
it maps to $\ker_\Gamma(Z(\lambda))$ by construction of $Z(\lambda)$.
Let now $f\in\ker_\Gamma(Z(\lambda))$.
We  claim that for any
$h\in C^\infty(B_\Gamma)$ we have $Z(\lambda)(hf)=0$.
Let $k$ be the order of the differential operator $Z(\lambda)$ which is a projection.
Then we have
\begin{eqnarray*}
Z(\lambda)(hf)&=& Z(\lambda)^{k+1}(hf)\\
&=&Z(\lambda)^{k}[Z(\lambda),h]f\\
&=&Z(\lambda)^{k-1}[Z(\lambda),[Z(\lambda),h]]f\\
&\dots&\\
&=&\underbrace{[Z(\lambda),\dots,[Z(\lambda),h]\dots]}_{k+1}f\\
&=&0\ .
\end{eqnarray*}
This shows that claim.

We conclude that $T_P(f)\in \ker_U(Z(\lambda))$ for
all $P\in\tilde P$.
By Lemma \ref{compat2} we find
$g_P\in B_U(\sigma_\lambda,\vp)$, $P\in\tilde\cP$,
such that $T_P(f)=i_{\sigma,\mu}^{U_P}(g_P)$.
Let $$f_0:=f-\sum_{P\in\tilde \cP} i_{\sigma,\mu}^\Gamma\circ T^P (g_P)\in
C_c^\infty(B_\Gamma,V_{B_\Gamma}(1_{\lambda+\mu},\pi_{\sigma,\mu}\otimes\vp))\ .$$
As in the proof of Lemma \ref{compat1},4. we can find
$g_0\in C_c^\infty(B_\Gamma,V_{B_\Gamma}(\sigma_\lambda,\vp))$ such that
$f_0=i^\Gamma_{\sigma,\mu}(g_0)$.
Thus $f= i^\Gamma_{\sigma,\mu}\left(\sum_{P\in\tilde \cP\cup\{0\}} g_P\right)$.

\subsubsection{}\label{neuj3006}

We show the first assertion of 4. and leave the second to the reader.
As in the proof of Lemma \ref{compat1}, 5. we construct
a  holomorphic family
$$J:C_c^\infty(B_\Gamma,V_{B_\Gamma}(1_{\lambda+\mu},\pi_{\sigma,\mu}\otimes\vp))
\rightarrow C_c^\infty(B_\Gamma,V_{B_\Gamma}(\sigma_\lambda,\vp))$$
as the composition $J:=j\circ (1-Z(\lambda))$,
where $$j:C_c^\infty(B_\Gamma,V_{B_\Gamma}(1_{\lambda+\mu},\pi_{\sigma,\mu}\otimes\vp))
\rightarrow C_c^\infty(B_\Gamma,V_{B_\Gamma}(\sigma_\lambda,\vp))$$
comes from a bundle homomorphism.
Then we define
$$\tilde j^\Gamma_{\sigma,\mu}:= J \chi_0 (1-Z(\lambda))+
\sum_{P\in\tilde \cP} \tilde T^P\circ j^{U_P}_{\sigma,\mu} \circ T_P\circ (1-Z(\lambda))\ .$$
Here $\tilde T^P$ is defined as $T^P$ but
using a cut-off function $\tilde\chi_p\in C^\infty(B_\Gamma)$, $P\in p$,
satisfying $\supp(\tilde\chi_p)\in B_p$ and $\tilde\chi_p\chi_p=\chi_p$
if $\tilde\cP^<$. If $P\in\tilde\cP^{max}$, then we set
$\tilde T^P:=T^P$. One checks that $\tilde j^\Gamma_{\sigma,\mu}\circ
i^\Gamma_{\sigma,\mu}=\id$.
Since $\tilde j^\Gamma_{\sigma,\mu}$ vanishes on
$\ker(1-Z(\lambda))$ we conclude that it coincides
with $j^\Gamma_{\sigma,\mu}$ for $\lambda\not\in I_\aaaa$.
We thus have shown that $j^\Gamma_{\sigma,\mu}$
extends to a  meromorphic family of continuous maps.
\hB

%

\subsection{Push-down}

\subsubsection{}

For the first part of the present subsection
we choose an Iwasawa decomposition
$G=KAN$ and a parabolic subgroup $P=MAN$, $M\subset K$.
Then we write $X=G/K$ and $\partial X=G/P=K/M$.

\subsubsection{}

For $f\in C^\infty(\partial X,V(\sigma_\lambda,\vp))$
we consider the push-down
$$\pi^\Gamma_*(f)\in C^\infty(B_\Gamma,B_{B_\Gamma}(\sigma_\lambda,\vp))\oplus \bigoplus_{P\in\tilde\cP^{max}} \cR_{U_P}(\sigma_\lambda,\vp)\ .$$
It is given by the sum
$$(\pi^\Gamma_*(f))_1:=\sum_{g\in\Gamma}\pi^{\sigma_\lambda,\vp}(g) (f_{|\Omega_\Gamma})\in{}^\Gamma C^\infty(\Omega_\Gamma,V(\sigma_\lambda,\vp))$$
in the first component. Its second component
$$(\pi^\Gamma_*(f))_2\in \bigoplus_{P\in\tilde\cP^{max}} \cR_{U_P}(\sigma_\lambda,\vp)$$
will be constructed below.
The goal of the present section is to show
that the push-down
converges for $\Ree(\lambda)<-\delta_\Gamma-\delta_\vp$ and defines
a holomorphic family of maps
$$\pi^\Gamma_*:C^\infty(\partial X,V(\sigma_\lambda,\vp))\rightarrow
B_{\Gamma}(\sigma_\lambda,\vp)\ .$$
\subsubsection{}

In the following lemma we consider points of $\partial X$ as subsets of $K$ using
the identification $\partial X\cong K/M$.
\begin{lem}\label{ggll}
If $W\subset\Omega_\Gamma$ is compact, then $\Gamma\cap WMA_+K$ is finite.
\end{lem}
\proof
The set $ WMA_+K$ is a precompact subset of $X\cup \Omega_\Gamma$. Since
$\Gamma$ acts properly discontinuously on $X\cup \Omega_\Gamma$ the intersection of
the orbit $\Gamma K$ with $WMA_+K$ is finite. \hB

\subsubsection{}

\begin{lem}\label{barpi}
For $\Ree(\lambda)<-\delta_\Gamma-\delta_\vp$
the sum $$\sum_{g\in\Gamma}\pi^{\sigma_\lambda,\vp}(g)f_{|\Omega_\Gamma}$$ converges in $C^\infty(\Omega_\Gamma,V(\sigma_\lambda,\vp))$ and defines a
holomorphic family of continuous maps 
$$C^\infty(\partial X,V(\sigma_\lambda,\vp))\to 
{}^\Gamma  C^\infty(\Omega_\Gamma,V(\sigma_\lambda,\vp))\ .$$
\end{lem}
\proof
Using Lemma \ref{ggll} we employ the same argument as in the
proof of \cite{MR1749869}, Lemma 4.2. \hB

\subsubsection{}

Next we study the behaviour of the sum $\sum_{g\in\Gamma}\pi^{\sigma_\lambda,\vp}(g)f_{|\Omega_\Gamma}$
near the cusps. Let $P\in\tilde\cP$. We choose a system $\Gamma^P$ of representatives of  $U_P\backslash \Gamma$
by taking in each class $[g]\in U_P\backslash \Gamma$ an element  $h\in [g]$ which minimizes
$\dist_X(\cO,h\cO)$.

\begin{lem}\label{wwe}
There is a neighbourhood $W\subset\partial X$ of $\infty_P$ such that
$WMA_+K\cap \Gamma^P$ is finite.
\end{lem}
\proof
Let $D(\cO,U_P):=\{x\in X|\dist_X(x,\cO)\le \dist_X(gx,\cO)\:\forall g\in U_P\}$ be the
Dirichlet domain of $U_P$. If $h\in \Gamma^P$, then $h\cO\in D(\cO,U_P)$.
Let $\pi:D(\cO,U_P)\rightarrow Y_{U_P}$ denote the projection and
consider $E:=\pi^{-1}(e_P(Y_p))$. Since
$e_P:Y_p\rightarrow Y_{U_P}$ is an isometry
we have
$\sharp(E\cap \Gamma^P\cO)\le 1$. Now the closure in $\bar X$ of
$D(\cO,U_P)\setminus E$ does not contain $\infty_P$ and is therefore
disjoint from a neighbourhood $W\subset \partial X$ of $\infty_P$.
Since $WMA_+K\cap (D(\cO,U_P)\setminus E)$ has a compact
closure inside $X$ we conclude that $WMA_+K\cap (D(\cO,U_P)\setminus E)K\cap
\Gamma^P$ is finite. Since $\Gamma^P\subset D(\cO,U_P)K$ we have shown the lemma.\hB

\subsubsection{}

Let $P\in \tilde \cP$ and $W_P$ be a neighbourhood of $\infty_P$ as constructed in Lemma \ref{wwe}. 
Using Lemma \ref{wwe} and the arguments of the proof of 
\cite{MR1749869}, Lemma 4.2, we show
\begin{lem}\label{bbarpi}
For $\Ree(\lambda) < -\delta_\Gamma-\delta_\vp$ and $f\in C^\infty(\partial X,V(\sigma_\lambda,\vp))$ the sum 
$$\bar\pi^P_*(f):=\sum_{g\in\Gamma^P}(\pi^{\sigma_\lambda,\vp}(g)f)_{|W_P\cup\Omega_\Gamma}$$
converges in $C^\infty(W_P\cup\Omega_\Gamma,V(\sigma_\lambda,\vp))$.
For varying $\lambda$ the maps $\bar\pi^P_*$ form a  holomorphic family of continuous maps.
\end{lem} \hB

\subsubsection{}

We choose a cut-off function $\kappa_P\in C^\infty_c(W_P)$
such that $\kappa_P$ is equal to one near $\infty_P$.  
Then $\kappa_P \bar\pi^P_*(f))\in C^\infty(\partial X,V(\sigma_\lambda,\vp))$,
and we can write
$$(\pi^\Gamma_*(f))_1=\left((\pi^{U_P}_*(\kappa_P \bar\pi^P_*(f)))_1 + \sum_{u\in U_P} \pi^{\sigma_\lambda,\vp}(u)(1-\kappa_P)
\bar\pi^P_*(f)\right)_{|\Omega_\Gamma}\ .$$
Here $$(\pi^{U_P}_*(\kappa_P \bar\pi^P_*(f)))_1=\pi^{U_P}_*(\kappa_P \bar\pi^P_*(f))$$
if the cusp associated to $U_P\subset P$ has smaller rank.
The second sum is locally finite and defines a smooth section on
$\Omega_\Gamma$. Moreover we have ($\bar \pi^P_*$ was defined in Lemma \ref{bbarpi})
$$\chi_P \sum_{u\in U_P} \pi^{\sigma_\lambda,\vp}(u)(1-\kappa_P)
\bar\pi^P_*(f)_{|\Omega_\Gamma}\in C^\infty_c(B_{U_P},V_{B_{U_P}}(\sigma_\lambda,\vp))\ .$$
 
\subsubsection{}
We define
$$(\pi^\Gamma_*(f))_2:=\bigoplus_{P\in\tilde\cP^{max}} [\pi^{U_P}_*]\circ AS (\kappa_P \bar\pi^P_*(f))\ .$$
This definition is independent of the choice of $\kappa_P$.
Finally we define $\pi^\Gamma_*(f)$ to be the sum  $(\pi^\Gamma_*(f))_1\oplus (\pi^\Gamma_*(f))_2$.

\subsubsection{}

\begin{lem}
For $\Ree(\lambda) < -\delta_\Gamma-\delta_\vp$ the push-down
$\pi^\Gamma_*$ induces a holomorphic family of continuous maps
$$\pi^\Gamma_*:C^\infty(\partial X,V(\sigma_\lambda,\vp))\rightarrow B_\Gamma(\sigma_\lambda,\vp)\ .$$
\end{lem}
\proof
This follows from Lemmas \ref{barpi}, \ref{bbarpi}, the observation that
\begin{equation}\label{xxxt}T_P(\pi^\Gamma_*(f))_1\in \chi_P \pi^{U_P}_*(\kappa_P \bar\pi^P_*(f)) +
 C_c^\infty(B_{U_P},V_{B_{U_P}}(\sigma_\lambda,\vp))\subset B_{U_P}(\sigma_\lambda,\vp)\end{equation}
for $P\in \tilde\cP^<$, and the definition of  $(\pi^\Gamma_*(f))_2$.
\hB

\subsubsection{}\label{neuj1003}

Let $\chi^\Gamma\in C^\infty(\Omega_\Gamma)$ be a smooth cut-off function
such that $\sum_{g\in\Gamma} g^*\chi^\Gamma \equiv 1$.
We choose $\chi^\Gamma$ such that it coincides with $\chi^{U_P}$
in a neighbourhood of $\infty_P$, and such that
the restriction of the projection $\Omega_\Gamma\rightarrow B_\Gamma$
to $\supp(\chi^\Gamma)$ is proper.
Then multiplication by
$\chi^\Gamma$ defines a holomorphic family of right-inverses 
of the push-down (see Lemma \ref{schwspl}) $$\{Q\}:\cS_\Gamma(\sigma_\lambda,\vp)\rightarrow
C^\infty(\partial X,V(\sigma_\lambda,\vp))$$
such that $\pi^{\Gamma}_*\circ \{Q\}=\id$.

\subsubsection{}
Occasionally we will need the partial push-down
$$\pi^{\Gamma/U_P}_*:B_{U_P,k}(\sigma_\lambda,\vp)\rightarrow B_{\Gamma,k}(\sigma_\lambda,\vp)$$
which is defined as the average over $\Gamma/U_P$ if the cusp associated to
$U_P\subset P$ has smaller rank.
If the cusp associated to $U_P\subset P$ has full rank,
then we define 
$$\pi^{\Gamma/U_P}_*(f\oplus v):=\pi^\Gamma_*(\{Q\}(f-(\pi^{U_P}_*\circ L\circ [Q](v))_1)+L \circ [Q](v))\ ,$$
where
$$f\oplus v\in C^\infty(B_{U_P},V_{B_{U_P}}(\sigma_\lambda,\vp))\oplus \cR_{U_P}(\sigma_\lambda,\vp)=B_{U_P}(\sigma_\lambda,\vp)\ .$$
If $P\in \tilde\cP^<$, then we have
$$\pi^{\Gamma/U_P}_*(f):=\pi^\Gamma_*(\{Q\}(f-\pi^{U_P}_*\circ L\circ [Q]\circ AS(f))+L \circ [Q] \circ AS (f))\ .$$
This formula together with (\ref{xxxt}) can be used to verify
that $\pi^{\Gamma/U_P}_*$ has the required mapping properties.
In particular we conclude that if
$P\in\tilde\cP^<$, then the partial
push-down converges for 
$\Ree(\lambda)<-\delta_\Gamma-\delta_\vp$ and depends
holomorphically on $\lambda$. If $P\in\tilde\cP^{max}$, then the partial push-down
is defined as a meromorphic family of continuous maps
for $\Ree(\lambda)<-\delta_\Gamma-\delta_\vp$.

\subsubsection{}

We collect the following useful identities
\begin{equation}\label{useid}\pi^{\Gamma/U_P}_*\circ \pi^{U_P}_*=\pi^\Gamma_*,\quad
\pi^{\Gamma/U_P}_* \circ T_P = \chi_p,\quad \pi^{\Gamma/U_P}_* \chi_P = T^P\ ,\end{equation}
where $P\in \tilde\cP^<$ for the last two equations
and multiplication by $\chi_p$ implicitly involves
a projection onto the first component.
If $P\in \tilde\cP^{max}$, then we have
\begin{equation}\label{useid2}
\pi^{\Gamma/U_P}_* \circ AS_P =AS_P 
\ .\end{equation}

\subsubsection{}

We define $q_{\tilde\sigma,\mu}:=q^{\{1\}}_{\tilde\sigma,\mu}$ as in Lemma \ref{compat3}, 4. (applied to the trivial group).
\begin{lem}\label{compat5}
We have the following commutative
diagrams of meromorphic families of maps:
\begin{enumerate}
\item
$$\begin{array}{ccccc}
C^\infty(\partial X,V(\sigma_\lambda,\vp))&\stackrel{i_{\sigma,\mu}}{\rightarrow}&C^\infty(\partial X,V(1_{\lambda+\mu},\pi_{\sigma,\mu}\otimes\vp))&\stackrel{Z(\lambda)}{\rightarrow}&
C^\infty(\partial X,V(1_{\lambda+\mu},\pi_{\sigma,\mu}\otimes\vp))\\
\downarrow \pi^\Gamma_*&&\downarrow \pi^\Gamma_*&&\downarrow \pi^\Gamma_*\\
B_\Gamma(\sigma_\lambda,\vp)&\stackrel{i^\Gamma_{\sigma,\mu}}{\rightarrow}&
\cB_\Gamma(1_{\lambda+\mu},\pi_{\sigma,\mu}\otimes\vp)&\stackrel{Z(\lambda)}{\rightarrow}&
\cB_\Gamma(1_{\lambda+\mu},\pi_{\sigma,\mu}\otimes\vp)\end{array}\ .$$
\item
$$\begin{array}{ccc}
C^\infty(\partial X,V(\tilde\sigma_{-\lambda},\tilde\vp))
&\stackrel{q_{\tilde\sigma,\mu}}{\rightarrow}&
C^\infty(\partial X,V(1_{-\lambda-\mu},\pi_{\tilde\sigma,\mu}\otimes\tilde\vp))\\
\downarrow \pi^\Gamma_*&&\downarrow \pi^\Gamma_*\\
B_\Gamma(\tilde\sigma_{-\lambda},\tilde\vp)
&\stackrel{p^\Gamma_{\tilde \sigma,\mu}}{\leftarrow}&
B_\Gamma(1_{-\lambda-\mu},\pi_{\tilde\sigma,\mu}\otimes\tilde\vp)\\
\end{array}\ .$$
\end{enumerate}
\end{lem}
\proof
The push-down is $\cZ(\gaaa)$-equivariant.
This implies commutativity of the right square of the diagram 1.
The commutativity of the
left square is obvious for the part involving
$(\pi^\Gamma_*)_1$, and for $(\pi^\Gamma_*)_2$ we invoke
Definition \ref{cv1}.

For the second diagram consider the enlarged diagram
$$\xymatrix{C^\infty(\partial X,V(1_{-\lambda+\mu^\prime},\pi_{\tilde\sigma,\mu^\prime}\otimes\tilde\vp))
\ar[d]^{\pi^\Gamma_*}&C^\infty(\partial X,V(\tilde\sigma_{-\lambda},\tilde\vp))\ar@/_-0.5cm/[r]^{q_{\tilde\sigma,\mu}}
\ar[l]^{\qquad\quad i_{\tilde \sigma,\mu^\prime}}
\ar[d]^{\pi^\Gamma_*}&C^\infty(\partial X,V(1_{-\lambda-\mu},\pi_{\tilde\sigma,\mu}\otimes\tilde\vp))
\ar[d]^{\pi^\Gamma_*}\ar[l]^{p_{\tilde \sigma,\mu}\quad}\\
B_\Gamma(1_{-\lambda+\mu^\prime},\pi_{\tilde\sigma,\mu^\prime}\otimes\tilde\vp)&
B_\Gamma(\tilde\sigma_{-\lambda},\tilde\vp)
\ar[l]^{\qquad i^\Gamma_{\tilde \sigma,\mu^\prime}}&
B_\Gamma(1_{-\lambda-\mu},\pi_{\tilde\sigma,\mu}\otimes\tilde\vp)
\ar[l] ^{p^\Gamma_{\tilde \sigma,\mu}}}\ .$$
The left square commutes by the Definition \ref{cv1} of the push-down
in the middle. The outer rectangle commutes by the naturality of the definition of the push-down.
\hB

\subsubsection{}

Let $G^n$ be one of $\{Spin(1,n),SO(1,n)_0,SU(1,n),Sp(1,n)\}$.
Recall the definition of $i^*_\Gamma$ given in Subsection \ref{embedd}.
\begin{lem}\label{comppp}
We have the following commutative diagram:
$$\begin{array}{ccc} 
C^\infty(\partial X^{n+m},V(1^{n+1}_{\lambda },\vp)) &\stackrel{i^*}{\rightarrow}&C^\infty(\partial X^n, V(1^n_{\lambda-\zeta},\vp)) \\
 \downarrow \pi^{\Gamma,n+1}_*&&\downarrow \pi^{\Gamma,n}_* \\
 B_\Gamma(1^{n+1}_{\lambda },\vp)&\stackrel{i_\Gamma^*}{\rightarrow }& B_\Gamma(1^n_{\lambda-\zeta},\vp)
\end{array}\ .$$
 \end{lem}
\proof
We use Proposition \ref{klopp}.
\hB

\subsection{Extension, restriction, and the scattering matrix}

\subsubsection{}

We now define the distribution spaces $D_{\Gamma}(\sigma_\lambda,\vp)$ associated to $\Gamma$.
These spaces are the domain of the extension $ext^\Gamma$ and the Eisenstein series (see Definition  \ref{neuj1001}).

\begin{ddd}
We define
$D_{\Gamma,k}(\sigma_\lambda,\vp)$ to be the dual space to
$B_{\Gamma,k}(\tilde\sigma_{-\lambda},\tilde\vp)$. Furthermore let
$D_{\Gamma}(\sigma_\lambda,\vp):=\bigcup_{k\in\nat_0} D_{\Gamma,k}(\sigma_\lambda,\vp)$.
\end{ddd}
Lemma \ref{spst} has the following consequence:
\begin{kor}
The spaces $\{D_{\Gamma,k}(\sigma_\lambda,\vp)\}_{\lambda\in\aca}$, form
a  trivial holomorphic bundle of dual Fr\'echet spaces
and fit into the split exact sequence
$$0\rightarrow \cR_{\Gamma,k}(\tilde\sigma_{-\lambda},\tilde\vp)^*\stackrel{AS^*}{\rightarrow} D_{\Gamma,k}(\sigma_\lambda,\vp)\rightarrow \cS_{\Gamma,k}(\tilde\sigma_{-\lambda},\tilde\vp)^*\rightarrow 0\ .$$
 The spaces $\{D_{\Gamma}(\sigma_\lambda,\vp)\}_{\lambda\in\aca}$,
are dual Fr\'echet and Montel spaces, and they
form a direct limit of  trivial bundles. Furthermore
we have the exact sequence
$$0\rightarrow \cR_{\Gamma}(\tilde\sigma_{-\lambda},\tilde\vp)^*\stackrel{AS^*}{\rightarrow} D_{\Gamma}(\sigma_\lambda,\vp)\rightarrow \cS_{\Gamma}(\tilde\sigma_{-\lambda},\tilde\vp)^*\rightarrow 0\ .$$
\end{kor}
Let $(\pi_{\sigma,\mu},V_{\pi_{\sigma,\mu}})$ be a finite-dimensional
representation of $G$ as in Subsection \ref{ttw}.

\subsubsection{}

The following Lemma states that the distribution spaces are compatible with twisting.

\begin{lem}\label{compat6}
\begin{enumerate}
\item
$D_\Gamma(1_{\lambda+\mu},\pi_{\sigma,\mu}\otimes\vp)$
is a $\cZ(\gaaa)$-module.
\item
There is a natural inclusion
$i^\Gamma_{\sigma,\mu}:D_\Gamma(\sigma_{\lambda},\vp)
\rightarrow D_\Gamma(1_{\lambda+\mu},\pi_{\sigma,\mu}\otimes\vp)$
which identifies
$D_\Gamma(\sigma_{\lambda},\vp)$ with $$\ker_\Gamma(Z(\lambda)):=\ker(Z(\lambda):
D_\Gamma(1_{\lambda+\mu},\pi_{\sigma,\mu}\otimes\vp)   \rightarrow      D_\Gamma(1_{\lambda+\mu},\pi_{\sigma,\mu}\otimes\vp))$$
for $\lambda\not\in I_\aaaa$.
\item
Moreover, the map
$$j_{\sigma,\mu}^\Gamma:D_\Gamma(1_{\lambda+\mu},\pi_{\sigma,\mu}\otimes\vp)
\stackrel{1-Z(\lambda)}{\longrightarrow} \ker_\Gamma(Z(\lambda))\stackrel{(i_{\sigma,\mu}^\Gamma)^{-1}}{\longrightarrow}
D_\Gamma(1_{\lambda+\mu},\pi_{\sigma,\mu}\otimes\vp)$$
extends to a tame meromorphic family of continuous maps.
\end{enumerate}
\end{lem}
\proof
We employ Lemma \ref{compat3}.
The map $i^\Gamma_{\sigma,\mu}$ is defined as the adjoint
of $$p^\Gamma_{\tilde\sigma,\mu}:B_\Gamma(1_{-\lambda-\mu},\pi_{\tilde\sigma,\mu}\otimes\tilde\vp)\rightarrow B_\Gamma(\tilde\sigma_{-\lambda},\tilde\vp)\ .$$
$j_{\sigma,\mu}^\Gamma$ is just the adjoint
of $$q^\Gamma_{\tilde\sigma,\mu}:B_\Gamma(\tilde\sigma_{-\lambda},\tilde\vp)\rightarrow B_\Gamma(1_{-\lambda-\mu},\pi_{\tilde\sigma,\mu}\otimes\tilde\vp)\ .$$
\hB

\subsubsection{}

\begin{ddd}For $\Ree(\lambda)>\delta_\Gamma+\delta_\vp$
we define the extension map 
$$ext^\Gamma:D_{\Gamma}(\sigma_\lambda,\vp)\rightarrow C^{-\infty}(\partial X,V(\sigma_\lambda,\vp))$$
to be the adjoint of the  push-down 
$$\pi^\Gamma_*:C^\infty(\partial X,V(\tilde\sigma_{-\lambda},\tilde\vp))\rightarrow B_{\Gamma}(\tilde\sigma_{-\lambda},\tilde\vp)\ .$$ 
\end{ddd}

The extension 
$ext^\Gamma$ is a holomorphic family of continuous maps and has values in $\Gamma$-invariant distributions.
We will also need the partial extensions
$ext^{\Gamma/U_P}:D_{\Gamma,k}(\sigma_\lambda,\vp)\rightarrow D_{U_P,k}(\sigma_\lambda,\vp)$,
which are defined as the adjoints of $\pi^{\Gamma/U_P}_*$.
Taking the adjoint of the relations (\ref{useid}) and (\ref{useid2})
we obtain
$$ext^{U_P}\circ ext^{\Gamma/U_P}=ext^\Gamma,\quad 
T_P^*\circ ext^{\Gamma/U_P}=\chi_p,\quad  \chi_P ext^{\Gamma/U_P}  = (T^P)^*$$
(where $P\in\tilde\cP^<$ for the last two, and multiplication by
$\chi_p$ implicitly involves projection onto the first component),
and
$$AS_P^*\circ ext^{\Gamma/U_P}=AS_P^*\ ,$$
if $P\in\tilde\cP^{max}$.

\subsubsection{}
The space $D_{\{1\},k}(\sigma_\lambda,\vp)$ defined
in Def.~\ref{neuj1002} depends on the choice of a parabolic subgroup $P$.
In the present subsection we write $D_{P,k}(\sigma_\lambda,\vp)$
for that space and let $$D_{\{1\},k}(\sigma_\lambda,\vp):=\bigcap_{p\in \cP_\Gamma,P\in p} D_{P,k}(\sigma_\lambda,\vp)\ .$$
Using the relation $ext^{U_P}\circ ext^{\Gamma/U_P}=ext^\Gamma$
we see that $$ext^\Gamma: D_{\Gamma,k}(\sigma_\lambda,\vp)
\rightarrow D_{\{1\},k}(\sigma_\lambda,\vp)\ ,$$ and that it is holomorphic
as a map $ext^\Gamma: D_{\Gamma,k}(\sigma_\lambda,\vp)
\rightarrow D_{\{1\},k_1}(\sigma_\lambda,\vp)$ for
any $k_1>k$ (and, of course, for $\Ree(\lambda)>\delta_\Gamma+\delta_\vp$).

\subsubsection{}\label{neuj3002}

From now on we consider the spherical case $\sigma=1$.
We introduce the restriction maps
$$res^\Gamma: D_{\{1\},k}(1_\lambda,\vp)\rightarrow D_{\Gamma,k_1}(1_\lambda,\vp)$$
defined for all $k$ and suitable $k_1\ge k$ depending on $k$.
The map $res^\Gamma$ is a refinement of the naive restriction
$$\{res^\Gamma\}: C^{-\infty}(\partial X,V(1_\lambda,\vp))\rightarrow
\cS_\Gamma(1_{-\lambda},\tilde\vp)^*$$ which is given
as the adjoint of $\{Q\}$ (see \ref{neuj1003}).
We define the meromorphic family of maps
$res^\Gamma$  
by
\begin{equation}\label{rreessdef}
res^\Gamma(f):=\pi^\Gamma_*(\chi^\Gamma\chi_0 f) +
\sum_{P\in\tilde\cP} T_P^* \circ res^{U_P}(f)\ ,\end{equation}
where we interprete $T_P^*$ as $AS_P^*$,
if $P\in\tilde\cP^{max}$. Here $\pi^\Gamma_*$ simply stands for the average over $\Gamma$.
It is well-defined and holomorphic for all $\lambda\in\aca$ since
$\chi^\Gamma\chi_0 f\in C_c^{-\infty}(\Omega_\Gamma,V(1_\lambda,\vp))$
and thus $\pi^\Gamma_*(\chi^\Gamma\chi_0 f)\in C_c^{-\infty}(B_\Gamma,V_{B_\Gamma}(1_\lambda,\vp))$.
Alternatively, one could write $\chi_0\{res^\Gamma\}(f)$ for  $\pi^\Gamma_*(\chi^\Gamma\chi_0 f)$.
Like the maps $res^{U_P}$ (see Definition \ref{neuj707}) the restriction map $res^\Gamma$ depends on choices.

\subsubsection{}

A priori the composition $res^\Gamma\circ ext^\Gamma$
is meromorphic and depends on many choices. The following lemma shows that the situation is much better. It generalizes \ref{neuj2001}.
\begin{lem}\label{weih114}
The composition 
$res^\Gamma\circ ext^\Gamma$
is regular and
coincides with the inclusion
of $$D_{\Gamma,k}(1_\lambda,\vp)\to D_{\Gamma,k_1}(1_\lambda,\vp)\ .$$
\end{lem}
\proof
Assume that $\lambda\in \aca$ is such that $res^\Gamma$ is regular.
Note that we still assume $\Ree(\lambda)>\delta_\Gamma+\delta_\vp$ in order ensure
convergence of $ext^\Gamma$.
Let $f\in D_{\Gamma,k}(1_\lambda,\vp)$.
It is easy to see that 
$\pi^\Gamma_*(\chi^\Gamma\chi_0 ext^\Gamma(f))=\chi_0 f$.
Furthermore, if $P\in\tilde\cP^<$, then
\begin{eqnarray*}
T_P^* \circ res^{U_P} \circ ext^\Gamma(f) &=&T_P^* \circ res^{U_P}
\circ ext^{U_P}\circ ext^{\Gamma/U_P}(f)\\
&=&T_P^* \circ ext^{\Gamma/U_P}(f)\\
&=&\chi_p f\ .
\end{eqnarray*}
If $P\in\tilde\cP^{max}$, then
we have
\begin{eqnarray*}
AS_P^* \circ res^{U_P} \circ ext^\Gamma(f) &=&AS_P^* \circ res^{U_P}
\circ ext^{U_P}\circ ext^{\Gamma/U_P}(f)\\
&=&AS_P^* \circ ext^{\Gamma/U_P}(f)\\
&=&AS_P (f)\ .
\end{eqnarray*}
Summing these equations over $\tilde\cP\cup\{0\}$, then we obtain
the desired identity $$res^\Gamma\circ ext^\Gamma(f)=f\ .$$ 
Since this equation holds true for generic $\lambda$, the assertion of the lemma follows. \hB

\subsubsection{}

We also have the following useful identity
$$T_P^*\circ res^{U_P}= \chi_p res^\Gamma$$
for all $P\in p \in \cP^<_\Gamma$.

\subsubsection{}\label{weih113}
Consider $P\in \tilde\cP$ and
recall the Definition \ref{neuj108} of $Ext_{U_P}(1_\lambda,\vp)$.
If $\lambda$ is admissible (for $P$) in the sense of \ref{neuj300}, then on this space the restriction $res^{U_P}$ is well-defined independent of choices. 

We now say that $\lambda\in \aca$ is admissible, if it is admissible for all $P\in \tilde\cP$.
Furthermore we define
$${}^\Gamma C^{-\infty}(\partial X,V(1_\lambda,\vp))^0:=
{}^\Gamma C^{-\infty}(\partial X,V(1_\lambda,\vp))\cap\bigcap_{P\in \tilde\cP}
Ext_{U_P}(1_\lambda,\vp)\ .$$
It immediately follows from (\ref{rreessdef}) that $res^\Gamma$ is well-defined on this space independent of choices. 
\subsubsection{}\label{neuj3001}

If we assume that $\lambda\not\in I_\aaaa$ or $\lambda>-\rho$, then 
$${}^\Gamma C^{-\infty}(\partial X,V(1_\lambda,\vp))^0={}^\Gamma C^{-\infty}(\partial X,V(1_\lambda,\vp))$$
by Prop.~\ref{maertins}, 2. Thus $res^\Gamma$ is pointwise well-defined at admissible $\lambda$
satisfying these additional conditions.

\subsubsection{}

\begin{ddd}
For $\Ree(\lambda)>\delta_\Gamma+\delta_\vp$
we define the scattering matrix 
$$S_\lambda^\Gamma:D_\Gamma(1_\lambda,\vp)\rightarrow D_\Gamma(1_{-\lambda},\vp)$$
as a meromorphic family of operators, which is given by 
$$S_\lambda^\Gamma:=res^\Gamma \circ J_\lambda\circ ext^\Gamma$$
provided $-\lambda$ is  non-integral and admissible in the sense of \ref{weih113}.
\end{ddd}
Using the scattering matrices for the cusps (Definition \ref{neuj2003}) we can rewrite the formula for $S^\Gamma_\lambda$ as follows:
\begin{equation}\label{spltt}
S_\lambda^\Gamma =\pi^\Gamma_* \circ\chi^\Gamma\chi_0 J_\lambda \circ ext^\Gamma + \sum_{P\in \tilde \cP} T^*_P\circ S^{U_P}_\lambda \circ ext^{\Gamma/U_P}\ .\end{equation}
Using this formula and the results of 
\ref{neuj2004}
we see that $S_\lambda^\Gamma$ is a well-defined
 meromorphic family of continuous maps.

\subsubsection{}\label{gaga}

As in the case of pure cusps (Lemma \ref{fundd}) we obtain a natural inclusion
$$B_\Gamma(1_\lambda,\vp)\hookrightarrow D_\Gamma(1_\lambda,\vp)$$
provided that 
$2\rho^{U_P}>l_{\vp_{|U_P}}$ for all $P\in  \tilde \cP^<$.
Note that $l_{\vp_{|U_P}}=2\delta_{\vp_{|U_P}}$.
%

\begin{lem}\label{restrf}
We assume that all cusps of $\Gamma$ have smaller rank.
Assume that $\delta_{\vp_{|U_P}}<\rho^{U_P}$ for all $P\in\tilde\cP$.
Then the restriction of the scattering matrix to $B_\Gamma(1_\lambda,\vp)$
induces  a meromorphic family of continuous maps
$$( S_\lambda^\Gamma)_{|B_\Gamma(1_\lambda,\vp)}:B_\Gamma(1_\lambda,\vp)\rightarrow B_\Gamma(1_{-\lambda},\vp)\ .$$
\end{lem}
\proof
First note that if $f\in B_\Gamma(1_\lambda,\vp)$ (considered as distribution),
then $ext^\Gamma(f)$ is smooth on $\Omega_\Gamma$.
It follows that $J_\lambda\circ ext^\Gamma(f)$ is smooth on $\Omega_\Gamma$,
and hence
$S^\Gamma_\lambda(f)=res^\Gamma\circ J_\lambda\circ ext^\Gamma(f)$
is represented by a smooth function.
We choose a cut-off function $\tilde\chi_P$ on $B_{U_P}$ which is
supported in $e_{P}(B_p)$ such that $\tilde\chi_P\chi_P=\chi_P$.
In order to see that $S^\Gamma_\lambda(f)\in
B_\Gamma(1_{-\lambda},\vp)$ we employ
Equation (\ref{spltt}). Let $P\in\tilde\cP$.
Then we have
\begin{eqnarray}
T_P\circ S^\Gamma_\lambda(f)&=&
 T_P(\chi^\Gamma \chi_0 J_\lambda\circ ext^\Gamma(f))\\
&&+T_P\circ T_P^*\circ S^{U_P}_\lambda\circ ext^{\Gamma/U_P}(f)\nonumber\\
&=&T_P(\chi^\Gamma \chi_0 J_\lambda\circ ext^\Gamma(f))\label{ppo0}\\
&&+\chi_P^2 S_\lambda^{U_P}\circ \tilde\chi_P ext^{\Gamma/U_P}(f)\label{ppo1}\\
&&+\chi_P^2  S_\lambda^{U_P}\circ(1- \tilde\chi_P) ext^{\Gamma/U_P}(f)\label{ppo2}\ .
\end{eqnarray}
The term (\ref{ppo0}) belongs to $C_c^\infty(B_\Gamma,V_{B_\Gamma}(1_\lambda,\vp))$.
We have $$\tilde\chi_P ext^{\Gamma/U_P}(f)=(\tilde T^P)^*(f)=\tilde T_P(f)\in
B_{U_P}(1_{-\lambda},\vp)\ ,$$ where $\tilde T_P, \tilde T^P$ are defined
as $T_P,T^P$, but using $\tilde\chi_P$ instead of $\chi_P$.
We can now apply Lemma \ref{inside} to (\ref{ppo1}) and Lemma \ref{offside} to (\ref{ppo2})
in order to see that these terms belong to  $B_\Gamma(1_{-\lambda},\vp)$.
\hB

\subsection{Vanishing results}

\subsubsection{}

If $\lambda\in\aca$ is admissible in the sense of \ref{weih113} and $\Ree(\lambda)>-\rho$,
then $res^\Gamma$ is 
well-defined on ${}^\Gamma C^{-\infty}(\partial X,V(1_\lambda,\vp))$
(see \ref{neuj3001}).

\subsubsection{}

Let $f\in {}^\Gamma C^{-\infty}(\partial X,V(1_\lambda,\vp))$.
The condition $\{res^\Gamma\}(f)=0$ (see \ref{neuj3002}) is equivalent to the condition that the support of $f$ as a distribution is contained in the limit set $\Lambda_\Gamma$.

In the case that $\Gamma$ is convex cocompact we know from \cite{MR1749869}, Thm.~4.7 that the space
$$\{f\in {}^\Gamma C^{-\infty}(\partial X,V(1_\lambda,\vp))|\{res^\Gamma\}(f)=0\}$$
is trivial if $\Ree(\lambda)>\delta_\Gamma+\delta_\vp$.
Already the example of pure cusps shows that this is not true in the presence of cusps.

Note that in the convex-cocompact case $\{res^\Gamma\}=res^\Gamma$.
It is at the heart of the matter that for geometrically finite groups $\Gamma$ and
$\Ree(\lambda)$ large the stronger condition
$res^\Gamma(f)=0$ implies $f=0$.
In \cite{MR1926489} we introduced a related condition ``$f$ is strongly supported on the limit set'' in order to show 
such vanishing results.

\subsubsection{}

\begin{prop}\label{vani}
We assume  that 
$\lambda$ is admissible in the sense of \ref{weih113}
and satisfies 
\begin{equation}\label{neuj3004}
\Ree(\lambda)>\max\left(\{\delta_\Gamma\}\cup\{\rho_{U_P}\mid P\in\tilde\cP\}\right)+\delta_\vp
\ .
\end{equation} 
Let $f\in {}^\Gamma C^{-\infty}(\partial X,V(1_\lambda,\vp))$. If
$res^\Gamma(f)=0$, then  $f=0$.
\end{prop}
\proof

\subsubsection{}

Let $\lambda$ be admissible and $f\in {}^\Gamma C^{-\infty}(\partial X,V(1_\lambda,\vp))$. From the defining formula (\ref{rreessdef}) for $res^\Gamma$
we conclude that $res^\Gamma(f)=0$ if and only if $f$ is supported on $\Lambda_\Gamma$
and $res^{U_p}(f)\in  C_c^{-\infty}(B_{U_P},V_{B_{U_P}}(1_\lambda,\vp))\subset D_{U_P}(1_\lambda,\vp)$ for all $P\in\tilde\cP$. The latter condition implies
that 
\begin{equation}\label{wuff} 
res^{U_P}(f)=\{res^{U_P}\}(f)
\end{equation} 
for all $P\in\tilde\cP$.

\subsubsection{}

If $\Ree(\lambda)>-\rho^{U_P}$ and $h\in C^{\infty}(\partial X,V(1_{-\lambda},\tilde\vp))$, then Proposition \ref{maertins}, 2, Lemma \ref{wauwau}, and Equation (\ref{wuff}) imply 
$$
\langle f,h\rangle =\langle ext^{U_P}\circ res^{U_P}(f),h\rangle 
=\langle \{res\}^{U_P}(f),\pi_*^{U_P}(h)\rangle $$
and therefore
\begin{equation}\label{woff}
\langle f,h\rangle =\sum_{u\in U_P} \langle\chi^{U_P}f,\tilde\vp(u)h(u^{-1}.)\rangle\ .
\end{equation}
Thus, if $res^\Gamma(f)=0$, then $f$ is supported on the limit set and satisfies (\ref{woff})
for all $P\in\tilde\cP$, $h\in C^{\infty}(\partial X,V(1_{-\lambda},\tilde\vp))$. These are
precisely the defining conditions for being ``strongly supported on the limit set'' in the sense of \cite{MR1926489}. The main result of \cite{MR1926489} now states that an element
$f\in {}^\Gamma C^{-\infty}(\partial X,V(1_\lambda,\vp))$ that is ``strongly supported on the limit set'' vanishes provided that (\ref{neuj3004}) holds. The proposition follows.
\hB

\subsubsection{}
 
\begin{kor}\label{wieder}
Assume $\lambda\in\aca$ satisfies (\ref{neuj3004}). If
$f_\mu\in {}^\Gamma C^{-\infty}(\partial X,V(1_\mu,\vp))$
is a germ of a meromorphic family at $\lambda$, then
$ext^\Gamma\circ res^\Gamma(f_\mu) = f_\mu$.
\end{kor}
\proof
We apply $res^\Gamma$, and we use $res^\Gamma\circ ext^\Gamma=\id$ and the injectivity
of $res^\Gamma$ for generic $\mu$ near $\lambda$ proved in Proposition \ref{vani}.
\hB.

\subsection{Meromorphic continuations and general $\sigma$}

\subsubsection{}

Recall that we have defined $ext^\Gamma$
for $\Ree(\lambda)>\delta_\Gamma+\delta_\vp$.
The scattering matrix  $S^\Gamma_\lambda$ was defined
on the same range of $\lambda$
for the trivial $M$-type
$1$.

In the present subsection we extend the definition
of the scattering matrix to general $M$-types.
Further we show that the extension
and the scattering matrix have meromorphic continuations
to all of $\aca$.
This finishes the proof of theorems \ref{t119}, \ref{t110}, and \ref{t113}.

\subsubsection{}

We start with the meromorphic continuation of $S^\Gamma_\lambda$ to a certain half-plane
$W\subset \aca$ in the spherical case. Under these conditions
we can construct a meromorphic family of parametrices
for $S_\lambda^\Gamma$ with finite-dimensional singularities.
In a second step we employ twisting in order to show
meromorphy of the extension and the scattering matrix on all of $\aca$
and for general $\sigma$.

\subsubsection{}

\begin{prop}\label{meroext}
The extension
$$ext^\Gamma:D_\Gamma(1_\lambda,\vp)\rightarrow C^{-\infty}(\partial X,V(1_\lambda,\vp))$$
and the scattering matrix $$S_\lambda^\Gamma:D_\Gamma(1_\lambda,\vp)\rightarrow
D_\Gamma(1_{-\lambda},\vp)$$ have meromorphic continuations to $$W:=\{\lambda\in\aca\:|\:
\Ree(\lambda)>-\rho+\beta\}\ ,$$ where $\beta=0$ for $X=\RH^n, \CH^n$,
and $\beta=2\alpha$ for $X=\HH^n$.
The family $ext^\Gamma$ has finite-dimensional singularities.
\end{prop}
\proof
\subsubsection{}
We first show the proposition under an additional assumption on $\Gamma$ and $\vp$.
\begin{lem}\label{wahn1}
We assume that
\begin{equation}\label{kuckuck}
\delta_\Gamma+\delta_\vp<-\max\left(\{0\}\cup\{\rho_{U_P}+\delta_\vp\mid P\in\tilde\cP\}\right)
\ .
\end{equation} 
Then the extension $ext^\Gamma$
and the scattering matrix $S_\lambda^\Gamma$ have
a meromorphic continuation to the half-plane $W$.
The family $ext^\Gamma$ has finite-dimensional singularities.
\end{lem}
\proof
\subsubsection{}
The additional assumption in particular implies that $\rho^{U_P}>\delta_\vp\ge 0$
for all $P\in\tilde\cP$.
Therefore
all cusps of $\Gamma$ have smaller rank, and 
there is an embedding
$$B_\Gamma(1_\lambda,\vp)\hookrightarrow D_\Gamma(1_\lambda,\vp)$$
(see \ref{gaga}).
We consider the non-empty open subset
$$U:=\{\lambda\in\aca\:|\:\max\{\delta_\Gamma+\delta_\vp,-\rho+\beta\}<\Ree(\lambda)<\rho-\beta\}\subset \aca\ .$$
Based on Lemma \ref{restrf} we first consider the restriction of the scattering matrix
to $B_\Gamma(1_\lambda,\vp)$.
For $\lambda\in U$ we construct a parametrix (see \ref{neuj3005} for $T_P$ and \ref{neuj3006} for $\tilde T^P$)
$$Q_{-\lambda}:B_\Gamma(1_{-\lambda},\vp)\rightarrow B_\Gamma(1_\lambda,\vp)$$
of $S^\Gamma_\lambda$ by
$$Q_{-\lambda} (f) := \pi^\Gamma_*(\chi^\Gamma \chi_0 J_{-\lambda} (\chi^\Gamma
\tilde\chi_0 f)) + \sum_{P\in\tilde\cP} T_P^* \circ S_{-\lambda}^{U_P}\circ (\tilde T^P)^* (f)\ .$$
Here $\tilde \chi_0$ is some cut-off function on $B_\Gamma$
of compact support such that $\chi_0\tilde \chi_0=\chi_0$, and
$\tilde T^P$ is defined in the same way as
$T^P$ but using a cut-off function $\tilde\chi_p$ with support
on $B_p$ such that $\chi_p\tilde\chi_p=\chi_p$.

Note that $Q_{-\lambda}$ has a continuous extension 
to a map
$$Q_{-\lambda}:D_\Gamma(1_{-\lambda},\vp)\rightarrow D_\Gamma(1_\lambda,\vp)\ .$$
This is clear for the second term $\sum_{P\in\tilde\cP} T_P^* \circ S_{-\lambda}^{U_P}\circ (\tilde T^P)^*$
since the scattering matrices $S_{-\lambda}^{U_P}$ extend to distributions.
The map $f\mapsto\chi^\Gamma\chi_0J_{-\lambda}(\chi^\Gamma\tilde\chi_0f)$
extends to a continuous map from $D_\Gamma(1_{-\lambda},\vp)$ to
$C^{-\infty}_c(\Omega_\Gamma,V(1_\lambda,\vp))$.
The push-down extends to a continuous map from
$C^{-\infty}_c(\Omega_\Gamma,V(1_\lambda,\vp))$
to $C^{-\infty}_c(B_\Gamma,V_{B_\Gamma}(1_\lambda,\vp))$.
This implies that the first term
$\pi^\Gamma_*(\chi^\Gamma \chi_0 J_{-\lambda} (\chi^\Gamma \tilde\chi_0 f))$
extends to distributions, too.

\subsubsection{}
We claim that for $|\Ree(\lambda)|<\rho-\beta$ the scattering matrix $S_\lambda^{U_P}$
has at most finite-dimensional singularities. To see the claim  we write
$S_\lambda^{U_P}=res^{U_P}\circ J_\lambda\circ ext^{U_P}$.
Since $J_\lambda$ is regular and bijective for these
$\lambda$, and $ext^{U_P}$ has finite-dimensional singularities,
the only term which may contribute infinite-dimensional singularities
is $res^{U_P}$. 

Since $\{res^{U_P}\}$ is always regular the singular part of $res^{U_P}$ has values
in the space $\cR_{U_P}(1_{\lambda},\tilde\vp)^*$. 
Now $S_\lambda^{U_P}$ is the continuous extension of its restriction
to the dense subspace $B_{U_P}(1_\lambda,\vp)$. 
Since  $S^{U_P}_\lambda$ maps $B_{U_P}(1_\lambda,\vp)$ to $B_{U_P}(1_{-\lambda},\vp)$
the range of the singular part of $S_\lambda^{U_P}$ does not
contain non-trivial elements of $\cR_{U_P}(1_{\lambda},\tilde\vp)^*$.

Let $D_{U_P}(1_\lambda,\vp)^0\subset D_{U_P}(1_\lambda,\vp)$
be a closed subspace of finite codimension on which $ext^{U_P}$
is regular and put $B_{U_P}(1_\lambda,\vp)^0:=
D_{U_P}(1_\lambda,\vp)^0\cap B_{U_P}(1_\lambda,\vp)$.
Then $B_{U_P}(1_\lambda,\vp)^0$ has finite codimension
in $B_{U_P}(1_\lambda,\vp)$, and the closure of
$B_{U_P}(1_\lambda,\vp)^0$ in $D_{U_P}(1_\lambda,\vp)$ is
$D_{U_P}(1_\lambda,\vp)^0$. 
Looking at Laurent expansions we conclude that the restriction of $S_\lambda^{U_P}$ to $B_{U_P}(1_\lambda,\vp)^0$ and hence to
$D_{U_P}(1_\lambda,\vp)^0$ is regular.
This proves the claim.

We conclude that for
$\pm\lambda\in U$ the family $Q_\lambda$ has at most finite-dimensional
singularities.

\subsubsection{}
We define the remainder
$$R_\lambda:=Q_{-\lambda}\circ S^\Gamma_\lambda-\id \ .$$
We are going to show that $R_\lambda$ is a meromorphic family of
smoothing operators with finite-dimensional singularities.
We start with
\begin{eqnarray*}
R_\lambda(f)&=&Q_{-\lambda}\circ S^\Gamma_\lambda(f)-f\\
&=&\pi^\Gamma_* (\chi^\Gamma \chi_0 J_{-\lambda}( \chi^\Gamma\tilde\chi_0 S^\Gamma_\lambda(f)))\\
&& + \sum_{P\in\tilde\cP} T_P^* S_{-\lambda}^{U_P}(\tilde T^P)^* S^\Gamma_\lambda(f)-f
\end{eqnarray*}
Next we insert the definition $S^\Gamma_\lambda=res^\Gamma\circ J_\lambda\circ
ext^\Gamma$ and the definition of $res^\Gamma$  (\ref{rreessdef}). We obtain 
\begin{eqnarray*}
R_\lambda(f)&=&
\pi^\Gamma_*(\chi^\Gamma \chi_0 J_{-\lambda} (\chi^\Gamma\tilde\chi_0 \pi^\Gamma_*(\chi^\Gamma\chi_0 J_\lambda\circ ext^{\Gamma}(f)))) \\
&&+\sum_{P\in\tilde\cP}\pi^\Gamma_* (\chi^\Gamma \chi_0 J_{-\lambda} (\chi^\Gamma\tilde\chi_0  T_P^*  res^{U_P}\circ J_\lambda \circ ext^{\Gamma}(f)))\\
&&+\sum_{P\in\tilde\cP} T_P^*  S_{-\lambda}^{U_P}
(\tilde T^P)^*\pi^\Gamma_*( \chi^\Gamma\chi_0 J_\lambda  \circ ext^\Gamma(f))\\
&&+\sum_{P,Q\in\tilde\cP} T_P^*  S_{-\lambda}^{U_P} (\tilde T^P)^*
T_Q^*  (res^{U_Q}\circ J_\lambda \circ ext^{\Gamma} (f)) - f\ .
\end{eqnarray*}
Now we employ 
\begin{eqnarray*}
\chi^\Gamma\tilde\chi_0 \pi^\Gamma_* \chi^\Gamma \chi_0&=&
\chi^\Gamma\chi_0\\
\chi^\Gamma\tilde\chi_0\sum_{P\in\tilde\cP} T_P^*\circ  res^{U_P}
&=& \tilde\chi_0 (1-\chi_0)\chi^\Gamma\\
\sum_{Q\in\tilde\cP}  (\tilde T^P)^*
T_Q^* res^{U_Q} + (\tilde T^P)^*\pi^\Gamma_* \chi^\Gamma\chi_0 & =&  \tilde \chi_P res^{U_P}\\
res^{U_P}\circ J_\lambda\circ ext^\Gamma &=& S^{U_P}_\lambda\circ ext^{\Gamma/U_P}
\end{eqnarray*} in order
to obtain 
\begin{eqnarray*}
R_\lambda(f)
&=&\pi^\Gamma_*(\chi^\Gamma \chi_0 J_{-\lambda} (\chi^\Gamma\tilde\chi_0 J_\lambda \circ ext^{\Gamma}(f))) \\
&&+\sum_{P\in\tilde\cP} T_P^*  S_{-\lambda}^{U_P} (\tilde\chi_P
S_{\lambda}^{U_P} \circ ext^{\Gamma/U_P}(f))- f\ .
\end{eqnarray*}
Using the functional equations of the intertwining operators and scattering matrices we can further write
\begin{eqnarray*}
R_\lambda(f)
&=&\pi^\Gamma_*(\chi^\Gamma \chi_0 ext^{\Gamma}(f)) \\
&&-\pi^\Gamma_*(\chi^\Gamma \chi_0 J_{-\lambda}( (1-\chi^\Gamma\tilde\chi_0) J_\lambda \circ ext^{\Gamma}(f))) \\
&&+\sum_{P\in\tilde\cP} T_P^*
 ext^{\Gamma/U_P}(f)\\
&&-\sum_{P\in\tilde\cP} T_P^*  S_{-\lambda}^{U_P} ((1-\tilde\chi_P)
S_{\lambda}^{U_P} \circ ext^{\Gamma/U_P}(f))- f\\
&=&-\pi^\Gamma_*(\chi^\Gamma \chi_0 J_{-\lambda}( (1-\chi^\Gamma\tilde\chi_0) J_\lambda\circ ext^{\Gamma}(f))) \\
&&-\sum_{P\in\tilde\cP} T_P^*  S_{-\lambda}^{U_P} ((1-\tilde\chi_P)
S_{\lambda}^{U_P}\circ  ext^{\Gamma/U_P}(f))\ .
\end{eqnarray*}
Note that $\chi^\Gamma \chi_0  (1-\chi^\Gamma\tilde\chi_0)=0$,
and that $J_{-\lambda}$ is off-diagonal smoothing. Furthermore note
that $\chi_P (1-\tilde\chi_P)=0$, and that $S_{-\lambda}^{U_P}$
is off-diagonal smoothing by Lemma \ref{offside}.
We conclude that
$R_\lambda$ extends to a meromorphic family of continuous maps
$$R_\lambda:D_\Gamma(1_\lambda,\vp)\rightarrow B_\Gamma(1_\lambda,\vp)$$
with finite-dimensional singularities.
Since $J_0=\id$ and $S_0^{U_P}=\id$ we have $R_0=0$.

\subsubsection{}

We define the  Banach space $\cS_{\Gamma,0,0}(1_\lambda,\vp)$ to be the subspace
of all $f\in C(B_\Gamma,V_{B_\Gamma}(1_\lambda,\vp))$
such that $T_P(f)\in \cS_{U_P,0,0}(1_\lambda,\vp)$ for all $P\in \tilde\cP$ (see \ref{weih259}).
Then we have compact inclusions
$$ B_\Gamma(1_\lambda,\vp)\hookrightarrow \cS_{\Gamma,0,0}(1_\lambda,\vp)
\hookrightarrow D_\Gamma(1_\lambda,\vp)$$
(in order to see the second inclusion use Lemma \ref{fundd}).
The family of Banach spaces $\{\cS_{\Gamma,0,0}(1_\lambda,\vp)\}_{\lambda\in\aca}$
forms a trivial holomorphic bundle.

\subsubsection{}
We now apply meromorphic Fredholm theory \cite{reedsimon78}, Thm. VI.13,
(or better its version for families of operators on Banach spaces)
to the family $$1+R_\lambda:\cS_{\Gamma,0,0}(1_\lambda,\vp)\rightarrow
\cS_{\Gamma,0,0}(1_\lambda,\vp)\ .$$
As in  \cite{MR1749869}, Lemma 3.6,  we conclude that $(1+R_\lambda)^{-1}$ exists as a meromorphic
family of maps on $\cS_{\Gamma,0,0}(1_\lambda,\vp)$ of the form
$(1+T_\lambda)$, where $T_\lambda$ is a meromorphic family of continuous
operators from $D_\Gamma(1_\lambda,\vp)$ to $B_\Gamma(1_\lambda,\vp)$
with finite-dimensional singularities.

\subsubsection{}
We see that the inverse of the scattering matrix is given on $U$ by $$(S^\Gamma_\lambda)^{-1} := (1+R_\lambda)^{-1} Q_{-\lambda}\ .$$
as a meromorphic family and with finite-dimensional singularities.

\subsubsection{}
We have $-U\cup \{\lambda\in\aca\:|\: \Ree(\lambda)>\delta_\Gamma+\delta_\vp\}=W$.
By our assumption the set  
$$W_0=\left\{\lambda\in\aca\:|\: \delta_\Gamma+\delta_\vp <
\Ree(\lambda)< -\delta_\vp-\max\left(\{\delta_\Gamma\}\cup\{\rho_{U_P}\mid P\in\tilde\cP\}\right)\right\}$$
is non-empty.
By Corollary \ref{wieder} we find for $\lambda\in W_0$ the meromorphic identity
\begin{eqnarray*}
S_{-\lambda}^\Gamma\circ S^\Gamma_\lambda&=&
res^\Gamma\circ J_{-\lambda}\circ ext^\Gamma\circ res^\Gamma\circ J_\lambda\circ ext^\Gamma\\
&=&res^\Gamma\circ J_{-\lambda}\circ J_\lambda\circ ext^\Gamma\\
&=&res^\Gamma\circ ext^\Gamma\\
&=&\id  \ .
\end{eqnarray*}
Thus defining $S^\Gamma_\lambda:=(S_{-\lambda}^\Gamma)^{-1}$ for $-\lambda\in U$
we obtain the continuation of $S^\Gamma_\lambda$ to $W$.

\subsubsection{}

We obtain the meromorphic continuation of $ext^\Gamma$ to $W$
with finite-dimensional singularities  using the identity
$$ext^\Gamma=J_{-\lambda}\circ ext^\Gamma\circ S_\lambda^\Gamma\ .$$
See the proof of \cite{MR1749869}, Lemma 5.11,
for the corresponding argument.
\hB

\subsubsection{}

The following lemma finishes the proof of Proposition \ref{meroext}. 

\begin{lem}\label{hjdhdjwqdwqdwqdwq}
Lemma \ref{wahn1} holds true without the assumption (\ref{kuckuck}).
\end{lem}
\proof
We are going to employ the embedding trick.
First we assume  that $G^n$ belongs to the list
$\{Spin(1,n),SO(1,n)_0,SU(1,n),Sp(1,n)\}$. 
If we consider $\Gamma\subset G^n$ as a subgroup of $G^{n+m}$,
then for sufficiently large $m$ the inequality (\ref{kuckuck})
is satisfied. Indeed, all ingredients of (\ref{kuckuck}) but $\delta_\Gamma$
are independent of $m$, while $\delta_\Gamma^{n+m}=\delta^n_\Gamma-m\zeta$, $\zeta=\rho^{n+1}-\rho^n$.

It thus suffices to show that a meromorphic continuation of $ext^{\Gamma,n+m}$
implies a meromorphic continuation of $ext^{\Gamma,n}$. 
Assume that 
we have a meromorphic continuation of $ext^{\Gamma,n+m}$
to $W^{n+m}:=\{\lambda\in\aca\:|\:
\Ree(\lambda)>-\rho^{n+m}+\beta\}$ with finite-dimensional singularities.
Then we employ the diagram
$$\begin{array}{ccc}  C^{-\infty}(\partial X,V(1^{n}_\lambda,\vp)) & 
 \stackrel{ j^*}{ \leftarrow} &C^\infty(\partial X^{n+m}, V(1^{n+m}_{\lambda-m\zeta},\vp))
 \\
 \uparrow ext^{\Gamma,n}_*&&\uparrow ext^{\Gamma,n+m}_* \\ D_\Gamma(1^{n}_\lambda,\vp)&\stackrel{i^\Gamma_*}{\rightarrow }&D_\Gamma(1^{n+m}_{\lambda-m\zeta},\vp)
  \end{array}$$
which apriori holds in the domain of convergence 
$\Ree(\lambda)>\delta_\Gamma^n+\delta_\vp^n$ by Lemma \ref{comppp}. Note that the left-inverse $j^*$ of $i^*$ was constructed in Lemma \ref{jjjjjjjjjjjtwzew} (or better its iterate for the $m$-fold embedding). 
We can define the meromorphic continuation to
$W^n=W^{n+m}+m\zeta$
of $ext^{\Gamma,n}$ using the identity $ext^{\Gamma,n}=j^*\circ ext^{\Gamma,n+m}\circ i^\Gamma_*$.

Let now $G$ be general. Then we can find a subgroup
$\Gamma^0\subset \Gamma$ of finite index to which we can apply
the concept of embedding (see \ref{ohu}). In particular we obtain a
meromorphic continuation of
$ext^{\Gamma^0}$ to $W$. If we identify
$D_\Gamma(1_\lambda,\vp)$ with the subspace
${}^{\Gamma/\Gamma^0} D_{\Gamma^0}(1_\lambda,\vp)$,
then we obtain the meromorphic continuation
of $ext^\Gamma$ using the identity
$ext^{\Gamma}=ext^{\Gamma^0}_{|D_\Gamma(1_\lambda,\vp)}$.

Using the meromorphic continuation of $ext^{\Gamma,n}$
we then obtain the meromorphic continuation of
the scattering matrix $S^\Gamma_\lambda$ by
$S_\lambda^\Gamma=res^\Gamma\circ J_\lambda\circ ext^\Gamma$.
\hB
 
\subsubsection{}

In Proposition \ref{meroext} we have obtained (in the spherical case) a meromorphic extension of
$ext^\Gamma$ and the scattering matrix $S_\lambda^\Gamma$
to a half space $W\subset \aca$. In the following argument, using twisting,  we find the meromorphic continuation of these objects to all of $\aca$.

\subsubsection{}
 
Let $\sigma$ be a Weyl invariant representation of $M$,
and $(\pi_{\sigma,\mu},V_{\pi_{\sigma,\mu}})$ be a finite-dimensional
representation of $G$ as in
\ref{weih144}.
The adjoint of the diagram
2. of Lemma \ref{compat5}
gives
the following commutative diagram of
meromorphic families of maps
\begin{equation}\label{dieletzte}\begin{array}{ccc}
D_\Gamma(\sigma_\lambda,\vp)
&\stackrel{i^\Gamma_{\sigma,\mu}}{\rightarrow}&D_\Gamma(1_{\lambda+\mu},\pi_{\sigma,\mu}\otimes\vp)\\
\downarrow ext^\Gamma&&\downarrow ext^\Gamma\\
C^{-\infty}(\partial X,V(\sigma_\lambda,\vp))
&\stackrel{j_{\sigma,\mu}}{\leftarrow}&
 C^{-\infty}(\partial X,V(1_{\lambda+\mu},\pi_{\sigma,\mu}
\otimes\vp))
\end{array}
\end{equation}
for $\Ree(\lambda)>\delta_\Gamma+\delta_\vp$.
Using Proposition \ref{meroext} we conclude that
$$ext^\Gamma:D_\Gamma(\sigma_\lambda,\vp)\rightarrow C^{-\infty}(\partial X,V(\sigma_\lambda,\vp))$$
is meromorphic on the set $W-\mu$.
Since we can choose $\mu$ arbitrary large
we obtain the meromorphic continutaion of $ext^\Gamma$ to all of $\aca$.

From (\ref{dieletzte}),
injectivity of $i_{\sigma,\mu}$,
and Proposition \ref{meroext} we conclude that $ext^\Gamma$
has finite-dimensional singularities. This finishes the proof of Theorem \ref{t110}.
The adjoint of the extension is the push-down.
Therefore Theorem \ref{t119} follows from \ref{t110}.

\subsubsection{}

Let for a moment be $\sigma=1$. Using the meromorphic continuation of
the extension we obtain the meromorphic continuation of the
scattering matrix to all of $\aca$
by the formula $S^\Gamma_\lambda=res^\Gamma\circ J_\lambda\circ ext^\Gamma$.
In order to proceed similarly for general $\sigma$ we first have to define the 
restriction map in this situation.

Let now $\sigma$ be a Weyl-invariant representation of
$M$
and $(\pi_{\sigma,\mu},V_{\pi_{\sigma,\mu}})$ be a finite-dimensional
representation of $G$ as in
\ref{weih144}.
For $k\in \nat$
we define restriction maps
$$res^\Gamma:{}^\Gamma C^{-k}(\partial X,V(\sigma_\lambda,\vp))
\rightarrow D_\Gamma(\sigma_\lambda,\vp)$$
using the diagrams
$$\begin{array}{ccc}
C^{-k}(\partial X,V(\sigma_\lambda,\vp))
&\stackrel{i_{\sigma,\mu}}{\rightarrow}&C^{-k}(\partial X,V(1_{\lambda+\mu},\pi_{\sigma,\mu}\otimes\vp))\\
\downarrow res^\Gamma&&\downarrow res^\Gamma\\
D_\Gamma(\sigma_\lambda,\vp)&\stackrel{j^\Gamma_{\sigma,\mu}}{\leftarrow}&
D_\Gamma(1_{\lambda+\mu},\pi_{\sigma,\mu}\otimes\vp)
\end{array}\ .$$
The right column depends on choices (including $\mu$), and so does
the left column. But it is easy to check that 
\begin{equation}\label{juhu}
res^\Gamma\circ ext^\Gamma =\id
\end{equation}

\subsubsection{}\label{hjshqwjsdwhqdwqdwqdwqdqw}

Note that (\ref{dieletzte}) implies the meromorphic identity $ext^\Gamma\circ j_{\sigma,\mu}^\Gamma=j_{\sigma,\mu}\circ ext^\Gamma$.
Let $f_\nu\in  {}^\Gamma C^{-\infty}(\partial X,V(\sigma_\nu,\vp))$
be a meromorphic family defined on {\em all} of $\aca$.
We compute for $\nu\gg 0$ using Corollary \ref{wieder}
and any set of choices for $res^\Gamma$
\begin{eqnarray}\label{jetztaber}
ext^\Gamma\circ res^\Gamma(f_\nu)&=&
ext^\Gamma\circ j_{\sigma,\mu}^\Gamma \circ res^\Gamma\circ i_{\sigma,\mu}(f_\nu)\nonumber\\
&=&j_{\sigma,\mu}\circ ext^\Gamma\circ res^\Gamma\circ i_{\sigma,\mu}(f_\nu)\nonumber\\
&=&j_{\sigma,\mu}\circ i_{\sigma,\mu}(f_\nu)\nonumber\\
&=&f_\nu\ .
\end{eqnarray}
Formula (\ref{juhu}) implies that $ext^\Gamma$ is injective
on families.
Thus we can conclude from (\ref{jetztaber}) that $res^\Gamma(f_\nu)$ is independent of choices.

\subsubsection{}

Now we can  define the scattering matrix for general
$\sigma$ by the same formula as in the spherical case
$$S_\lambda^\Gamma:=res^\Gamma\circ J_\lambda\circ ext^\Gamma\ .$$
It is meromorphic on $\aca$.

In order to finish the proof of Theorem \ref{t113}
it remains to show the functional equation.
Using (\ref{jetztaber}) we find
\begin{eqnarray*}
S^\Gamma_{-\lambda}\circ S_\lambda^\Gamma&=&
res^\Gamma\circ J_{-\lambda} \circ ext^\Gamma\circ
res^\Gamma\circ J_\lambda\circ ext^\Gamma\\
&=&res^\Gamma\circ J_{-\lambda}\circ J_\lambda\circ ext^\Gamma\\
&=&res^\Gamma\circ ext^\Gamma\\
&=&\id\ .
\end{eqnarray*}

\bibliographystyle{plain}


\end{document}